\documentclass[twoside, a4paper,reqno]{amsart}
\usepackage[T1]{fontenc}
\usepackage[utf8]{inputenc}
\usepackage{lmodern}
\usepackage{graphicx}
\usepackage[all]{xy}
\usepackage{relsize}
\usepackage[british]{babel}
\usepackage[babel, final]{microtype}
\usepackage{mathtools,booktabs,multicol}
\usepackage[pdftex, dvipsnames]{xcolor}
\usepackage{caption}
\usepackage{subcaption}
\usepackage{multirow}
\usepackage{amssymb}
\usepackage[shortlabels]{enumitem}

\usepackage[margin=30mm]{geometry}
\usepackage[breaklinks,%
    pdftex,%
    final,%
    colorlinks=true,%
    linkcolor=NavyBlue,%
    citecolor=blue,
    filecolor=NavyBlue,%
    menucolor=NavyBlue,%
    urlcolor=NavyBlue,%
    bookmarks=true,%
    bookmarksdepth=2,%
    bookmarksnumbered=true,%
    bookmarksopen=true,%
    bookmarksopenlevel=2,%
    unicode
]{hyperref}
\usepackage{footnotebackref}

\usepackage{overpic}
\usepackage{tikz,tikz-cd}
\usetikzlibrary{tikzmark}
\tikzset{
    math to/.tip={Glyph[glyph math command=rightarrow]},
    loop/.tip={Glyph[glyph math command=looparrowleft, swap]},
    weird/.tip={Glyph[glyph math command=Rrightarrow, glyph length=1.5ex]},
    pi/.tip={Glyph[glyph math command=pi, glyph length=1.5ex, glyph axis=0pt]},
}

\newtheorem*{theoremA}{Theorem A}
\newtheorem*{theoremB}{Theorem B}

\newcommand{\wt}{\widetilde}
 \newcommand{\ol}{\overline}
\newcommand{\Z}{\mathbb{Z}}
 \newcommand{\homeo}{\approx}
 \newcommand{\km}{\mathrm{km}}
 
 \newcommand{\R}{\mathbb{R}}
\DeclareMathOperator{\Int}{Int}
 \newcommand{\hooklongrightarrow}{\lhook\joinrel\longrightarrow}
 \newcommand{\tw}{\mathrm{tw}}
 \newcommand{\spin}{\mathrm{spin}}
 \newcommand{\Arf}{\operatorname{Arf}}
 \newcommand{\calS}{\mathcal{S}}
 \newcommand{\calN}{\mathcal{N}}
 \DeclareMathOperator{\cone}{cone}
 \DeclareRobustCommand\sm{\mathbin{\mathpalette\smaux\relax}}
\newcommand\smaux[2]{\mspace{-4mu}
\raisebox{\rsmraise{#1}\depth}{\rotatebox[origin=c]{-25}{$#1\smallsetminus$}}
 \mspace{-4mu}
}
\newcommand\rsmraise[1]{%
  \ifx#1\displaystyle .5\else
    \ifx#1\textstyle .5\else
      \ifx#1\scriptstyle .3\else
        .45%
      \fi
    \fi
  \fi}

\newtheorem{theorem}{Theorem}[section]

\newtheorem{corollary}[theorem]{Corollary}

\newtheorem{proposition}[theorem]{Proposition}

\newtheorem{step}{Step}

\newtheorem*{rep@theorem}{\rep@title}
\newcommand{\newreptheorem}[2]{
\newenvironment{rep#1}[1]{
\def\rep@title{#2 \ref{##1}}
\begin{rep@theorem}}
{\end{rep@theorem}}}
\newreptheorem{theorem}{Theorem}
\newreptheorem{lemma}{Lemma}
\newreptheorem{proposition}{Proposition}
\newreptheorem{corollary}{Corollary}

\theoremstyle{definition}

\newtheorem{definition}[theorem]{Definition}

\newtheorem{exercise}{Exercise}

\theoremstyle{remark}
\newtheorem{remark}[theorem]{Remark}

\numberwithin{equation}{section}
\numberwithin{exercise}{subsection}

\title{Constructing locally flat surfaces in 4-manifolds}
\author{Arunima Ray}
\address{The University of Melbourne, Grattan St, Parkville VIC 3010, Australia.}
\email{aru.ray@unimelb.edu.au}

\usepackage[capitalise, noabbrev]{cleveref}
\Crefname{step}{step}{steps}
\Crefname{step}{Step}{Steps}
\Crefname{exercise}{exercise}{exercises}
\Crefname{exercise}{Exercise}{Exercises}

\begin{document}

\begin{abstract}
    There are two main approaches to building locally flat embedded surfaces in $4$-manifolds: direct methods which geometrically manipulate a given map of a surface, and more  indirect methods using surgery theory. Both rely on Freedman--Quinn's disc embedding theorem. In this expository article, we give an introduction to these methods by sketching proofs of the following results: every primitive second homology class in a closed, simply connected $4$-manifold is represented by a locally flat embedded torus (Lee--Wilczy\'{n}ski~\cite{LW97}); and every Alexander polynomial one knot in~$S^3$ is topologically slice (Freedman--Quinn~\cite{FQ}). 
\end{abstract}

\maketitle

\section{Introduction}
\label{sec:intro}
Embedded surfaces in $4$-manifolds have been instrumental in the study of $4$-manifolds in general. For example, Freedman's \emph{disc embedding theorem}~\cite{F} (\Cref{thm:DET,thm:DET-full} below) was the key ingredient in his celebrated proof of the topological $4$-dimensional Poincar\'e conjecture. Embedded surfaces are also used in numerous constructions of novel $4$-manifolds. For instance, \emph{knot surgery} and \emph{torus surgery} on embedded tori can be used to produce infinitely many distinct smooth structures on a given $4$-manifold~\cite{donaldson:h-cobordism,okonek-vandeven:infinite-log-transform,friedman-morgan:algebraicI,friedman-morgan:algebraicII,fintushel-stern:knot-surgery} and \emph{Gluck twists} on embeddings~$S^2\hookrightarrow S^4$ provide a possible means of constructing exotic smooth structures on $S^4$, which would disprove the smooth $4$-dimensional Poincar\'e conjecture. The minimal genus of an embedded surface representing an element of second homology, encoded in the \emph{genus function}, is a powerful invariant for $4$-manifolds. The study of embedded surfaces in $4$-manifolds is also a natural analogue of classical knot theory. Therefore it is not surprising that there is a lot of interest in the construction of embedded surfaces in~$4$-manifolds. 

Four is the lowest dimension in which there are manifolds that do not admit any smooth structure. Locally flat embedded surfaces are therefore the most we can hope to find in an arbitrary $4$-manifold, which may well be non-smoothable. In addition, there is a remarkable disparity between the smooth and topological settings in dimension four, in particular related to the behaviour of embedded surfaces. Thus, even in a smooth $4$-manifold, it is interesting to consider locally flat surfaces, e.g.~in order to detect when an invariant or phenomenon is `purely smooth' vs `purely topological'. 

In this expository article, we give an overview of the tools and techniques available to construct such purely topological embeddings, with the hope of emboldening more researchers to attack some of the many interesting open problems about locally flat surfaces in topological $4$-manifolds. 

Broadly speaking there are two types of techniques in this setting. First there are direct and hands-on methods, where we draw explicit pictures and modify them, keeping careful track of how intersection points are created or removed. This includes the manoeuvres in the constructive part of the proof of the disc embedding theorem (\Cref{thm:DET,thm:DET-full}). Describing these manoeuvres is the first goal of this paper. Specifically we will see them deployed in practice to prove the following theorem, due to Lee and Wilczy\'{n}ski.

\begin{theoremA}[{\cite[{$d=1$ case of Theorem~1.1}]{LW97}}]\label{thm:primitive-torus}
    Let $M$ be a closed, simply connected $4$-manifold. Then every primitive class in $H_2(M;\Z)$ is represented by a locally flat embedded torus. 
\end{theoremA}

Here a class is said to be \emph{primitive} if it is not a nonzero multiple of another class. The original proof of Lee and Wilczy\'{n}ski is not especially direct. We will give a more geometric proof from~\cite{KPRT:sigmet}. Theorem A is a shared special case of two distinct general results, from~\cite{LW97} and \cite{KPRT:sigmet}; we will state both in~\Cref{sec:direct}. 

In the second half of this paper we will use more abstract techniques from surgery theory to find locally flat embedded surfaces. The disc embedding theorem is the key reason why surgery theoretic techniques are available in dimension four, and they do not apply in the smooth setting. We will see how surgery theory can be used to show the following result due to Freedman and Quinn.  

\begin{theoremB}[{\cite[{Theorem~11.7B}]{FQ}}]\label{thm:alex-poly-one}
    Every knot $K\colon S^1\hookrightarrow S^3$ with Alexander polynomial one is $($topologically$)$ slice. 
\end{theoremB}

The treatment in this portion will be somewhat less detailed than the first half, relegating many ingredients to the exercises.   

\subsection*{Conventions}  Homeomorphism of manifolds is denoted by the symbol $\homeo$. Manifolds are \emph{not} assumed to be smooth. By definition submanifolds are locally flat (see \Cref{sec:defs-fund-tools}). Starting from \Cref{sec:direct}, all embeddings are assumed to be locally flat, although we will continue to specify this on occasion to try to avoid confusion.

Given surfaces $\Sigma$ and $\Sigma'$ intersecting transversely in some ambient $4$-manifold, we denote the set of intersection points by $\Sigma \pitchfork \Sigma'$, and the cardinality of this set by~$\lvert\Sigma \pitchfork \Sigma'\rvert$.

\subsection*{Outline}
In \Cref{sec:defs-fund-tools} we recall the precise definition of locally flat embeddings, before reviewing fundamental tools such as topological transversality and the existence of normal bundles. We also explain how to visualise embedded and immersed surfaces in $4$-manifolds, and introduce finger and Whitney moves. In \Cref{sec:DET} we give the first version of the disc embedding theorem along with some remarks on its history and more general versions. A detailed proof of Theorem A comprises the bulk of \Cref{sec:direct}, which ends with the statements of two more general results, of which Theorem A is a shared special case. \Cref{sec:alex-poly-one} gives the proof of Theorem B after providing the requisite background, including a $0$-surgery characterisation of sliceness, equivariant intersection and self-intersection numbers, and the sphere embedding theorem, and ends with a discussion of the surgery sequence and more general results proven with similar techniques. After a short conclusion in \Cref{sec:conclusion} a list of exercises is given in \Cref{sec:exercises}. 

\subsection*{Associated lectures} This paper was written to accompany a minicourse at the Georgia Topology Summer School in 2024, but may be read entirely independently. Indeed, many details and subtleties in the written account were elided in the lectures for convenience and brevity. The order of topics has also been slightly modified. Interested readers may find videos of the lectures online.  

\subsection*{Acknowledgements}
I would like to thank Akram Alishahi, Eduardo Fern\'{a}ndez Fuertes, 
David Gay, and Gordana Mati\'{c}, for organising the Georgia Topology Summer School 2024, where this minicourse took place, as well as the audience for their many questions, encouraging feedback, including poetry, and eager participation. I am grateful to Daniel Hartman, who was the TA for the accompanying problem sessions. My warm thanks go also to
Daniel Hartman, 
Patrick Orson, 
Mark Pencovitch,
and 
Mark Powell, 
as well as 
the anonymous referee and their anonymous PhD student,
for comments on previous drafts,  
and to Elise Brod and Megan Fairchild for their help with some of the figures.

\section{Definitions and fundamental tools}\label{sec:defs-fund-tools}

We begin this section by recalling the precise definition of locally flat embeddings. Next we review fundamental results for topological $4$-manifolds, such as the existence of normal bundles and topological transversality, due to Quinn~\cite{quinn:endsIII}, without which working in this setting would be nigh impossible. We next consider generic immersions, along with the \emph{immersion lemma}, which allows us to replace an arbitrary continuous map of a surface to a $4$-manifold by a generic immersion. We explain how to visualise locally flat embedded  and generically immersed surfaces in $4$-manifolds next. Finally, we briefly review Whitney moves and regular homotopies in the topological setting, stating \Cref{thm:generic-immersions-bijection} which gives a relationship between homotopy and regular homotopy for surfaces in a $4$-manifold.

The main results in this section (\Cref{thm:normal-bundles-transversality,thm:immersion-lemma}) were proven by Quinn~\cite{quinn:endsIII} and Freedman--Quinn~\cite{FQ}, using Freedman's disc embedding theorem (\Cref{thm:DET,thm:DET-full}) from~\cite{F}. We will not go into their proofs, which are quite intricate. Instead, we will be glad that these tools exist and use them freely in the rest of this paper. Analogous results hold for smooth maps of surfaces in smooth $4$-manifolds. These are often covered in introductory differential topology courses and the reader may well use them automatically without much thought. The takeaway of this section is that, at least with respect to normal bundles, transversality, and immersions, we can also be similarly confident about locally flat embedded or generically immersed surfaces in topological $4$-manifolds.

\subsection{Locally flat embeddings}\label{sec:loc-flat-embeddings}
For $m \geq 0$, let $\R^m_+ := \{(x_1,\dots,x_m) \in \R^m \mid x_1 \geq 0\}$. 

\begin{figure}[tb]
    \centering
    \begin{tikzpicture}
        \node[anchor=south west,inner sep=0] at (0,0){\includegraphics[width=10cm]{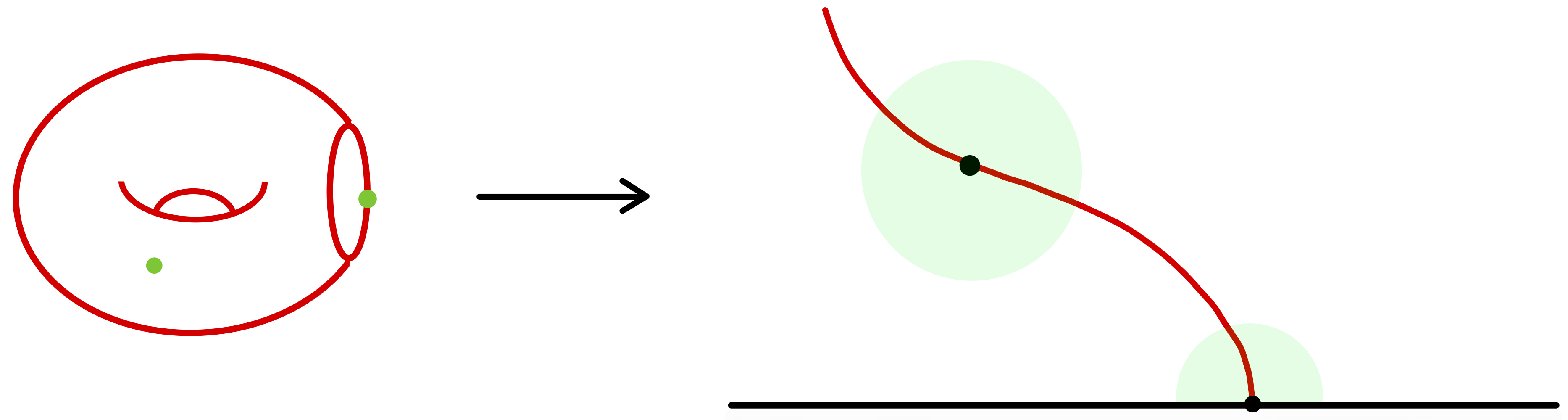}};
		\node at (1.15,1) {$x$};
		\node at (2.5, 1.4) {$y$};
		\node at (1.15,0.35) {$\Sigma$};
		\node at (2.4, 0.85) {$\partial \Sigma$};
		\node at (9,2) {$M$};
		\node at (5, 0.3) {$\partial M$};
		\node at (6.1,1.4) {$f(x)$};
		\node at (8.05, -0.15) {$f(y)$};
            \node at (7.1,2) {$U_x$};
		\node at (8.6, 0.5) {$U_y$};
            \node at (3.6,1.7){$f$};
            \node at (7.8,1.25){$f(\Sigma)$};
	\end{tikzpicture}  
    \caption{A locally flat embedding $f$ of a surface $\Sigma$ in a $4$-manifold $M$. Here we would have homeomorphisms $(U_x,U_x\cap f(\Sigma))\homeo (\R^4,\R^2)$ and $(U_y,U_y\cap f(\Sigma))\homeo (\R^4_+,\R^2_+)$.}
    \label{fig:loc-flat}
\end{figure}
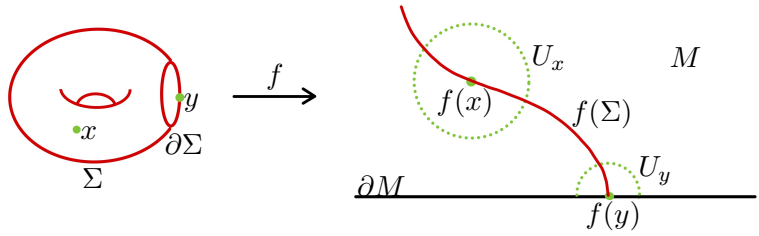

\begin{definition}\label{def:loc-flat-embedding}
Let $X$ be a $k$-manifold and let $M$ be an $n$-manifold with $k\leq n$. An embedding $f\colon X\hookrightarrow M$, i.e.~a continuous map which is a homeomorphism onto its image, is said to be \emph{locally flat} if for all $x\in X$ there is a neighbourhood $U\subseteq M$ of~$f(x)$ such that $(U,U\cap f(X))$ is homeomorphic to either $(\R^n,\R^k)$, if $x\in \Int{X}$, or to one of $(\R^n,\R^k_+)$ or $(\R^n_+, \R^k_+)$, if $x\in \partial X$, depending on whether $f(x)$ lies in $\Int{M}$ or $\partial M$, respectively. See the schematic in \Cref{fig:loc-flat}.

A locally flat embedding $f\colon X\hookrightarrow M$ is said to be \emph{proper} if its image is a closed subset of $M$, in the sense of point-set topology, and $f^{-1}(\partial M)= \partial X$.\footnote{In point-set topology, a continuous map of topological spaces is said to be \emph{proper} if the pre-image of every compact set is compact. This is distinct from our definition, which is sometimes called \emph{neat}. A relationship between the two notions is described in~\cite{lee:MSE-propervsproper}.}
\end{definition}

We point out that a locally flat embedding $f\colon X\hookrightarrow M$ is required to map interior points of $X$ to interior points of~$M$, but it is possible that $\partial X$ is mapped to $\Int{M}$. A proper locally flat embedding must map boundary points to boundary points.

For smooth $4$-manifolds one usually considers smooth embeddings. In the case of~$4$-manifolds which might not admit smooth structures, locally flat embeddings are the correct analogue. In particular, by definition, (proper) submanifolds of a topological manifold are images of (proper) locally flat embeddings. There do exist embeddings which are not locally flat (Exercise~\ref{ex:emb-not-loc-flat}). However, these lack some very useful properties enjoyed by locally flat embeddings, which we review next. We primarily discuss embeddings in $4$-manifolds, relegating other dimensions to \Cref{rem:other-dims}. For a more detailed survey, we direct the reader to \cite{FNOP:4dguide,DET-book-flowchart}.

\begin{theorem}[{\cite[Theorem~2.4.1]{quinn:endsIII},~\cite[Theorem]{quinn:transversality}, and~\cite[Theorems~9.3 and 9.5A]{FQ}}]\label{thm:normal-bundles-transversality}
    Let $M$ be a $4$-manifold. 
    \begin{enumerate}
        \item\label{item:normal-bundles} (Existence of normal vector bundles) Every proper submanifold of $M$ has a normal vector bundle, which is unique up to bundle isomorphism and ambient isotopy.
        \item\label{item:top-transversality} (Topological transversality) Let $X_1$ and $X_2$ be proper submanifolds of $M$. Then there is an ambient isotopy of $M$ taking $X_1$ to some $X_1'$ such that $X_1'$ and $X_2$ intersect transversely.
    \end{enumerate}
\end{theorem}

Without going into too many details, we recall the definition of a normal vector bundle. 

\begin{definition}\label{def:normal-bundle-submanifold}
    Let $M$ be a $4$-manifold and let $(X,\partial X)\subseteq (M,\partial M)$ be a $k$-dimensional proper submanifold. A \emph{normal vector bundle} of $X$ in $M$ is a pair~$(E, p\colon E \to X)$ with the following properties.
    \begin{enumerate}
        \item The total space $E$ is a neighbourhood of $X$ in $M$ and a codimension zero proper submanifold of~$M$;
        \item the map $p\colon E \to X$ is an $(4-k)$-dimensional vector bundle such that $p(x) = x$ for all $x\in X$;
        \item the boundary $\partial E=p^{-1}(\partial X)\subseteq \partial M$; and 
        \item the data above are \emph{extendable}, i.e.~given any $(4-k)$-dimensional vector bundle~$q\colon F \to X$, any radial homeomorphism from an open convex disc bundle of~$F$ to $E$ can be extended to a homeomorphism from all of $F$ to a neighbourhood of $E$ in $M$. 
    \end{enumerate}
\end{definition}
Given a vector bundle $q\colon F\to X$, an \emph{open convex disc bundle} of $F$ is a bundle arising from restricting each fibre to an open convex disc. 

The purpose of the first three properties in \Cref{def:normal-bundle-submanifold} is for the normal vector bundle to mimic the notion of an open tubular neighbourhood in the smooth setting. There is a technical problem that the closure of such an open neighbourhood might have undesirable self-intersections. The fourth property of extendability is designed to avoid this.

We now recall the definition of transversality. 

\begin{definition}
Let $(X_1,\partial X_1),(X_2,\partial X_2)\subseteq (M,\partial M)$ be proper submanifolds of a~$4$-manifold $M$, of dimension $k_1$ and $k_2$ respectively. We say that $X_1$ and $X_2$ \emph{intersect transversely} if the following conditions hold. 
    \begin{enumerate}
        \item If $k_1+k_2<4$ then $X_1$ and $X_2$ are disjoint. 
        \item If $k_1+k_2=4$, then $\partial X_1$ and $\partial X_2$ are disjoint and each point $x\in X_1\cap X_2$ has a neighbourhood $U\subseteq M$ such that 
    \[
    (U,U\cap X_1,U\cap X_2)\homeo
         \big(\R^4, \R^{k_1}\times \{0\}, \{0\}\times \R^{k_2}\big).
    \]
    Note that this homeomorphism maps $x\mapsto 0$.
        \item If $k_1+k_2> 4$, then each point $x\in X_1\cap X_2$ has a neighbourhood $U\subseteq M$ such that 
    \[
    \leavevmode
    (U,U\cap X_1,U\cap X_2)\homeo\left\{
    \begin{array}{@{}lr@{}}
        \big(\R^4, \R^{k_1}\times \{0\}, \{0\}\times \R^{k_2}\big),    &\text{if }x\in \Int{M};\\
        \big(\R^4_+,\R^1_+\times \R^{k_1-1}\times \{0\} , \R^1_+\times \{0\}\times \R^{k_2-1}\big),  &\text{if }x\in \partial M.
    \end{array}
    \right.
    \]
    Note that these homeomorphisms map $x\mapsto 0$.
    \end{enumerate}
\end{definition}

\begin{figure}[tb]
    \centering
    \begin{tikzpicture}
        \node[anchor=south west,inner sep=0] at (0,0){ \includegraphics[width=7cm]{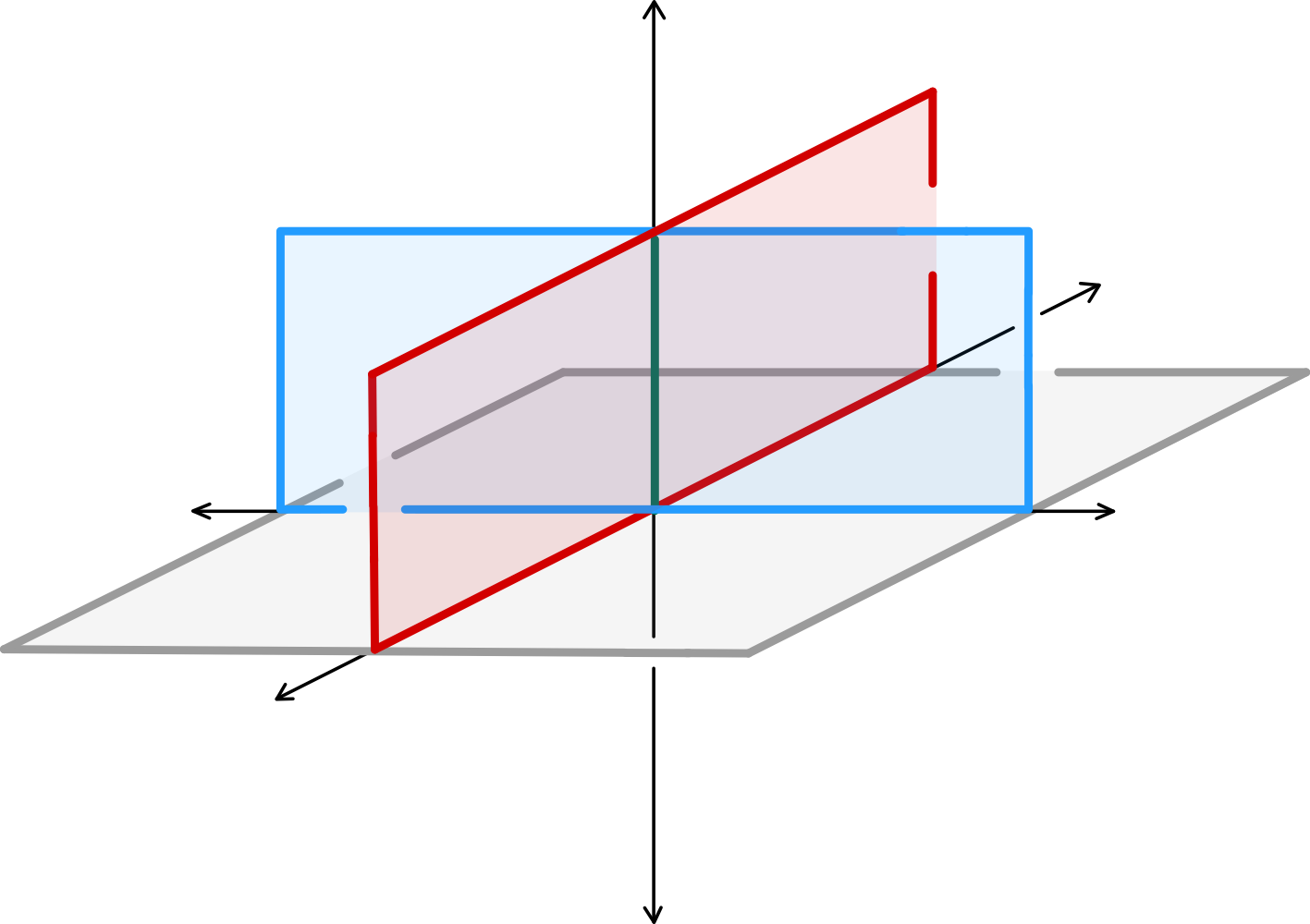}};
            \node at (3.75,4.8) {$x_1$};
		\node at (1.15, 2.35) {$x_2$};
		\node at (1.7,1.1) {$x_3$};
    \end{tikzpicture}  
    \caption{Within the upper half space $\R^3_+:=\{(x_1,x_2,x_3)\mid x_1\geq 0\}$, the half planes $\R^1_+\times \R\times 0$ and $\R^1\times 0\times \R$ are shown in blue and red respectively. They intersect transversely along $\R^1_+$. We show the intersection in green.}
    \label{fig:transverse-at-boundary}
\end{figure}

For the final condition above where $x\in \partial M$, see \Cref{fig:transverse-at-boundary}.
Intuitively, you should think of topological transversality (\Cref{thm:normal-bundles-transversality}\,\eqref{item:top-transversality}) as saying that, given a pair of submanifolds of a $4$-manifold, we can isotope one of them so that it now intersects the other one in the smallest possible dimension. 

\begin{remark}\label{rem:other-dims}
    An oft-repeated slogan is that topological $n$-manifolds, with $n\geq 4$ behave like high-dimensional smooth manifolds\footnote{To lightly paraphrase Andrew Ranicki~\cite{ranicki:orsay-slides}: Topological manifolds of dimension $\geq 4$ are just like DIFF and PL manifolds, only more so!} (whereas smooth $4$-manifolds do not). Indeed, there are situations where topological $4$-manifolds are even better behaved than higher-dimensional topological manifolds. As an example of this, we note that (locally flat) submanifolds of high-dimensional manifolds do not necessarily have normal vector bundles. For more on this, see~\cite[Section~9.4]{FQ} and \cite[Example 6.23]{FNOP:4dguide}. 
    
    On the other hand, topological transversality holds in all dimensions and codimensions, but the definition is much more complicated to parse, in particular due to the unavailability of normal vector bundles. The main results in this setting are due to Marin~\cite{marin:transversality} (cf.~\cite[Essay~III, Section~1, p.~83]{KS}) and Quinn~\cite{quinn:endsIII,quinn:transversality} (see also~\cite[Section~9.5]{FQ} and \cite[Chapter~10]{FNOP:4dguide}). 

    We also remark that working in dimension four gives us access to a variety of geometric tools (not to mention intuition), unavailable to our colleagues studying higher-dimensional manifold theory. Broadly and loosely speaking, working on topological $4$-manifolds allows us to use the algebraic machinery of higher-dimensional manifolds, often with more tractable algebraic and homotopy-theoretic problems to solve, while also benefitting from lower-dimensional geometric technology.    
    \end{remark}

\subsection{Topological generic immersions}
In addition to locally flat embeddings, there is also a useful notion of a \emph{generic immersion}, and a result saying that continuous maps can be approximated by these. To state this precisely, we first give the definition of an immersion of manifolds in the topological setting. For $k\leq n$, we have the following standard inclusions, which we used implicitly in \Cref{sec:loc-flat-embeddings}.
\begin{align}
  \iota \colon \R^k_{\phantom{+}} &= \R^k_{\phantom{+}} \times \{0\} \hooklongrightarrow \R^k_{\phantom{+}} \times \R^{n-k} = \R^n; \\
  \iota_\mathsmaller{+} \colon \R^k_+ &= \R^k_+ \times \{0\} \hooklongrightarrow \R^k_{\phantom{+}} \times \R^{n-k} = \R^n; \text{ and }\\
  \iota_\mathsmaller{++} \colon \R^k_+ &= \R^k_+ \times \{0\} \hooklongrightarrow \R^k_+ \times \R^{n-k} = \R^n_+.
\end{align}

\begin{definition}\label{def:top-generic-immersion}
    Let $X$ be a $k$-manifold and let $M$ be an $n$-manifold with $k\leq n$. A continuous map $f \colon X \to M$ is said to be a (topological or locally flat) \emph{immersion} if for each point $x\in X$ there is a chart~$\varphi$ around $x$ and a chart $\Psi$ around $f(x)$ fitting into one of the following commutative diagrams. The first diagram is for $x \in \Int X$ and $f(x) \in \Int M$; the second diagram is for $x \in \partial X$ and $f(x) \in \Int M$; and the third is for $x \in \partial X$ and $f(x) \in \partial M$. In particular $f$ is required to map interior points of $X$ to interior points of~$M$, but it is possible that $\partial X$ is mapped to $\Int{M}$. 
    \begin{equation}\label{diagram:above}
        \begin{gathered}
            \xymatrix{
            \R^k \ar[r]^{\iota} \ar[d]^{\varphi} & \ar[d]^{\Psi}  \R^n   \\
            X \ar[r]^f & M }
            \hspace{2em}
            \xymatrix{
            \R^k_+ \ar[r]^{\iota_+} \ar[d]^{\varphi} & \ar[d]^{\Psi}  \R^n   \\
            X \ar[r]^f & M }
            \hspace{2em}
            \xymatrix{
            \R^k_+ \ar[r]^{\iota_{++}} \ar[d]^{\varphi} & \ar[d]^{\Psi}  \R^n_+   \\
            X \ar[r]^f & M }
        \end{gathered}
    \end{equation}
    
    The \emph{singular set} of an immersion $f \colon X \to M$ is the set \[\mathcal{S}(f) := \{m \in M \mid |f^{-1}(m)| \geq 2\}.\]
\end{definition}

In other words, an immersion is a local, locally flat embedding.
As in the smooth setting, there is a notion of normal vector bundles for immersions. We recall the definition next.

\begin{definition}\label{def:normal-bundle-immersion}
Let $X$ be a $k$-manifold and let $M$ be an $n$-manifold with $k\leq n$. A \emph{normal vector bundle} for an immersion $f\colon X \to M$ is an $(n-k)$-dimensional real vector bundle $\pi \colon \nu f \to X$, together with an immersion $\wt{f} \colon \nu f \to M$ that restricts to $f$ on the zero section $s_0$, i.e.\ $\wt f\circ s_0 = f$, and such that each point $x\in X$ has a neighbourhood $U$ such that $\wt f|_{\pi^{-1}(U)}$ is a locally flat embedding. As in \Cref{def:normal-bundle-submanifold}, we further require these data to be extendable.
\end{definition}

Now we restrict to the case of surfaces mapping to $4$-manifolds, which is most relevant for us. 

\begin{definition}\label{def:gen-immersion}
     Let $\Sigma$ be a surface and let $M$ be a $4$-manifold. A continuous map~$f \colon \Sigma \to M$ is said to be a \emph{generic immersion}, denoted by $f\colon \Sigma\looparrowright M$, if it is an immersion and the singular set $\calS(f)$ is a closed, discrete subset of $M$ consisting only of transverse double points, each of whose preimages lies in the interior of $\Sigma$. In particular, whenever $m \in \mathcal{S}(f)$, there are exactly two points $p_1, p_2 \in \Sigma$ with $f(p_1)=m=f(p_2)$, and there are disjoint charts $\varphi_i$ around $p_i$, for $i=1,2$, where~$\varphi_1$ is as in the left-most diagram of \eqref{diagram:above}, and $\varphi_2$ is the same, with respect to the same chart $\Psi$ around $m$, but with $\iota$ replaced by
    \[
    \iota' \colon \R^2 = \{0\} \times \R^2 \hooklongrightarrow \R^2 \times \R^2 = \R^4.
    \]
\end{definition}

Given a generic immersion $f\colon \Sigma\looparrowright M$ of a surface $\Sigma$ in a $4$-manifold $M$, it is customary to use the symbol $f$ to denote the image of $\Sigma$. Note that a generic immersion $f\colon \Sigma\looparrowright M$ necessarily restricts to a locally flat embedding on $\partial \Sigma$. 

As mentioned previously, as in the smooth category, arbitrary continuous maps of a surface to a $4$-manifold may be replaced by generic immersions, by the following result.
 
\begin{theorem}[{Immersion lemma~\cite[{Corollary~9.5C}]{FQ}}]\label{thm:immersion-lemma}
    Let $\Sigma$ be a surface and let $M$ be a $4$-manifold. Every continuous map $f\colon \Sigma\to M$ is homotopic to a generic immersion. 
    
    If $f$ is already a generic immersion in a neighbourhood of $\partial \Sigma$, then the homotopy can be chosen to be constant on $\partial \Sigma$. 
\end{theorem}

We remark that we allow generic immersions to map the boundary of a surface to the interior of a $4$-manifold, since we will often apply the immersion lemma to find generically immersed Whitney discs, whose boundaries usually lie in the interior of the ambient $4$-manifold. 

Generic immersions admit particularly nice normal vector bundles, as we see in the following result. We will need their existence in \Cref{sec:proof-primitive-torus} in some of our geometric manoeuvres on Whitney discs. 

\begin{theorem}[{\cite[{Theorem~2.4}]{KPRT:sigmet}}]\label{thm:plumbed-normal-bundle}
Let $\Sigma$ be a surface and let $M$ be a~$4$-manifold. A generic immersion $f \colon \Sigma \looparrowright M$ has a normal vector bundle $\pi\colon \nu f\to \Sigma$ as in \Cref{def:normal-bundle-immersion}, with the additional property that the immersion $\wt f\colon \nu f\to M$ is an embedding outside a neighbourhood of $f^{-1}(\mathcal{S}(f))$, and near the double points plumbs two coordinate regions $\pi^{-1}(\varphi_i(\R^2))\homeo \varphi_i(\R^2)\times \R^2$, $i=1,2$, together i.e.~$\wt f\circ ( \varphi_1(x),y) = \wt f\circ (\varphi_2(y),x)$, using the notation of \Cref{def:gen-immersion}.
\end{theorem}

\subsection{Visualising surfaces in 4-manifolds}

In \Cref{sec:direct}, we will primarily modify generically immersed surfaces directly by hand. Therefore it will be crucial for us to visualise them. We will generally draw schematic pictures, but we begin with a few concrete ones.

By definition locally flat and generically immersed surfaces in an arbitrary $4$-manifold are standard in small coordinate charts, which we can draw precisely. Since each chart in a $4$-manifold is a copy of $\R^4=\R^3\times \R$, we can draw a sequence of pictures of $\R^3$, and see how our surfaces show up within them. 

Let $x$, $y$, and $z$ denote the usual Cartesian coordinates in $\R^3$, and let $t$ denote the fourth coordinate in $\R^4$. This fourth coordinate is usually thought of as representing time, so that the corresponding copies of $\R^3$ can be `played', either backwards or forwards, like in a movie. 
\begin{figure}[tb]
    \centering
    \begin{tikzpicture}
        \node[anchor=south west,inner sep=0] at (0,0){\includegraphics[width=10cm]{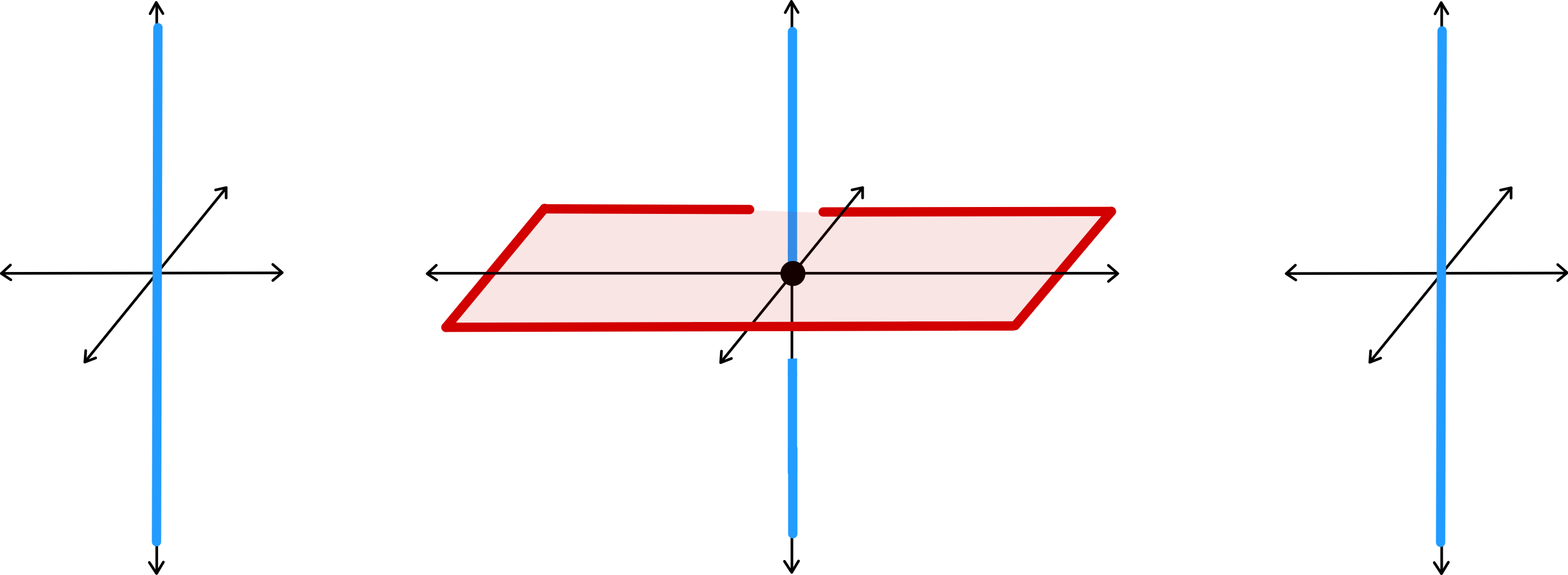}};
            \node at (1,-0.35) {$t=-\varepsilon$};
            \node at (5.1,-0.35) {$t=0$};
            \node at (9.2,-0.35) {$t=\varepsilon$};
            \node at (4.5,1.4) {$x$};
            \node at (7,1.75) {$y$};
            \node at (5.2,3.5) {$z$};
            \node at (0.425,1.4) {$x$};
            \node at (1.6,1.75) {$y$};
            \node at (1.15,3.5) {$z$};
            \node at (8.625,1.4) {$x$};
            \node at (9.8,1.75) {$y$};
            \node at (9.35,3.5) {$z$};
	\end{tikzpicture}  
    \caption{For each value of $t$, we see a corresponding copy of $\R^3$. The axes are shown in grey. The $xy$-plane is shown in red, and the $zt$-plane is shown in blue.}
    \label{fig:transverse-movie-1}
\end{figure}
In \Cref{fig:transverse-movie-1}, we depict a region of $\R^4$ centred at the origin. Note that we get a copy of $\R^3$ in each subfigure, as $t$ varies from $-\varepsilon$ to $\varepsilon$. The red plane in the central subfigure is the $xy$-plane. The $zt$-plane is depicted in blue -- for each value of $t$, we only see a line in the corresponding copy of $\R^3$, however these lines trace out the entire plane as we move backwards and forwards in time. Note that the blue and red planes intersect at a unique point, namely the origin, as expected. 

By definition, given a surface $\Sigma$, a $4$-manifold $M$, and a locally flat embedding~$f\colon \Sigma \hookrightarrow M$, for every point $p\in \Sigma$, there is an open set $U\subseteq M$ with $f(p)\in U$ and a homeomorphism $U\homeo \R^4$ such that $f(\Sigma)\cap U$ is mapped to the $xy$-plane. Similarly, given surfaces $\Sigma_1,\Sigma_2\subseteq M$ intersecting transversely at some point $q\in M$, by definition there is an open set $U\subseteq M$ with $q\in U$ and a homeomorphism $U\homeo \R^4$, mapping $\Sigma_1\cap U$ to the $xy$-plane and $\Sigma_2\cap U$ to the~$zt$-plane. In particular, the point of intersection $q$ is mapped to the origin in $\R^4$. The same holds for a self-intersection of a generically immersed surface. In other words, \Cref{fig:transverse-movie-1} gives a concrete, if local, picture of either a generically immersed surface in a $4$-manifold, or a pair of transversely intersecting locally flat surfaces in a $4$-manifold. 

\begin{remark}
    Now that we have explicit local pictures of generically immersed and transversely intersecting surfaces in a $4$-manifold, we can use these to concretely describe manoeuvres on them. In \Cref{sec:proof-primitive-torus} we will draw several such `movie pictures'. In the literature it is customary to only draw the middle time slice, leaving it to the reader to draw (or imagine) the extensions to the remaining time slices. We will try to draw all time slices in general, but we bend to tradition in \Cref{fig:boundary-pushoff,fig:transfer-move,fig:BFF-move}. 
\end{remark}

\begin{figure}[tb]
    \centering
    \begin{tikzpicture}
        \node[anchor=south west,inner sep=0] at (0,0){\includegraphics[width=9cm]{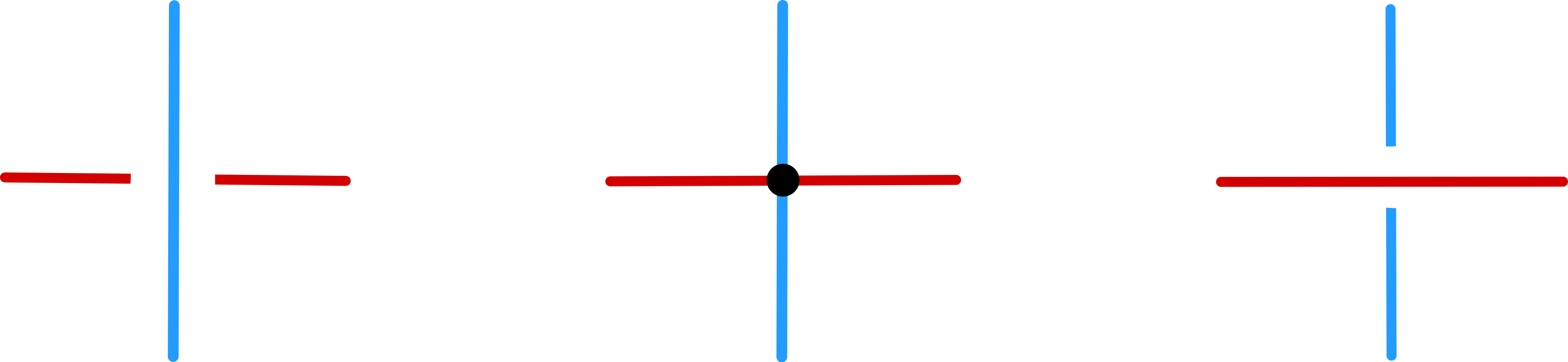}};
            \node at (1,-0.35) {$t=-\varepsilon$};
            \node at (4.5,-0.35) {$t=0$};
            \node at (8,-0.35) {$t=\varepsilon$};
	\end{tikzpicture}  
    \caption{Two surfaces, shown in red and blue respectively, intersect transversely at a single point.}
    \label{fig:transverse-movie-2}
\end{figure}

One might rightfully complain that \Cref{fig:transverse-movie-1} is not especially symmetric, since one surface is shown entirely in a single time slice, while the second surface is smeared across multiple times. A more symmetric (local) depiction of a transverse point of intersection between two locally flat surfaces in a $4$-manifold (or potentially a generic self-intersection of a single surface) is shown in \Cref{fig:transverse-movie-2}. 
In this case both surfaces, shown in red and blue respectively, appear as a single line in each time slice. As the `movie' is played, these lines trace out the corresponding surfaces.

\begin{figure}[tb]
    \centering
    \begin{tikzpicture}
        \node[anchor=south, inner sep=0] at (5.5,0){\includegraphics[width=10.25cm]{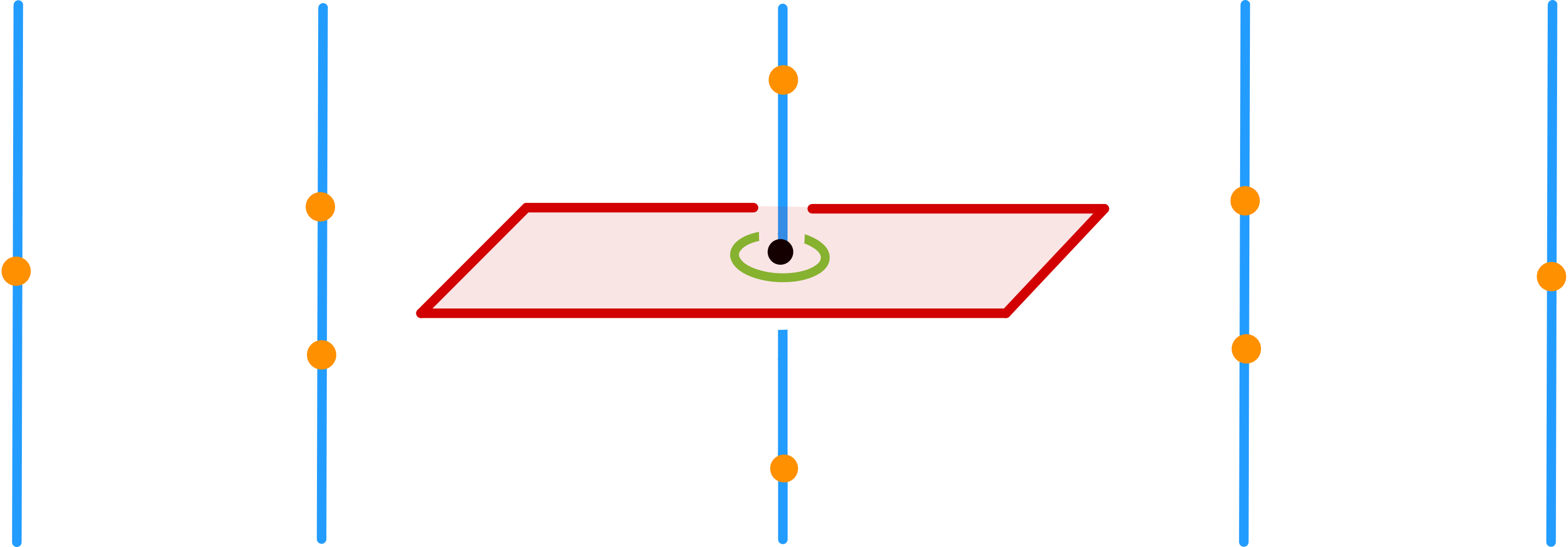}};
            \node at (0.4,-0.35) {$t=-\varepsilon$};
            \node at (2.5,-0.35) {$t=-\tfrac{\varepsilon}{2}$};
            \node at (5.5,-0.35) {$t=0$};
            \node at (8.6,-0.35) {$t=\tfrac{\varepsilon}{2}$};
            \node at (10.6,-0.35) {$t=\varepsilon$};
	\end{tikzpicture} 
    \caption{The unit circle in the $xy$-plane is shown in green and the unit circle in the~$zt$-plane is shown in orange.}
    \label{fig:movie-unit-circles}
\end{figure}

\begin{figure}[tb]
    \centering
    \begin{tikzpicture}
        \node[anchor=south,inner sep=0] at (5.5,0){\includegraphics[width=11cm]{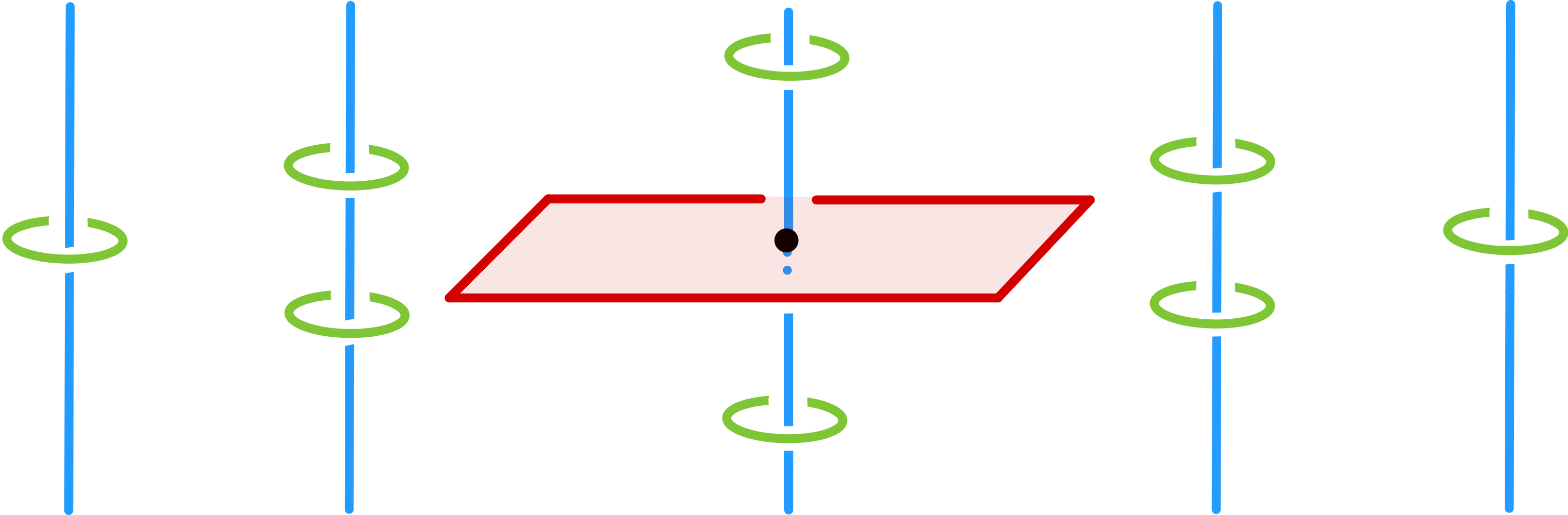}};
            \node at (0.4,-0.35) {$t=-\varepsilon$};
            \node at (2.5,-0.35) {$t=-\tfrac{\varepsilon}{2}$};
            \node at (5.5,-0.35) {$t=0$};
            \node at (8.6,-0.35) {$t=\tfrac{\varepsilon}{2}$};
            \node at (10.6,-0.35) {$t=\varepsilon$};
            \end{tikzpicture} 
    \caption{The Clifford torus is shown in green.}
    \label{fig:movie-clifford-torus}
\end{figure}

Let us take a moment to find the \emph{Clifford torus} in \Cref{fig:transverse-movie-1}, since we will use it in the proof of \Cref{prop:whitney-by-DET}. By definition, the Clifford torus is the product (in~$\R^4$) of the unit circle in the $xy$-plane with the unit circle in the $zt$-plane. Now that we have a concrete picture of $\R^4$ in \Cref{fig:transverse-movie-1}, we can draw the Clifford torus easily. We do so in \Cref{fig:movie-unit-circles,fig:movie-clifford-torus}. 

We already argued that \Cref{fig:transverse-movie-1} is a concrete picture of a small neighbourhood of a transverse point of intersection between two surfaces in a $4$-manifold, or a self-intersection of a generically immersed surface. Therefore, we can now find a Clifford torus at any such intersection point.

\subsection{Finger moves and Whitney moves}\label{sec:finger-whitney}

\begin{figure}[tb]
    \centering
    \begin{tikzpicture}
        \node[anchor=south,inner sep=0] at (5.75,8.5){\includegraphics[width=11.5cm]{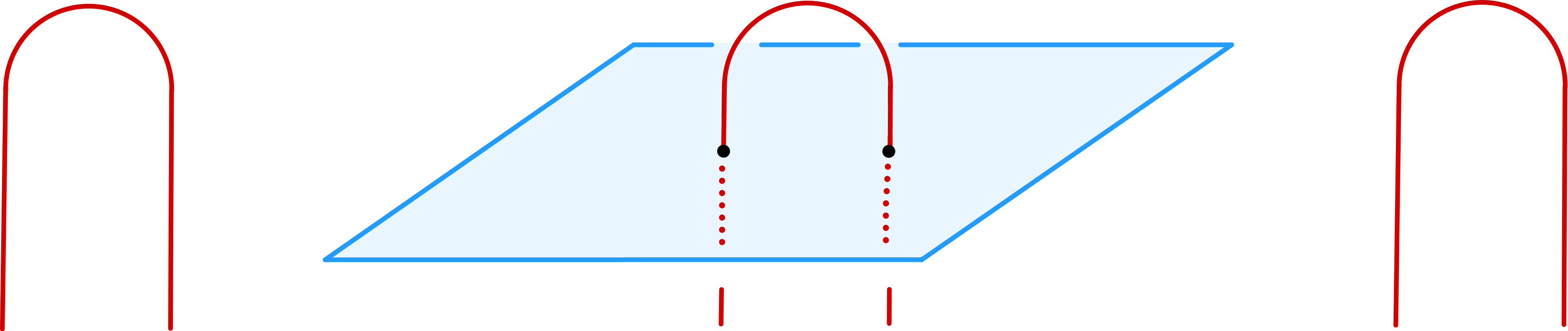}};
        \node[anchor=south,inner sep=0] at (5.75,4.25){\includegraphics[width=11.5cm]{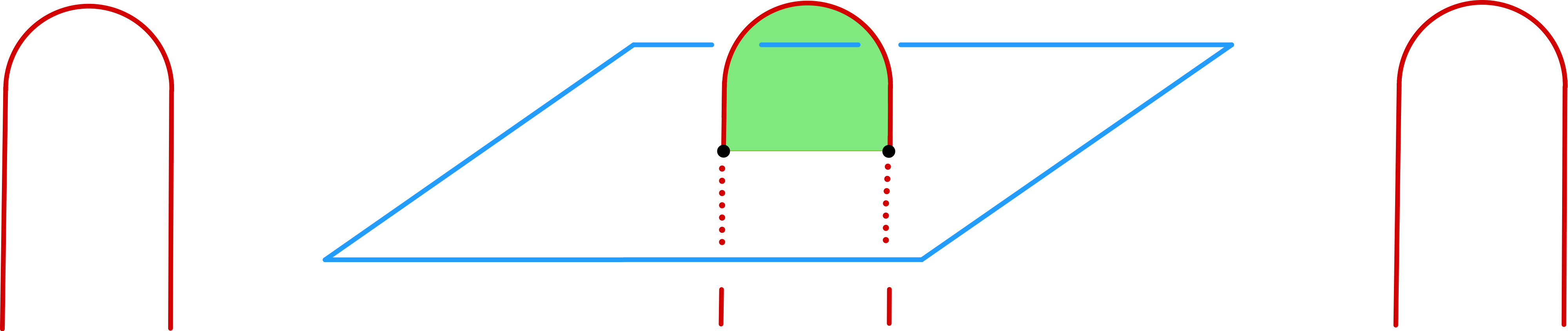}};
        \node[anchor=south,inner sep=0] at (5.75,0){\includegraphics[width=11.5cm]{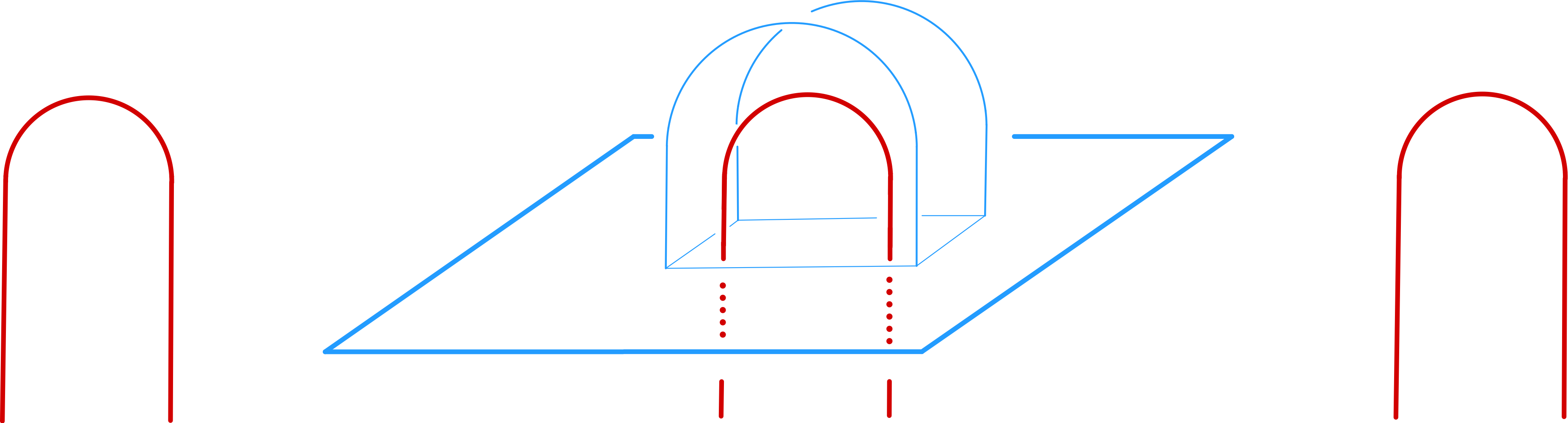}};
		\node at (5.75,7.75) {(a)};
		\node at (0.625,8.25) {$t=-\varepsilon$};
		\node at (5.75,8.25) {$t=0$};
		\node at (10.875,8.25) {$t=\varepsilon$}; 
		\node at (2.75,8.75) {$f$};    
            \node at (0.2,8.7) {$g$};
		\node at (5.45,8.7) {$g$};
		\node at (10.4,8.7) {$g$}; 
            \node at (5.75,3.5) {(b)};
		\node at (0.625,4) {$t=-\varepsilon$};
		\node at (5.75,4) {$t=0$};
		\node at (10.875,4) {$t=\varepsilon$}; 
		\node at (2.75,4.5) {$f$};    
            \node at (0.2,4.45) {$g$};
		\node at (5.45,4.45) {$g$};
		\node at (10.4,4.45) {$g$}; 
            \node at (5.9,6) {$W$};
		\node at (5.9,5.45) {{\small $A$}};
		\node at (6.25,6.75) {{\small $B$}};
            \node at (5.75,-0.75) {(c)};
		\node at (0.625,-0.25) {$t=-\varepsilon$};
		\node at (5.75,-0.25) {$t=0$};
		\node at (10.875,-0.25) {$t=\varepsilon$}; 
		\node at (2.75,0.25) {$f$};    
            \node at (0.2,0.2) {$g$};
		\node at (5.45,0.2) {$g$};
		\node at (10.4,0.2) {$g$};   
	\end{tikzpicture}  
    \caption{(a) Two surfaces $f$, shown in blue, and $g$, shown in red, intersect at two points marked in black. (b) A Whitney disc $W$ is shown in green. The Whitney arcs are $A$, along $f$, and $B$ along $g$. (c) The Whitney move on $f$ along $W$ has replaced a small strip neighbourhood of $A$ on $f$ with the union of two pushed off copies of $W$, along with a strip whose core is parallel to $B$. We do not shade the surface $f$ here for clarity. }
    \label{fig:whitney-move}
\end{figure}

We will not describe finger moves and Whitney moves in detail, referring instead to existing sources in the literature, such as \cite[Chapter~1]{FQ} and \cite{DET-book-DETintro}. Since Whitney discs will be the main object of our various geometric manoeuvres in \Cref{sec:direct}, we describe them briefly, relying primarily on \Cref{fig:whitney-move}. 
As depicted in the figure, we consider two transverse intersection points of opposite sign between generically immersed connected surfaces $f$ and $g$ in some ambient $4$-manifold $M$, where possibly $f=g$. We see that the two points can be joined by two arcs, called \emph{Whitney arcs}, one lying on $f$, denoted by $A$, and one on $g$, denoted by $B$. The union~$A\cup B$ is called a \emph{Whitney circle}. A disc in $M$ bounded by the Whitney circle is called a \emph{Whitney disc}. Suppose we have a generically immersed Whitney disc~$W$. Then the normal vector bundle of $W$, being a bundle on a contractible space, is trivial. Consider the restriction of this trivial $2$-plane bundle to the Whitney circle. We can define a $1$-plane subbundle by choosing vectors in the $f$-direction along~$A$, and vectors normal to $g$ along $B$. Let $s$ denote a section of this subbundle. The Whitney disc $W$ is said to be \emph{untwisted} if $s$ admits a nonvanishing extension to the entire normal vector bundle over $W$.\footnote{In the literature such Whitney discs are sometimes called \emph{framed}. We do not like this terminology since in general a trivial bundle is said to be framed if a trivialisation has been chosen. Note that bundles over discs are trivial and uniquely trivialisable, since discs are contractible.} Not all Whitney discs that we consider will be untwisted. The \emph{twisting number} of $W$, denoted by $\tw(\partial W)$, is the signed count of zeros of $s$ when extended over the normal vector bundle over all of $W$.\footnote{This is the relative euler number of the normal vector bundle of $W$ with respect to the \emph{Whitney framing} on the boundary described earlier.} 

The \emph{Whitney move} on $f$ along $W$ consists of replacing a small strip neighbourhood of $A$ on $f$ with two copies of $W$, pushed off along the sections $s$ and $-s$ respectively, union a small strip whose core is parallel to $B$. This procedure is described in \Cref{fig:whitney-move}. Note that if $W$ is embedded, untwisted, and has interior disjoint from $f$ and $g$, then the Whitney move on $W$ removes the two intersection points between~$f$ and $g$ being paired by $W$ without creating any new intersections. If a Whitney disc fails any one of these three conditions (no self-intersections, zero twisting, and having interior disjoint from $f$ and $g$), then one may still perform the Whitney move. However, while the original pair of intersections of $f$ and $g$ will be removed, new, undesirable intersections will be created. The bulk of the proof of Theorem A consists of arranging for these three conditions on Whitney discs. 

\subsection{Regular homotopies}
Recall that in the smooth setting a \emph{regular homotopy} is by definition a homotopy through immersions. A smooth regular homotopy of a generically immersed surface~$f$ in a $4$-manifold is generically a concatenation of (smooth) isotopies, finger moves, and Whitney moves along untwisted, embedded, and disjoint Whitney discs, with interiors disjoint from $f$ (see e.g.~\cite[{Section~III.3}]{GoGu}). This fact inspires the definition of the topological analogue. 

\begin{definition}\label{def:reg-htpy}
    A (topological) \emph{regular homotopy} of a generically immersed surface $f$ in a $4$-manifold is a concatenation of (topological) isotopies, finger moves, and Whitney moves along untwisted, embedded, and disjoint Whitney discs, with interiors disjoint from $f$.
\end{definition}

We take this opportunity to state a theorem relating homotopies and regular homotopies. 

\begin{theorem}[{\cite[special case of Theorem~2.32]{KPRT:sigmet}}]\label{thm:generic-immersions-bijection}
	Let $\Sigma$ be a closed, connected surface and let $M$ be an orientable $4$-manifold. Then the inclusion of the subspace of generic immersions $\mathrm{GenImm}(\Sigma,M)$ in the space $[\Sigma,M]$ of all continuous maps from $\Sigma$ to $M$ induces a map
	\[
	\frac{\mathrm{GenImm}(\Sigma,M)} { \{\text{regular homotopy}\}} \xrightarrow{\phantom{5}i\phantom{5}} [\Sigma,M], \]
    with the following properties.
	\begin{enumerate}
		\item\label{thm:gib-i} By \Cref{thm:immersion-lemma}, the map $i$ is surjective.
		\item\label{thm:gib-ii} The fibres of $i$ are related by \emph{interior twisting}. More precisely, suppose that $f$ and $g$ are homotopic generic immersions. Then we can perform interior twisting on $f$ and on $g$, to obtain $f'$ and $g'$ respectively, such that $f'$ and $g'$ are regularly homotopic. See \Cref{sec:proof-primitive-torus} for a description of interior twisting. 
		\item\label{thm:gib-iii} For every $f\in [\Sigma,M]$, there is a bijection
			\[i^{-1}(f)\cong \begin{cases}
			2\Z,&\text{if }w_2(\nu\wt f)=0;\\[1ex]
			2\Z+1,&\text{if } w_2(\nu\wt f)=1;\\[1ex]
		\end{cases}\]
		where $\nu\wt f$ is a normal vector bundle for $\wt f$, a generic immersion in $i^{-1}(f)$.
		This bijection is given by \[\wt f\mapsto e(\nu \wt{f}),\]
        where $e(\nu \wt{f})$ denotes the euler number of $\nu\wt{f}$. 
	\end{enumerate}
\end{theorem}

See~\cite[Theorem~2.32]{KPRT:sigmet} for a more general statement, including when $\Sigma$ has nonempty boundary and when $M$ is non-orientable.

\section{The disc embedding theorem}
\label{sec:DET}

The fundamental breakthrough in the study of topological $4$-manifolds, and surfaces within them, was due to the \emph{disc embedding theorem}. We begin by stating the simplest version of the theorem, and address more general versions in subsequent remarks. We also give the most general known statement later in \Cref{thm:DET-full}. Below we use $\cdot$ to denote the homological intersection pairing, on either a pair of absolute second homology classes, or a pair consisting of one absolute and one relative second homology class, in a $4$-manifold.

\begin{theorem}[Disc embedding theorem, simplest version; \cite{F} and {\cite[Theorem~5.1A]{FQ}}]\label{thm:DET}
    Let $M$ be a simply connected topological $4$-manifold. Suppose we have a generic immersion
    \[
    \begin{tikzcd}[
        /tikz/column 1/.append style={anchor=base east, column sep=-10pt, inner xsep =0pt}, 
        /tikz/column 3/.append style={column sep=0pt}, 
        /tikz/column 4/.append style={anchor=base west, column sep=0pt, inner xsep =0pt}
        ]
    f\colon &D^2\ar[loop-math to,r]   &M\\
    &\partial D^2\ar[hook,u]\ar[r,hook]  &\partial M\ar[hook,u],
    \end{tikzcd}
    \]
    where $f\vert_{\partial D^2}$ is a locally flat embedding and the vertical maps are inclusions. Suppose further that there is a generic immersion $g\colon S^2\looparrowright M^4$, such that 
    \begin{enumerate}[label=(\roman*)]
        \item the generic immersion $g$ has trivial normal vector bundle and trivial self-intersection, i.e.~$g\cdot g=0$; and 
        \item the maps $f$ and $g$ are \emph{algebraically dual} i.e.~$f\cdot g=1$.
    \end{enumerate}
    Then there is a locally flat embedding $\ol{f}\colon D^2\hookrightarrow M$ such that $f\vert_{\partial D^2}=\ol{f}\vert_{\partial D^2}$ and~$g$ is homotopic to a generic immersion $\ol{g}$, such that 
    $\ol{f}$ and $\ol{g}$ are \emph{geometrically dual}, i.e.~$\ol{f}$ and $\ol{g}$ intersect each other transversely and at a single point.
\end{theorem}

\begin{remark}\label{rem:DET-finite-collection}
    There is a version of the theorem for finite collections of discs~\cite[{Theorem~5.1A}]{FQ} (see \Cref{thm:DET-full}). The proof is essentially the same. There is a complicated generalisation to infinite collections of discs, called the \emph{disc deployment lemma}, which is significantly harder to prove~\cite[{Lemma~3.2}]{quinn:endsIII}.
\end{remark}

\begin{remark}\label{rem:DET-nontrivial-pi1}
    The theorem also holds for ambient $4$-manifolds with more general fundamental group~\cite[{Theorem~5.1A}]{FQ} (see \Cref{thm:DET-full}). Specifically, there is a version of the theorem for so-called \emph{good} groups, whose definition we will not go into (see instead~\cite{Freedman-Teichner:1995-1,DET-book-goodgroups}). For virtually all applications, it suffices to know that the class of good groups contains groups of subexponential growth~\cite{Freedman-Teichner:1995-1, Krushkal-Quinn:2000-1}, and is closed under subgroups, quotients, extensions, and colimits~\cite[{p.\ 44}]{FQ}. In particular, all finite groups and all solvable groups are good. It is not known whether non-abelian free groups are good. 

    In the case of non-trivial fundamental groups, the self-intersection number of~$g$ and the intersection between $f$ and $g$ is no longer just the signed count of intersections, but rather an \emph{equivariant} version, with values lying in (a quotient of)~$\Z[\pi_1(M)]$. We will define these in \Cref{sec:lambda-mu}.

    We will use \Cref{thm:DET} in \Cref{sec:direct}, and the version for a finite collection of discs in a $4$-manifold with infinite cyclic fundamental group in \Cref{sec:alex-poly-one}. 
\end{remark}

\begin{remark}
    The disc embedding theorem is the key ingredient in the proof of the $4$-dimensional topological $s$-cobordism theorem for good fundamental groups~\cite[{Theorem~7.1A}]{FQ}. The disc embedding theorem also implies the \emph{sphere embedding theorem} (\Cref{thm:sphere-emb-thm}), which is the key ingredient in proving the exactness of the topological surgery sequence in dimension four for good fundamental groups~\cite[{Theorem~11.3A}]{FQ} (see also~\cite{DET-book-surgery} and \Cref{sec:surgery-sequence}). These are powerful tools that are central, for example, in proving classification results for topological $4$-manifolds up to homeomorphism.
\end{remark}

\begin{remark}
     Historically, the first version of the disc embedding theorem was proven by Freedman for a finite collection of discs in an arbitrary smooth, simply connected $4$-manifold. This was the ingredient needed by Quinn in~\cite{quinn:endsIII} to prove many fundamental results, such as those mentioned in \Cref{sec:defs-fund-tools}. Using these tools, Freedman's proof could be repeated, but now in a topological ambient $4$-manifold. The techniques of the proof were also further developed by Freedman and Quinn to now apply to ambient $4$-manifolds with good fundamental group. This was the proof given in~\cite{FQ} and then explained further in~\cite{DETbook}.
\end{remark}

\section{Representing primitive homology classes by locally flat tori}\label{sec:direct}

In this section we prove Theorem A, and then briefly state the more general results of \cite{LW97,KPRT:sigmet}. 

\subsection{Proof of Theorem A}\label{sec:proof-primitive-torus}
We recall the statement of Theorem A for the convenience of the reader. 

\begin{theoremA}
    Let $M$ be a closed, simply connected $4$-manifold. Then every primitive class in $H_2(M;\Z)$ is represented by a locally flat embedded torus. 
\end{theoremA}
\begin{proof}
Let $\alpha\in H_2(M;\Z)$ be a primitive class. We split up the proof in a number of steps. 

\begin{step}\label{step:generic-immersions}
    Represent $\alpha$ by a generic immersion $f\colon S^2\looparrowright M$ with a geometrically dual sphere $g\colon S^2\looparrowright M$, i.e. $f$ and $g$ intersect each other transversely, and only at a single point. 
\end{step}

First we use that $M$ is simply connected. By the Hurewicz theorem~${\pi_2(M)\cong H_2(M;\Z)}$, so every class in $H_2(M;\Z)$ can be represented by a map of a sphere with a given orientation. Then by the immersion lemma (\Cref{thm:immersion-lemma}) we can assume further that this map is a generic immersion. Since $M$ is simply connected, it is orientable. Fix an orientation on $M$.

By Poincar\'{e} duality, we know that the intersection form of $M$ is unimodular. Therefore, since $\alpha$ is a primitive class in $H_2(M;\Z)$, it has a dual class. In other words, there is some $\beta\in H_2(M;\Z)$ such that $\alpha\cdot \beta=1$. Again by the immersion lemma (\Cref{thm:immersion-lemma}), the class $\beta$ can be represented by a generic immersion~$g\colon S^2\looparrowright M$, along with an orientation on $g$, such that $f\cdot g=1$, where this is both the homological intersection number and the signed count of intersections between $f$ and $g$. Here note that we need orientations on $f$ and $g$, as well as the orientation on $M$, to precisely talk about the signs of the intersection points, and to determine the intersection form on $M$. Next, by topological transversality (\Cref{thm:normal-bundles-transversality}\,\eqref{item:top-transversality}) we can assume, after an isotopy, that $f$ and $g$ intersect transversely. Note that this is not a direct application of the theorem, since $f$ and $g$ are not embeddings. But since $f$ and $g$ are generic immersions, we can restrict to subsets of the domain where the restrictions are in fact embeddings, and apply the theorem there. A more technical version of the topological transversality theorem~\cite[{Theorem~9.5A}]{FQ} then allows us to patch those local isotopies together. A final step ensures that double points of $f$ and of $g$ do not coincide with double points between $f$ and $g$. 

At this point, the spheres $f$ and $g$ are algebraically dual, but not necessarily geometrically dual. To arrange for them to be geometrically dual, we will use the \emph{geometric Casson lemma}, which we leave to the reader as an advanced exercise (Exercise~\ref{ex:geometric-Casson-lemma}).\footnote{Solving the exercise will be an excellent method for field-testing the various geometric manoeuvres described in the current proof.} This lemma says that we can perform a regular homotopy to remove a pair of algebraically cancelling intersections between $f$ and $g$, at the cost of more self-intersections of $f$ or of $g$. This is not a significant price for us, since we have no control on the self-intersections of $f$ and $g$ at this stage anyway. By repeated applications of the lemma we arrange that $f$ and $g$ are geometrically dual as desired. If we were being very precise, we would use new notation for the maps produced by applying the lemma. However, as is customary, we will keep using the  original symbols $f$ and $g$. 

\begin{step}\label{step:signed-count-zero}
    Arrange that the signed count of self-intersections of $f$ is zero. 
\end{step}

We will use \emph{interior twisting}. This procedure is best described pictorially (see \Cref{fig:interior-twisting}). 
\begin{figure}[tb]
\centering
    \begin{tikzpicture}
        \node[anchor=south west,inner sep=0] at (0,0){\includegraphics[width=6cm]{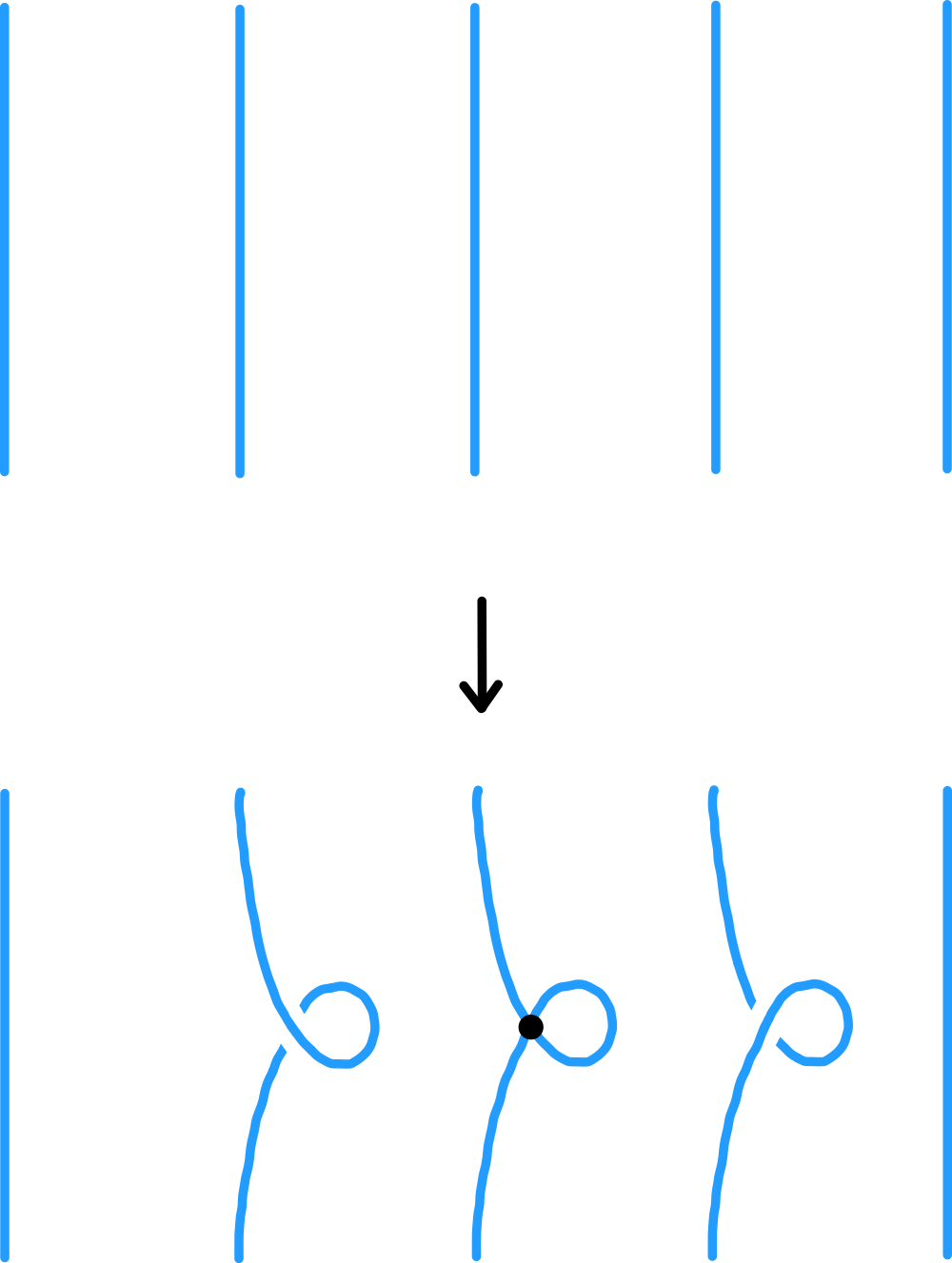}};
           \node at (0,-0.25) {$t=-\varepsilon$};
            \node at (1.5,-0.25) {$t=-\tfrac{\varepsilon}{2}$};
            \node at (3,-0.25) {$t=0$};
            \node at (4.5,-0.25) {$t=\tfrac{\varepsilon}{2}$};
            \node at (6,-0.25) {$t=\varepsilon$};
            \node at (0,4.75) {$t=-\varepsilon$};
            \node at (1.5,4.75) {$t=-\tfrac{\varepsilon}{2}$};
            \node at (3,4.75) {$t=0$};
            \node at (4.5,4.75) {$t=\tfrac{\varepsilon}{2}$};
            \node at (6,4.75) {$t=\varepsilon$};
	\end{tikzpicture} 
    \caption{Interior twisting. A small patch of a generically immersed surface $f$ in an ambient $4$-manifold is shown in blue on the top. Note the patch has no self-intersections of $f$. The procedure of interior twisting replaces the patch on the top with the patch on the bottom. Note the patch on the bottom contains a transverse self-intersection.}
    \label{fig:interior-twisting}
\end{figure}
Since this is our first explicit geometric construction, let us take a moment to describe it properly. In the figure on the top, we have a described a small patch on $f$ in a movie picture. In other words, the blue vertical lines can be stitched together to give a small patch on $f$, specifically a region with no double points. The boundary of the patch consists of the leftmost and rightmost time slices, as well as the boundaries of the intermediate time slices. As expected, these pieces glue together to give a rectangle on the boundary of the patch. 

The figure on the bottom describes a modified patch. Notice that the original and the modified patches agree on their boundaries, so you could imagine taking out the small patch on $f$ shown on top, and gluing in the surface on the bottom, like a band-aid.\footnote{A better analogy would be that we take out one piece of a jigsaw puzzle and replace it with another one, which of course is only allowed if the boundaries are identical.} The procedure of replacing an arbitrary patch on $f$ by this band-aid, or its mirror image, is called \emph{interior twisting}. 

The key property of the band-aid is that it contains a transverse double point singularity in the middle time slice. Using the mirror image of the patch results in a double point singularity of the opposite sign. Therefore, by enough interior twisting of the appropriate sign we can arrange that the signed count of self-intersections of~$f$ is zero. To be precise, the procedure of interior twisting changes $f$ by a homotopy, but the result is still a generically immersed sphere, which we continue to refer to as $f$. By doing the procedure away from $g$, we can assume that $f$ and $g$ remain geometrically dual. 

What is the sign of the intersection point created in \Cref{fig:interior-twisting}? We leave this as an exercise for the motivated reader (Exercise~\ref{ex:interior-twisting-sign}).

\begin{step}\label{step:pair-gen-whitney}
    Pair up the points in $f\pitchfork f$ by generically immersed Whitney discs. 
\end{step}

Since the signed count of self-intersections of $f$ is trivial, we can arbitrarily pair up points with opposite sign. For each such pair, the two constituent points can be joined by two arcs, one on each sheet.  
The union of these arcs is a circle in $M$. Recalling again that $\pi_1(M)=1$, we note that each such circle is null homotopic in~$M$. Applying \Cref{thm:immersion-lemma} and \Cref{thm:normal-bundles-transversality}\,\eqref{item:top-transversality}, we can assume that these circles bound a collection of generically immersed discs $\{W_i\}$, which intersect one another,~$f$, and~$g$ only transversely and only in the interiors, except along the boundary circles. 

\medskip
Let us take a brief hiatus from the proof to describe what we would like to be true for this collection $\{W_i\}$. In the ideal situation, we would be able to do the Whitney move on $f$ along $\{W_i\}$, resulting in an embedding. This would complete the proof of Theorem A, in fact producing an embedded sphere rather than a torus as claimed. For the Whitney move to produce an embedding, we would need the Whitney discs to be locally flat embedded, pairwise disjoint, have interiors disjoint from $f$, and induce the correct framing on the boundary. However, at present, we can guarantee none of these features. In other words, \textit{a priori} we have four distinct families of obstructions to being able to do the Whitney move on $\{W_i\}$: the intersections amon the interiors~$\{\mathring{W}_i\}$, including self-intersections; the intersections among~$\{\partial W_i\}$, including self-intersections; intersections between $\{\mathring{W}_i\}$ and $f$; and finally, for each $i$, the difference, denoted by $\tw(\partial W_i)\in \Z$, between the Whitney framing on $\partial W_i$ and the framing induced by the normal vector bundle of $W_i$. (For a few more details about the twisting numbers $\tw(\partial W_i)$ see \Cref{sec:finger-whitney}.) We summarise these obstructions in \Cref{tab:table}. While at first glance they may seem independent of one another, in fact they are related. Moreover, we have a toolbox of geometric manoeuvres, which allows us to trade problems of one sort for those of a different sort in a precise way, as indicated in the table. We note that most of the manoeuvres have an associated cost, so we cannot simply assume away all the obstructions. But we can still apply these moves cleverly and in the right order and hope for the best. We will see that in many (but not all) situations we can in fact assume that all the obstructions vanish (see Exercise~\ref{ex:spin-primitive-sphere}). 

\begin{table}[tb]
\caption{Problems, their solutions, and associated costs}
\label{tab:table}
\centering
\begin{tabular}{cp{20mm}p{45mm}p{55mm}}\toprule
Type    &Problem &Solution   &Cost\\
\midrule
\addlinespace[5pt]
1   &$\mathring{W}_i\pitchfork \mathring{W}_j$   &Disc embedding theorem &\begin{tabular}[c]{@{}l@{}}None if all else solved 
(\Cref{prop:whitney-by-DET})
\end{tabular}\\
\midrule
\addlinespace[5pt]
\multirow{4}{*}{2}  &\multirow{4}{*}{$\tw(\partial W_i)$} &Interior twisting &\begin{tabular}[c]{@{}l@{}}$\tw(\partial W_i)\mapsto \tw(\partial W_i)\pm 2$\\
\addlinespace[5pt]
$\left\lvert \mathring{W}_i\pitchfork \mathring{W}_i\right\rvert  \mapsto \left\lvert \mathring{W}_i\pitchfork \mathring{W}_i\right\rvert +1$
\end{tabular}\\
\cmidrule{3-4}
\addlinespace[5pt]
    &    &Boundary twisting   &\begin{tabular}[c]{@{}l@{}}$\tw(\partial W_i)\mapsto \tw(\partial W_i)\pm 1$\\
    \addlinespace[5pt]
    $\left\lvert \mathring{W}_i\pitchfork f\right\rvert  \mapsto \left\lvert \mathring{W}_i\pitchfork f\right\rvert +1$
\end{tabular}\\
\midrule
\addlinespace[5pt]
3   &$\partial W_i\pitchfork \partial W_j$   & Boundary pushoff  &\begin{tabular}[c]{@{}l@{}}$\left\lvert \partial W_i\pitchfork \partial W_j\right\rvert\mapsto \left\lvert\partial W_i\pitchfork \partial W_j\right\rvert-1$\\
\addlinespace[5pt]
$\left\lvert\mathring{W}_i\pitchfork f\right\rvert \mapsto \left\lvert\mathring{W}_i\pitchfork f\right\rvert +1$
\end{tabular}\\
\midrule
\addlinespace[5pt]
\multirow{4}{*}{4}  &\multirow{4}{*}{$\mathring{W}_i\pitchfork f$ } &Tubing into $g$ &\begin{tabular}[c]{@{}l@{}}$\left\lvert\mathring{W}_i\pitchfork f\right\rvert \mapsto \left\lvert\mathring{W}_i\pitchfork f\right\rvert -1$\\
\addlinespace[5pt]
$\tw(\partial W_i)\mapsto \tw(\partial W_i)+e(\nu g)$\\
\addlinespace[5pt]
$\left\lvert\mathring{W}_i \pitchfork \mathring{W}_j\right\rvert$ uncontrolled
\end{tabular}\\
\cmidrule{3-4}
\addlinespace[5pt]
    &    &Transfer move   &\begin{tabular}[c]{@{}l@{}}
    $\left\lvert\mathring{W}_i\pitchfork f\right\rvert \mapsto \left\lvert\mathring{W}_i\pitchfork f\right\rvert+1$\\
    \addlinespace[5pt]
    $\left\lvert\mathring{W}_j\pitchfork f\right\rvert \mapsto \textit{}\left\lvert\mathring{W}_j\pitchfork f\right\rvert+1$
\end{tabular}\\\addlinespace[5pt]
\bottomrule
\end{tabular}
\end{table}

Let us now work through the techniques mentioned in \Cref{tab:table}. First, we justify our statement in the table that problems of type 1 can be solved at no cost, if all other problems have also been solved, by applying the disc embedding theorem. 

\begin{proposition}\label{prop:whitney-by-DET}
    Let $\Sigma$ be a surface and let $M$ be a $4$-manifold. Let $f\colon \Sigma \looparrowright M$ be a generic immersion, such that all the double points of $f$ are paired up by generically immersed Whitney discs $\{W_i\}$. Assume further that $f$ has a geometrically dual sphere~$g$. Suppose that $\tw(\partial W_i)=0$ and $\partial W_i\pitchfork \partial W_j=\emptyset=\mathring{W}_i\pitchfork f$ for all $i,j$. Then there exists $\{\ol{W}_i\}$, a collection of locally flat embedded and disjoint Whitney discs pairing all the points in $f\pitchfork f$, with trivial twisting numbers, and with interiors disjoint from $f$.
\end{proposition}

\begin{proof}
    We will apply the disc embedding theorem (\Cref{thm:DET}) to $\{W_i\}$ in~${N:=M\sm \mathring{\nu}f}$, where $\mathring{\nu}f$ is an open tubular neighbourhood of $f$, which exists by \Cref{thm:plumbed-normal-bundle}. To be precise, we need the version for a finite collection of discs; see \Cref{thm:DET-full}. We have to check that the hypotheses hold. First we need that $\pi_1(N)=1$. This follows from Exercise~\ref{ex:pi1-negligible}. We also need algebraically dual spheres. For each $W_i$, let $T_i$ denote the Clifford torus at one of the two double points of $f$ paired by $W_i$. As we see in \Cref{fig:clifford-torus-to=alg-dual}, each $T_i$ lies in $N$ and is geometrically dual to $W_i$. Furthermore it satisfies~$T_i\pitchfork W_j=\emptyset$ if $i\neq j$. We will modify each $T_i$ to a sphere. Note that a meridional disc for $T_i$ intersects $f$ at a single point. Tube the meridional disc to $g$, to get a disc bounded by a meridian of $T_i$ lying entirely in $N$. Compressing $T_i$ along two copies of this meridional disc, using the framing induced by $T_i$, produces a sphere $S_i$ with trivial normal vector bundle. We need to check that this collection of spheres satisfies~$S_i\cdot W_j=\delta_{ij}$ and $S_i\cdot S_j=0$ for all $i,j$ -- both follow from the fact that each compression was along two copies of a fixed meridional disc, with opposite orientations. This shows that the hypotheses of the disc embedding theorem are satisfied for $\{W_i\}$ and $\{S_i\}$ in $N$. Therefore, the theorem provides the desired embedded and disjoint Whitney discs, with trivial twisting number, and with interiors disjoint from~$f$. 
    \end{proof}

\begin{figure}[tbp]
    \centering
    \begin{tikzpicture}
        \node[anchor=south,inner sep=0] at (5.45,6){\includegraphics[width=10.9cm]{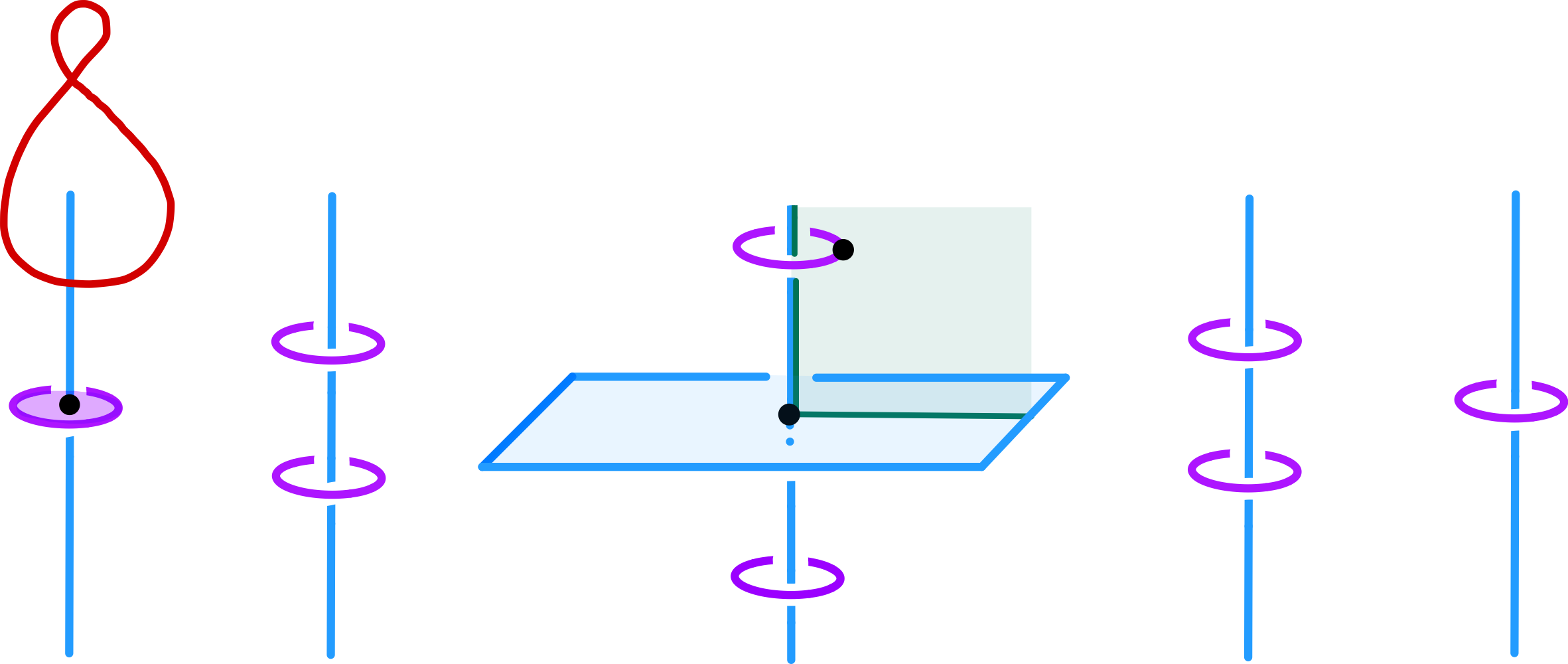}};
        \node[anchor=south,inner sep=0] at (5.4,0){\includegraphics[width=11cm]{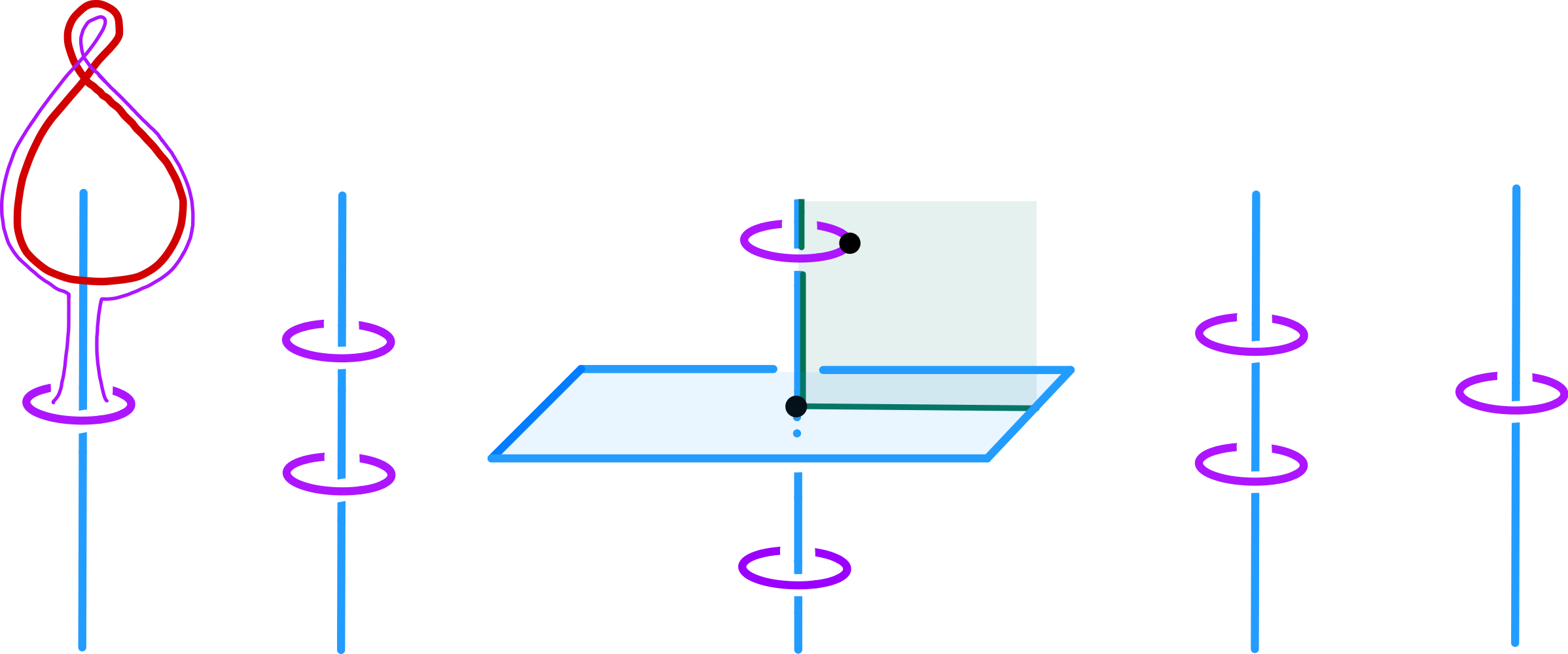}};
            \node at (5.5,5) {(a)};
            \node at (0.5,5.7) {$t=-\varepsilon$};
            \node at (2.3,5.7) {$t=-\tfrac{\varepsilon}{2}$};
            \node at (5.5,5.7) {$t=0$};
            \node at (8.7,5.7) {$t=\tfrac{\varepsilon}{2}$};
            \node at (10.55,5.7) {$t=\varepsilon$};
            \node at (-0.2,9) {$g$};       
            \node at (4,7.1) {$f$};  
            \node at (0.3,6.4) {$f$};
            \node at (5,9) {$T$};
            \node at (6.5,8.5) {$W$};
            \node at (5.5,-1) {(b)};
            \node at (0.5,-0.25) {$t=-\varepsilon$};
            \node at (2.3,-0.25) {$t=-\tfrac{\varepsilon}{2}$};
            \node at (5.5,-0.25) {$t=0$};
            \node at (8.7,-0.25) {$t=\tfrac{\varepsilon}{2}$};
            \node at (10.55,-0.25) {$t=\varepsilon$};
            \node at (-0.2,3) {$g$};       
            \node at (4,1.1) {$f$};  
            \node at (0.3,0.4) {$f$};
            \node at (5,3) {$T$};
            \node at (6.5,2.5) {$W$};
	\end{tikzpicture} 
    \caption{Using the Clifford torus to produce an algebraically dual sphere. (a) The surface $f$ is shown in blue, in a small coordinate patch close to a transverse self-intersection, shown in black. In other words, we see two sheets of $f$, one showing up as the flat sheet in the central time slice, and the other smeared out across the time direction, showing up as a single line in each time slice. The Whitney disc $W$ at this point of intersection is shown in green. The geometrically dual sphere $g$ is shown on the leftmost time slice -- as required it intersects $f$ at a single point. The Clifford torus $T$ is shown in purple. Note that $T$ intersects $W$ precisely once, in the central time slice. A meridional disc for $T$ is shaded in purple in the leftmost time slice. It intersects $f$ at a single point. (b) We show how tubing into a parallel copy of $g$ produces a meridional disc for $T$ which is disjoint from~$f$.}
    \label{fig:clifford-torus-to=alg-dual}
\end{figure}

\Cref{prop:whitney-by-DET} shows that if we can solve all the problems of type 2, 3, and 4, then the problems of type 1 can also be solved. Then we can do the Whitney move on $f$ along the resulting Whitney discs to produce a locally flat embedded sphere which is homotopic to $f$. Note that we can do the Whitney move along locally flat discs, since they have normal vector bundles, by \Cref{thm:normal-bundles-transversality}\,\eqref{item:normal-bundles}.

\begin{figure}[tb]
    \centering
    \begin{tikzpicture}
        \node[anchor=south west,inner sep=0] at (0,0){\includegraphics[width=11cm]{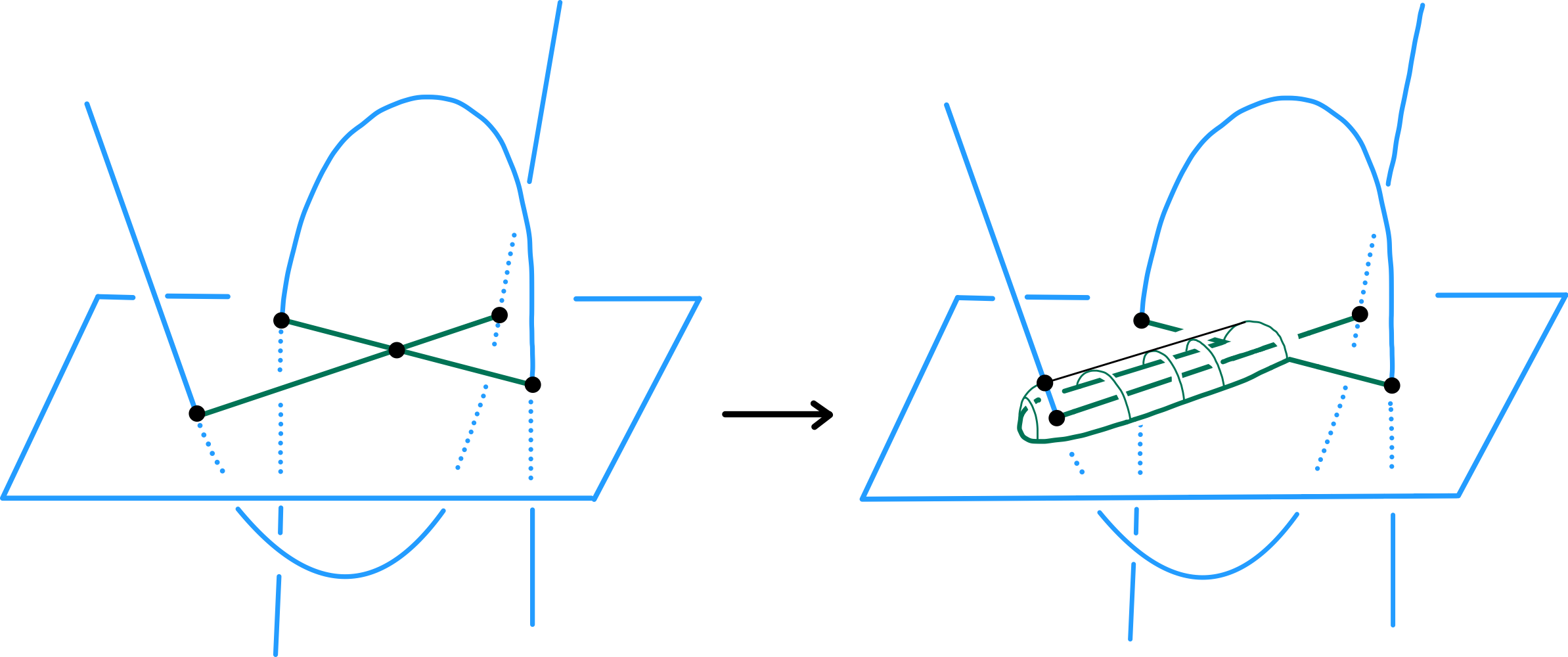}};
            \node at (1.75,0.2) {$f$};
            \node at (0.4,3.7) {$f$};               \node at (0.5,0.8) {$f$};
            \node at (7.7,0.2) {$f$};
            \node at (6.4,3.7) {$f$};               \node at (6.5,0.8) {$f$};
	\end{tikzpicture} 
    \caption{Boundary pushoff. The surface $f$ is shown in blue, along with two pairs of algebraically cancelling intersection points paired by Whitney discs. The Whitney discs are not shown, but a Whitney arc for each is shown in green. Note that they intersect in a single point. The procedure of boundary pushoff moves one of the Whitney arcs off the other, by pushing towards the boundary of the arc, as shown on the right. Note that this procedure creates an intersection between $f$ and a Whitney disc.}
    \label{fig:boundary-pushoff}
\end{figure}

Next we describe the geometric manoeuvres mentioned in \Cref{tab:table}. We already saw interior twisting in \Cref{step:signed-count-zero}. We also have the operation of \emph{boundary twisting}, described in \Cref{fig:boundary-twisting}. Determining the effect of interior and boundary twisting on the various problems in \Cref{tab:table} comprises Exercise~\ref{ex:twisting-euler}. The reader might wonder why we need two solutions to problems of type 2. So we remark that interior twisting is \textit{a priori} less effective than boundary twisting, since it can only change the twisting number by even numbers, rather than arbitrary integers. But interior twisting is cheap -- it only creates problems of type 1, which can be solved  `for free' by \Cref{prop:whitney-by-DET}. In contrast, boundary twisting is much more expensive -- it creates problems of type 4, which are in general much harder to fix. For example, solving a problem of type 4 by tubing into $g$ creates problems of type 2, which are what we were trying to solve in the first place. So with boundary twisting one is in danger of getting stuck in a loop of circular reasoning.

\begin{figure}[tb]
    \centering
    \begin{tikzpicture}
        \node[anchor=south,inner sep=0] at (5.5,7){\includegraphics[width=5cm]{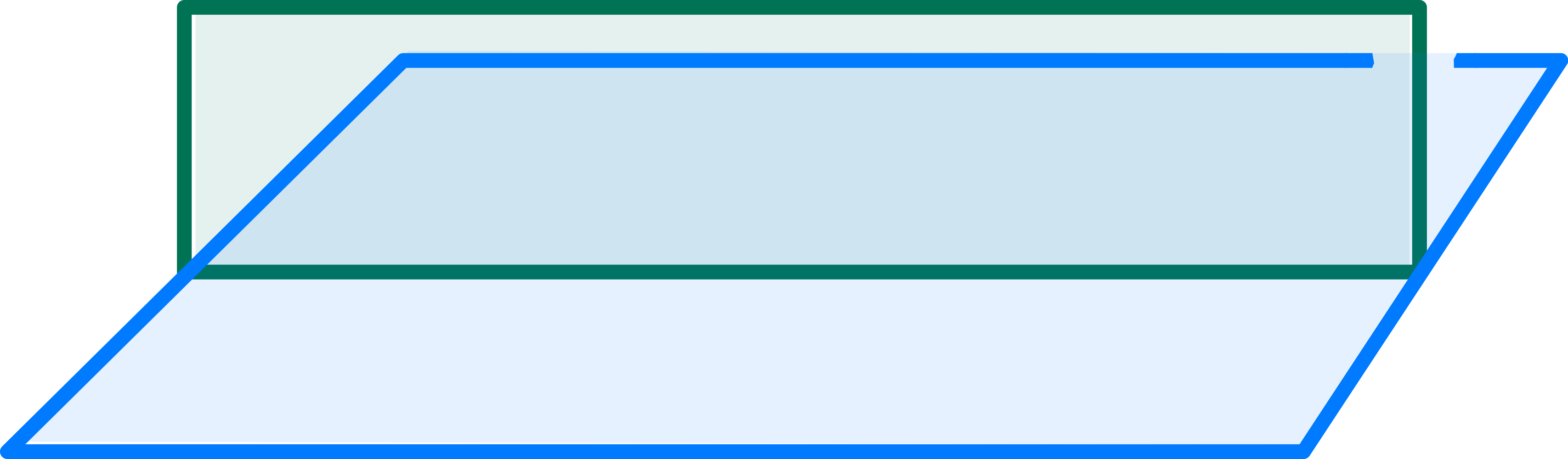}};
        \node[anchor=south west,inner sep=0] at (0,0){\includegraphics[width=11cm]{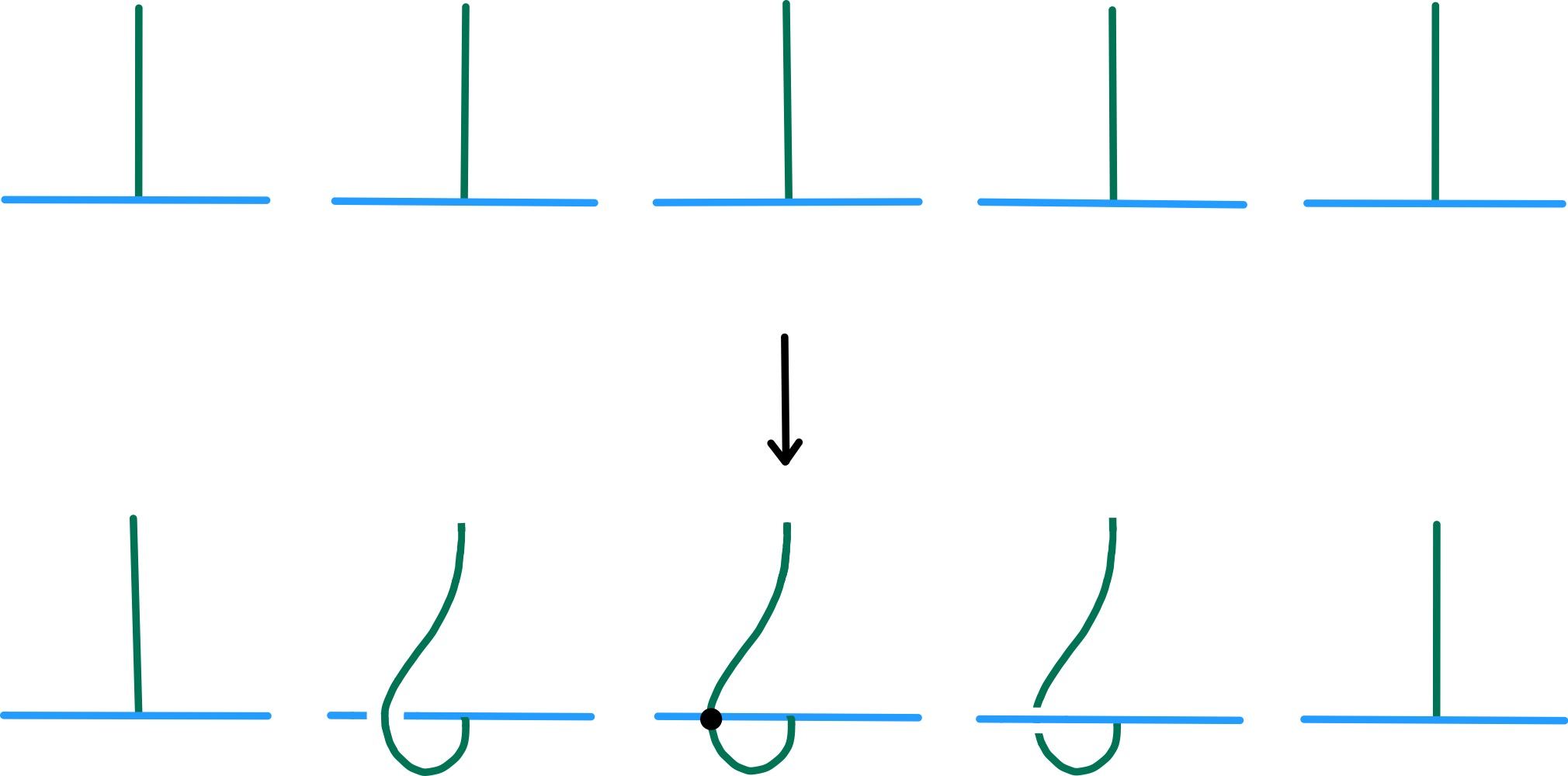}};
        \node at (5.5,6.5) {(a)};     
        \node at (3.5,6.75) {$f$};
        \node at (5.5,8.75) {$W$};
        \node at (5.5,-0.25) {(b)};  
        \node at (1.8,0.7) {$f$};
        \node at (4,0.7) {$f$};
        \node at (6.2,0.7) {$f$};
        \node at (8.6,0.7) {$f$};
        \node at (10.7,0.7) {$f$};
        \node at (1.2,1.6) {$W$};
        \node at (3.5,1.6) {$W$};
        \node at (5.8,1.6) {$W$};
        \node at (8.1,1.6) {$W$};
        \node at (10.3,1.6) {$W$};
        \node at (1.8,4.3) {$f$};
        \node at (4,4.3) {$f$};
        \node at (6.2,4.3) {$f$};
        \node at (8.6,4.3) {$f$};
        \node at (10.7,4.3) {$f$};
        \node at (1.2,5) {$W$};
        \node at (3.5,5) {$W$};
        \node at (5.8,5) {$W$};
        \node at (8.1,5) {$W$};
        \node at (10.3,5) {$W$};
	\end{tikzpicture} 
    \caption{Boundary twisting. (a) The procedure will take place in a small neighbourhood of a Whitney arc, as depicted here. A generically immersed surface $f$ in an ambient $4$-manifold is shown in blue. A small portion of a Whitney disc $W$ with boundary on $f$ is shown in green. (b) The region depicted in (a) is now split up into multiple time slices in a coordinate neighbourhood in the ambient manifold. The procedure of boundary twisting involves changing the local picture above to the local picture below, by twisting a boundary collar of $W$ around $f$, as shown. Note that this procedure creates a new intersection between $W$ and $f$.}
    \label{fig:boundary-twisting}
\end{figure}

We also have the boundary pushoff operation shown in \Cref{fig:boundary-pushoff}. The reader might rightly complain that we could have chosen the Whitney arcs originally so that they do not intersect. However, we include Whitney arc obstructions in our list in \Cref{tab:table} since some upcoming geometric constructions will create them, so it will be useful to know how to solve them and at what price. 

The next operation in \Cref{tab:table} is to tube intersections of some $\mathring{W}_i$ with $f$ into the geometric dual $g$. We already saw this operation in the proof of \Cref{prop:whitney-by-DET}, but we give a few more details here. Suppose we have a generically immersed connected surface $A$. Let $B$ and $B'$ be two other generically immersed surfaces intersecting $A$ transversely at points $p$ and $p'$ respectively with $p\neq p'$. Let $C$ be an embedded arc in $A$ joining $p$ and $p'$, and not passing through any double points of $A$. The normal vector bundle of $A$ restricted to $C$ is trivial (since $C$ is contractible). In other words, there is a copy of $C\times D^2$ which intersects $B$ and $B'$ in small discs about $p$ and $p'$ respectively, and only intersects $A$ along $C$. Cut out these discs from $B\cup B'$ and glue on the rest of the boundary $\partial (C\times D^2$) to $B\cup B'$ minus the discs. In other words, we are gluing in a meridional annulus for $C$. This process is called \emph{tubing $B$ into $B'$ (along $A$ (or $C$))}, and is described in \Cref{rfig:tubing}. 
\begin{figure}[tb]
    \centering
    \begin{tikzpicture}
        \node[anchor=south west,inner sep=0] at (0,0){\includegraphics[width=11cm]{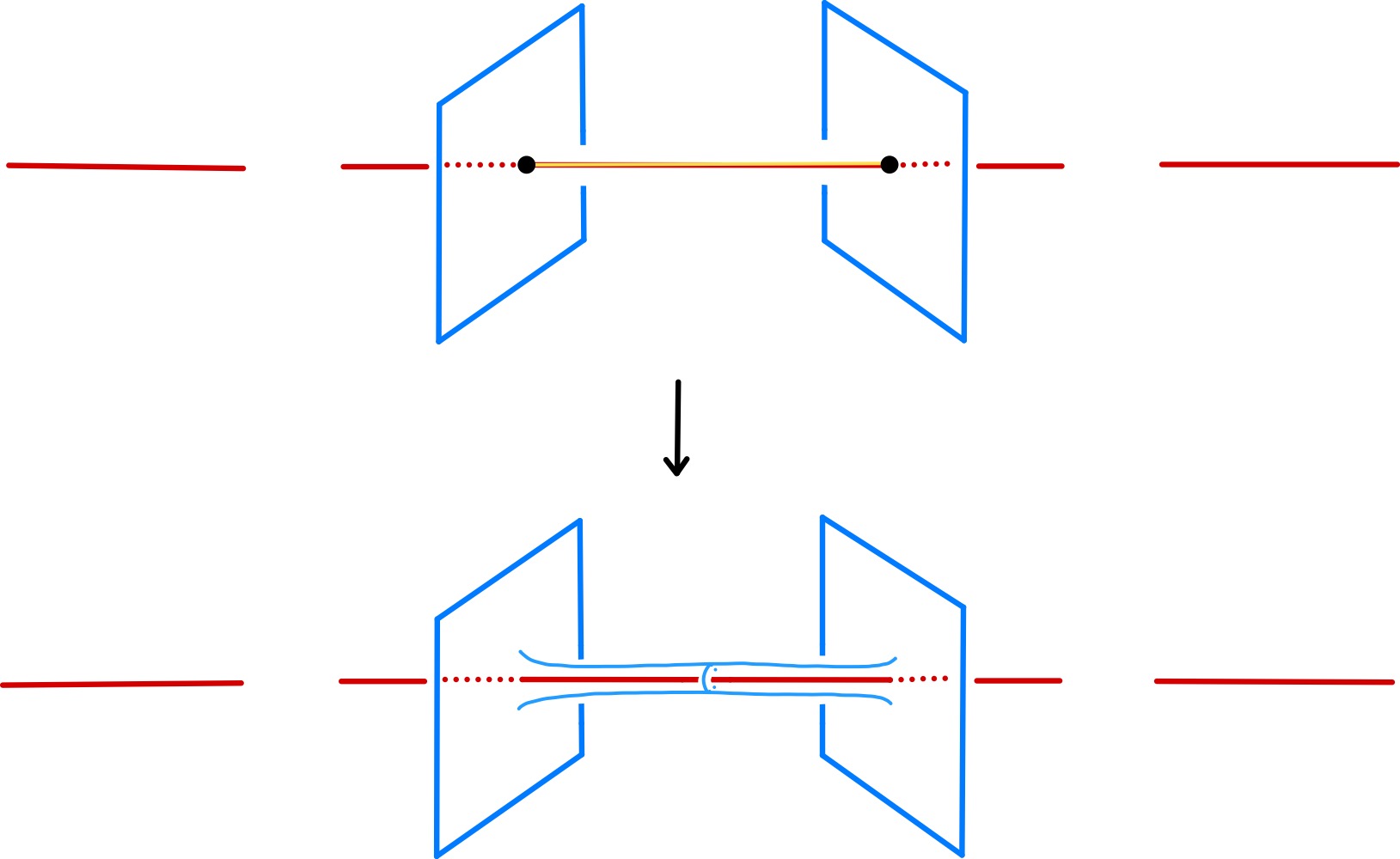}};
        \node at (1,-0.25) {$t=-\varepsilon$};
        \node at (5.5,-0.25) {$t=0$}; 
        \node at (10,-0.25) {$t=\varepsilon$};
        \node at (0.15,5.75) {$A$};
        \node at (2.9,5.75) {$A$};
        \node at (9.25,5.75) {$A$};
        \node at (3.2,4.2) {$B$};
        \node at (7.8,4.2) {$B'$};
        \node at (4.1,5.2) {$p$};
        \node at (7,5.2) {$p'$};        
        \node at (1,4.2) {$t=-\varepsilon$};
        \node at (5.5,4.2) {$t=0$}; 
        \node at (10,4.2) {$t=\varepsilon$};
        \node at (5.5,5.75) {$C$};
        \node at (0.15,1.6) {$A$};
        \node at (2.9,1.6) {$A$};
        \node at (9.25,1.6) {$A$};
        \node at (3.2,0.2) {$B$};
        \node at (7.8,0.2) {$B'$};        
	\end{tikzpicture} 
    \caption{Tubing. Top: two surfaces $B$ and $B'$, both shown in blue, intersect a third surface $A$, shown in red, at points $p$ and $p'$ respectively. An arc $C$ on $A$ joining $p$ and $p'$ is shown in yellow. Bottom: the result of tubing $B$ to $B'$ along $A$ is shown.}
    \label{rfig:tubing}
\end{figure}
Usually we do not tube into $B'$ but rather a pushoff thereof -- this allows us to tube multiple times. For instance if~$B'$ is embedded, geometrically dual to $A$, and has trivial normal vector bundle, then all intersections of some $B$ with $A$ can be removed by tubing into (distinct) pushoffs of~$B'$, without creating any additional intersections. Note that when a surface $B$ is tubed into a surface $B'$, the euler number of the normal vector bundle of the result is the sum of the two euler numbers of the original $B$ and $B'$. Similarly, if a Whitney disc is tubed into a generically immersed sphere $g$, the twisting number of the result is the sum of the original twisting number with the euler number of the normal vector bundle of $g$.

The final manoeuvre in \Cref{tab:table} is the \emph{transfer move}. This is described in \Cref{fig:transfer-move}. 
\begin{figure}[tb]
    \centering
    \begin{tikzpicture}
        \node[anchor=south west,inner sep=0] at (0,9.5){	\includegraphics[width=7cm]{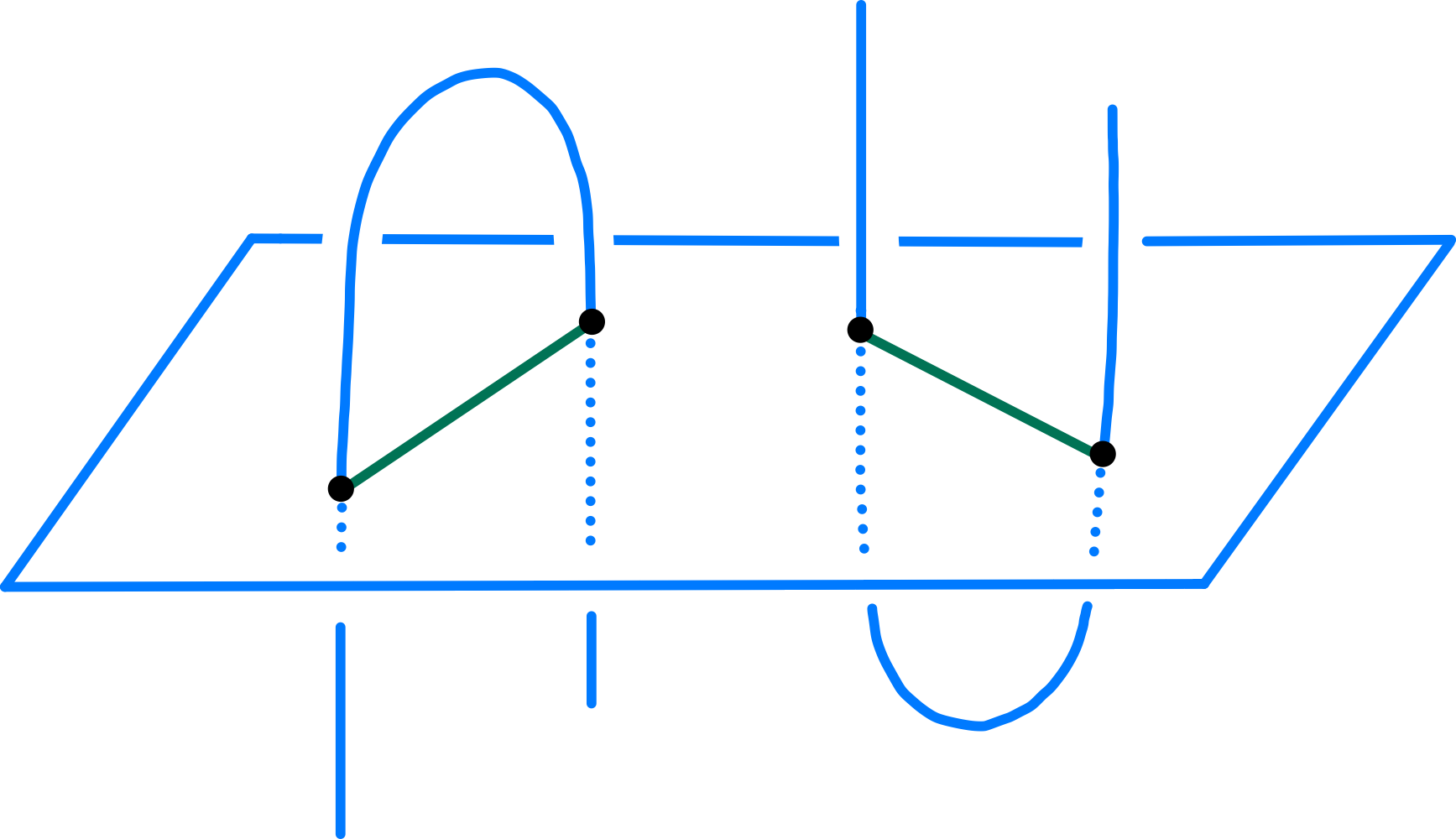}};
        \node[anchor=south west,inner sep=0] at (0,4.75){	\includegraphics[width=7cm]{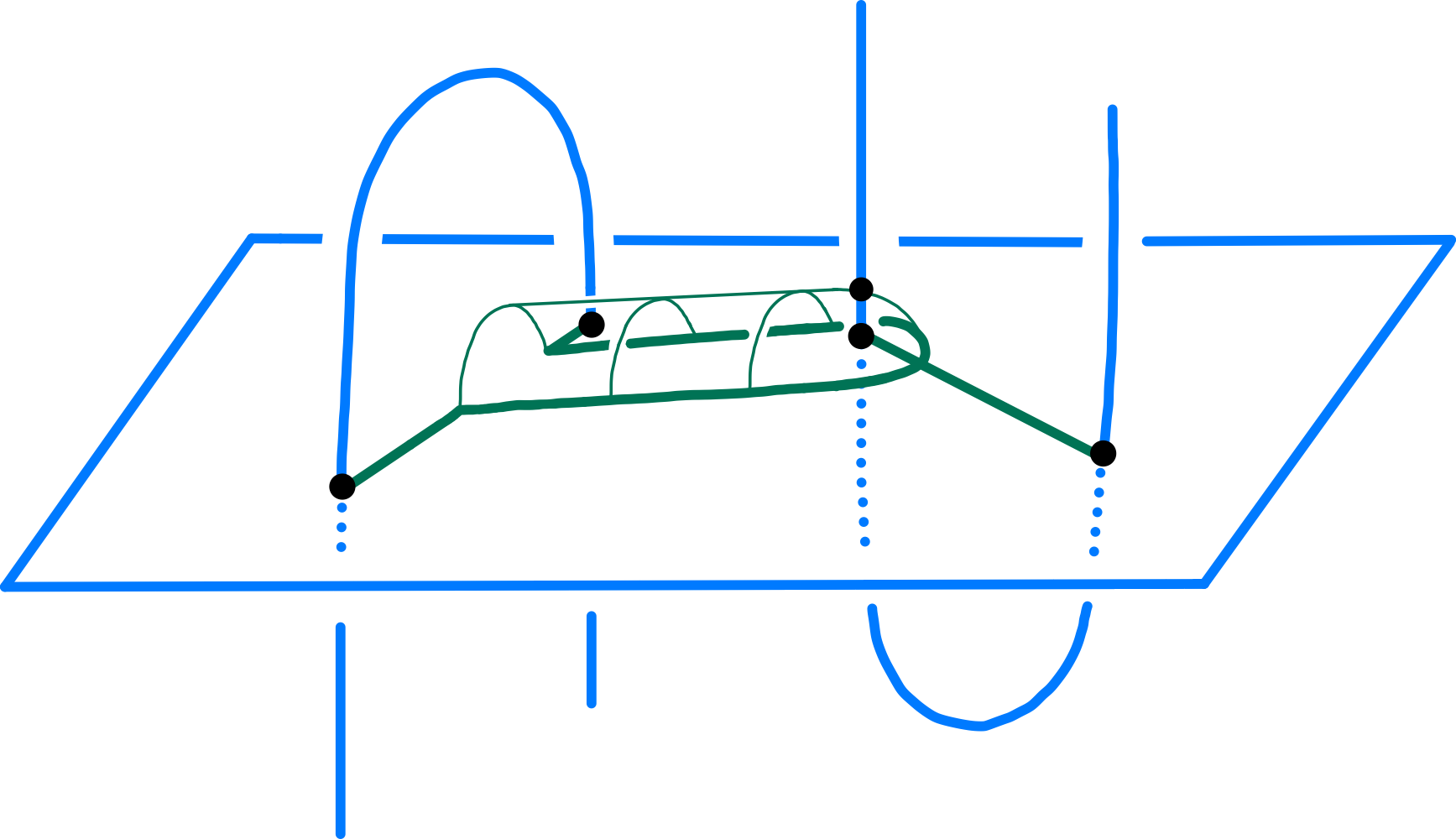}};
        \node[anchor=south west,inner sep=0] at (0,0){	\includegraphics[width=7cm]{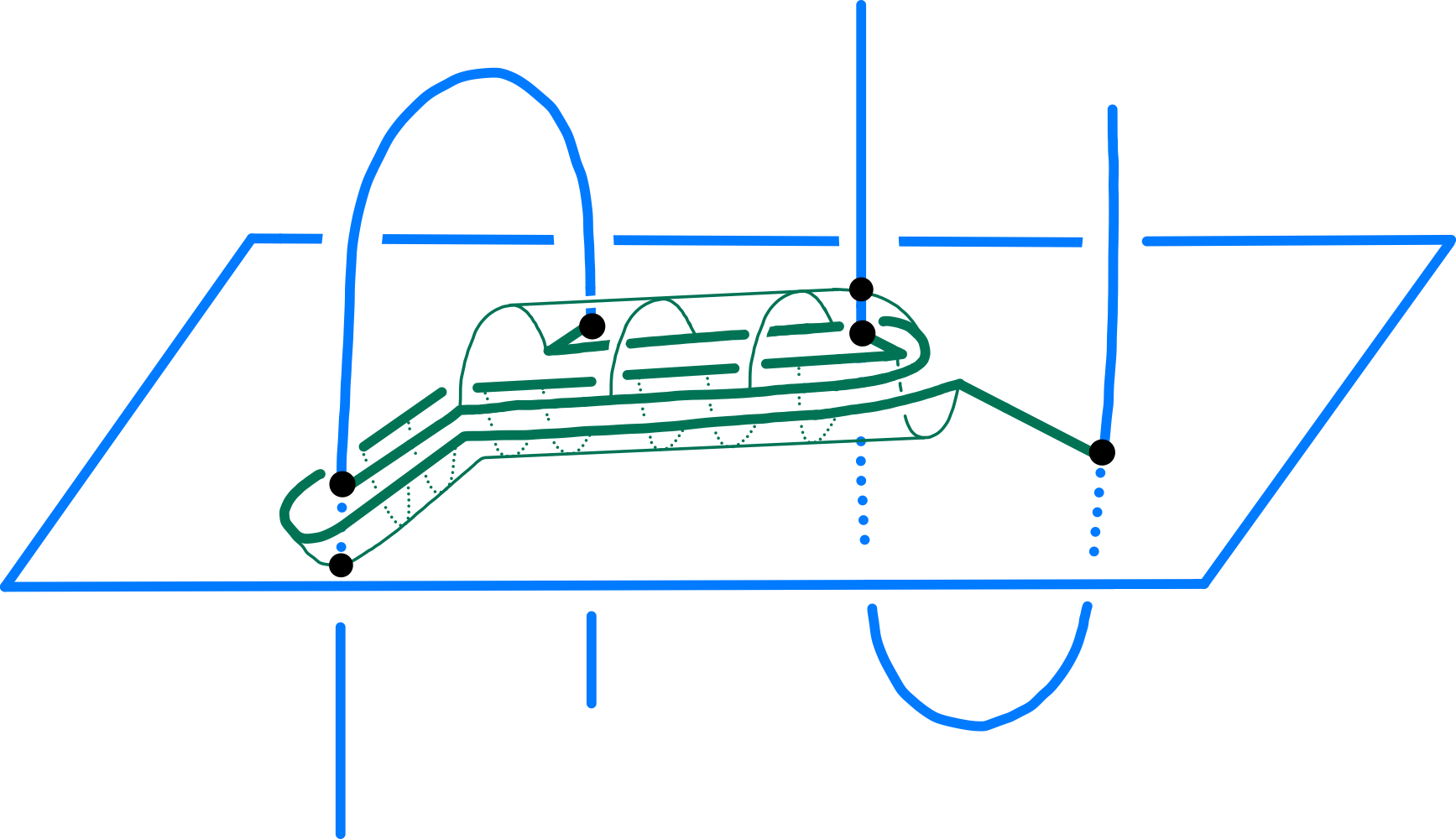}};        
        \node at (3.5,9.25) {(a)};
        \node at (0.5,10.4) {$f$};
        \node at (1.4,9.75) {$f$};
        \node at (3.9,13.25) {$f$};
        \node at (3.5,4.5) {(b)};
        \node at (0.5,5.65) {$f$};
        \node at (1.4,5) {$f$};
        \node at (3.9,8.5) {$f$};
        \node at (3.5,-0.25) {(c)};
        \node at (0.5,0.9) {$f$};
        \node at (1.4,0.25) {$f$};
        \node at (3.9,3.75) {$f$};
	\end{tikzpicture} 
    \caption{Transfer move. (a) A generically immersed surface~$f$ in an ambient $4$-manifold is pictured in blue, along with two Whitney circles. One Whitney arc in each Whitney circle is shown in green. (b) One green Whitney arc is changed, creating a new intersection between a Whitney disc and $f$. We also create an intersection between the two Whitney arcs. (c) The intersection between Whitney arcs is removed by a boundary pushoff operation, creating one more intersection between a Whitney disc and~$f$.}
    \label{fig:transfer-move}
\end{figure}
This move requires two Whitney discs on the same connected generically immersed surface $f$. As we see in the figure, the operation consists of first changing one of the Whitney arcs to create an intersection between two Whitney arcs. Next we remove this new intersection by the boundary pushoff operation. In this way, both of the two Whitney discs we started off with gain an intersection with $f$ in the interior. 

Having described the contents of \Cref{tab:table}, we now return to the proof of Theorem A.

\begin{step}\label{step:geometric-tricks}
    Use geometric manoeuvres to remove all type 2 and 3 problems, as well as all but at most one type 4 problem, i.e.~we arrange that $\{\mathring{W}_i\}$ and $f$ intersect in at most one point. 
\end{step}

Use boundary pushoff to solve all type 3 problems, creating more type 4 problems. To do this properly, first enumerate the Whitney arcs. Then work on the arcs in order. For the $i$th arc, push other arcs with index greater than $i$ off the $i$th arc, starting with one of the arcs closest to the endpoint, until the $i$th arc is disjoint from all other arcs. At the end of the process, all Whitney
arcs, and therefore Whitney circles, are mutually disjoint. Next, tube elements of $\{W_i\}$ into $g$ to remove all problems of type~4, creating new type 2 problems. Now we only have problems of type 1 and 2.

We would now like to remove all the problems of type 2, without creating new problems of type 3 and 4 in the process.\footnote{For example, if the current twisting numbers are all even, we could solve all type 2 problems using interior twisting, creating only new problems of type 1, which would complete the step.} For every $i$, perform interior twisting on~$W_i$ to arrange that $\tw(\partial W_i)$ is either $0$ or $1$. Only new problems of type 1 are created. The Whitney discs with trivial twisting number at this stage are ignored until the next step. Consider the Whitney discs with twisting number equal to one; call this set $\{\wt{W}_i\}_{i=1}^N$. Boundary twist each $\wt{W}_i$ to arrange that $\tw(\wt{W}_i)=0$, while creating a single intersection point of its interior with $f$. In case $N$ is even, pair up all the elements of $\{\wt{W}_i\}_{i=1}^N$. If $N$ is odd, set aside $\wt{W}_N$, and pair up the rest. Do the transfer move on each of the pairs we just assigned.  Now each element of~$\{\wt{W}_i\}_{i=1}^N$, except possibly $\wt{W}_N$, has two intersections with $f$. Tube $\wt{W}_i$ to $g$ at these new intersections. This solves all the type 4 problems within $\{\wt{W}_i\}_{i=1}^N$, except possibly for a single one in $\wt{W}_N$, while changing $\tw(\wt{W}_i)$ by $2e(\nu g)$ for each $i$, except possibly~$i=N$. Since~$2e(\nu g)$ is even, we can use interior twisting to solve these new type 2 problems, creating only type~1 problems in the process. This completes this step. Note that we have only type 1 problems and at most one type 4 problem left to solve. 

\begin{step}\label{step:BFF}
    If there are only type 1 problems left, proceed to the next step. If there is a type 4 problem remaining, stabilise to change the domain of $f$ to a torus, then do two \emph{band-fibre-finger moves} to remove the type 4 problem at the expense of adding in four new double points in $f$.
\end{step}

In this step we assume that we only have a single type 4 problem left to solve. In other words, the (generically immersed) Whitney discs $\{W_i\}$ have trivial twisting numbers as well as embedded and disjoint boundaries, and moreover $\{\mathring{W}_i\}\pitchfork f$ consists of a single point. By relabelling, we can assume that this intersection is with~$\mathring{W}_1$. 

Perform a trivial stabilisation of $f$. This means take the pairwise connected sum of $(M,f)$ with $(S^4,\Sigma)$, where $\Sigma$ is the standard, unknotted torus. Note that the meridian and longitude of $\Sigma$ bound embedded discs $V'_1$ and $V'_2$ in $S^4$, with interiors disjoint from $\Sigma$ and with $\partial V'_1\cap \partial V'_2$ a single point. After taking the pairwise connected sum, we can assume that these discs lie in $M$ as well. Let $D'_1$ and $D'_2$ be two embedded discs on $f$, with interiors pushed slightly in the normal direction. Construct the ambient connected sum of $V'_1$ with $D'_1$, and of $V'_2$ with $D'_2$, along embedded arcs in the ambient $4$-manifold. The result is a pair of embedded annuli,~$B_1$ and $B_2$, with boundaries $\partial V'_1\cup \partial D'_1$ and $\partial V'_2\cup \partial D'_2$ lying on $f$. Now we will need our final geometric manoeuvre, the \emph{band-fibre-finger move}, described next. 

Given a generically immersed surface $f$ in a $4$-manifold $M$ and an annulus $B\subseteq M$ with $\partial B$ lying in $f$, the band-fibre-finger move consists of doing a self-finger move on $f$ along one of the fibres in the annulus (\Cref{fig:BFF-move}). 
\begin{figure}[tb]
    \centering
    \begin{tikzpicture}
        \node[anchor=south west,inner sep=0] at (0,0){\includegraphics[width=11cm]{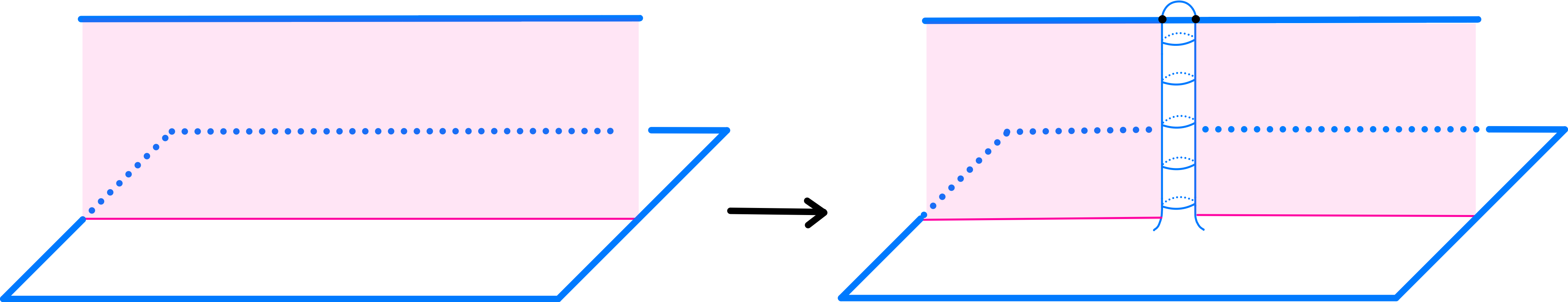}};
            \node at (0.25,-0.25) {$f$};
            \node at (0.75,2.25) {$f$};
            \node at (3.75,1.75) {$B$};
            \node at (6.24,-0.25) {$f$};
            \node at (6.74,2.25) {$f$};
            \node at (9.74,1.75) {$W_B$};
	\end{tikzpicture} 
    \caption{Band-fibre-finger move. Left: a generically immersed surface $f$ in an ambient $4$-manifold is shown in blue. Note that there are two sheets of the surface -- one appearing as a flat plane and the other as a line. A portion of an annulus $B$ with boundary on $f$ is shown in pink. Right: the band-fibre-finger move consists of doing a self-finger move on $f$ as shown. Two new double points of $f$ are created. These are paired by a Whitney disc $W_B\subseteq B$ which is shown in pink. }
    \label{fig:BFF-move}
\end{figure}
The two new double points created in this procedure are naturally paired by a trivial Whitney disc. However, under certain conditions, we get an alternate Whitney disc from the band $B$ minus a small strip along the finger move arc, as shown in \Cref{fig:BFF-move}. This is always the case when $M$ and $f$ are orientable. For the more general case, see \cite[{Construction~7.2}]{KPRT:sigmet}. 

Now we return to the proof. Do the band-fibre-finger move both along $B_1$ and along $B_2$. This changes $f$ by a regular homotopy, creating four new double points, and now the double points of $f$ are paired up by Whitney discs $\{W_i\}\cup \{V_1,V_2\}$,
where $V_1$ and $V_2$ come from $B_1$ and $B_2$ respectively. By construction, these Whitney discs have trivial twisting number and embedded boundaries. We also know that 
\[
\big(\{\mathring{W}_i\}\cup \{\mathring{V}_1,\mathring{V}_2\}\big)\pitchfork f
\]
is a single point in $\mathring{W}_1$. The boundaries are also disjoint, except $\partial V_1\pitchfork \partial V_2$ is a single point. So perform the boundary pushoff operation, to trade the Whitney arc intersection for an intersection between $V_1$ and $f$. Now the entire set $\{\mathring{W}_i\}\cup \{\mathring{V}_1,\mathring{V}_2\}$ intersects twice with $f$. The transfer move applied to $W_1$ and $V_1$ arranges that each (new) $W_1$ and $V_1$ intersects $f$ twice. Tube $W_1$ and $V_1$ into the geometric dual $g$ to remove these intersections. Each of $\tw(\partial W_1)$ and $\tw(\partial V_1)$ now equals an even number, namely $2e(\nu g)$. These twisting numbers can thus be changed back to zero by interior twisting, paying only the price of type 1 intersections. Now we have finally arrived at a collection of Whitney discs for $f\pitchfork f$ satisfying the hypotheses of \Cref{prop:whitney-by-DET}. 

\begin{step}\label{step:final}
    Apply \Cref{prop:whitney-by-DET} to $\{W_i\}$, then do the Whitney move on $f$ along the new Whitney discs $\{\ol{W}_i\}$. 
\end{step}

The hypotheses of \Cref{prop:whitney-by-DET} are satisfied by all of our previous work. The discs $\{\ol{W}_i\}$ produced by \Cref{prop:whitney-by-DET} are by construction locally flat embedded, with disjoint, embedded boundaries and with interiors disjoint from $f$. As previously discussed, doing the Whitney move on $f$ along these discs removes all the double points of $f$ and therefore results in a locally flat embedding, as desired. Note that under certain conditions we can bypass \Cref{step:BFF}, so we can obtain an embedded sphere rather than a torus. However, a torus is the best we can do in the general case. 
\end{proof}

\subsection{More general results} We end this section by stating the more general theorems that were proven by Lee--Wilczy\'{n}ski~\cite[{Theorem~1.1}]{LW97} and by  Kasprowski, Powell, Teichner, and the author in~\cite[{Theorem~1.2}]{KPRT:sigmet}. For the statement below, we remark that an embedding is said to be \emph{simple} if the fundamental group of the complement is abelian. The divisibility of a class $x\in H_2(M;\Z)$ is the least integer $d$ such that~$x=dy$ for some $0\neq y\in H_2(M;\Z)$. The signature of a $4$-manifold $M$ is denoted by $\sigma(M)$.

\begin{theorem}[{\cite[{Theorem~1.1}]{LW97}}]\label{thm:LW-general}
    Let $M$ be a compact, oriented, simply connected $4$-manifold whose boundary is a disjoint and possibly empty union of integral homology spheres. Suppose $x\in H_2(M;\Z)$ is a nonzero class of divisibility $d$. Then there exists a simple, locally flat embedding $\Sigma\hookrightarrow M$ representing $x$ by an oriented surface of genus $g>0$ if and only if 
    \[
    b_2(M)+2g\geq \max_{0\leq j\leq d} \left\lvert \sigma(M)-\frac{2j(d-j)}{d^2}x\cdot x\right\rvert.
    \]
\end{theorem}

Note that Theorem A is the case of $d=1$. This is a very powerful result, applicable in a variety of situations. There is a companion theorem~\cite[{Theorem~1.2}]{LW97} providing one further obstruction in the genus zero case, given by the Kervaire--Milnor condition relating the intersection number $x\cdot x$ to the Kirby--Siebenmann invariant of $M$ and the Rochlin invariant of the boundary $\partial M$. However, the condition on the fundamental group of the complement is essential, as is the requirement that the ambient $4$-manifold is either closed or has boundary a disjoint union of homology spheres. Roughly speaking, this is required due to the surgery theoretic approach used by Lee--Wilczy\'{n}ski. 

Now we state the result of \cite{KPRT:sigmet}. To do this we need to define the \emph{Kervaire--Milnor invariant}. 

\begin{definition}[{\cite[Definition~1.4]{KPRT:sigmet}}]
    Let $\Sigma$ be a surface and let $M$ be a $4$-manifold. Let~${F\colon (\Sigma,\partial \Sigma)\looparrowright (M,\partial M)}$ be a generic immersion. By definition, the \emph{Kervaire--Milnor invariant}, $\km(F)\in \Z/2$, vanishes if and only if, after finitely many finger moves taking $F$ to some $F'$, there is a collection of generically immersed Whitney discs $\{W_i\}$, pairing all the double points of $F'$, such that the boundaries are disjoint and embedded, the twisting numbers are trivial, and the interiors are disjoint from $F'$.
\end{definition}

\begin{theorem}[{\cite[{Theorem~1.2}]{KPRT:sigmet}}]\label{thm:KPRT-sigmet}
    Let $\Sigma$ be a compact surface with connected components $\{\Sigma_i\}_{i=1}^m$ and let $M$ be a connected $4$-manifold. Let~$F\colon (\Sigma,\partial \Sigma)\looparrowright (M, \partial M)$ be a generic immersion with components $\{f_i\colon (\Sigma_i,\partial \Sigma_i)\looparrowright (M, \partial M)\}_{i=1}^m$. 
    Suppose that $\pi_1(M)$ is a good group and that $F$ has algebraically dual spheres $G$, with components $\{g_i\colon S^2\looparrowright M\}_{i=1}^m$. In other words, $\lambda(f_i,g_j)=\delta_{ij}$. Then the following statements are equivalent.

    \begin{enumerate}[label=(\roman*)]
        \item \label{item-i-SET}  The intersection numbers $\lambda(f_i,f_j)$ for all $i < j$, the self-intersection numbers~$\mu(f_i)$ for all $i$, and the Kervaire--Milnor invariant $\km(F)\in\Z/2$, all vanish.
        \item \label{item-ii-SET} There is an embedding $\ol{F}=\{\ol{f}_i\}_{i=1}^m \colon (\Sigma, \partial\Sigma)  \hookrightarrow (M,\partial M)$, regularly homotopic to $F$ relative to $\partial \Sigma$, with geometrically dual spheres $\ol{G}=\{\ol{g}_i\colon S^2\looparrowright M\}_{i=1}^m$ such that $[\ol{g}_i]= [g_i]\in \pi_2(M)$ for all $i$.
    \end{enumerate}
\end{theorem}

If $\pi_1(M)$ is trivial, the intersection and self-intersection numbers in the theorem above are integers, obtained as a signed count. For more general fundamental groups, we have to use the equivariant versions (\Cref{sec:lambda-mu}) as mentioned in \Cref{rem:DET-nontrivial-pi1}. For the most general setting, where $\Sigma$ might have positive genus, see \cite[{Section~2}]{KPRT:sigmet} for the definitions of the equivariant intersection and self-intersection numbers. 

Unlike \Cref{thm:LW-general}, in \Cref{thm:KPRT-sigmet} we do not need to restrict to simple embeddings, or to $4$-manifolds with empty or homology sphere boundary. On the other hand, the algebraically dual sphere is essential, as is the condition that the ambient $4$-manifold has good fundamental group. This is due to our use of the disc embedding theorem. 

A helpful fact about \Cref{thm:KPRT-sigmet} is that we can often force the Kervaire--Milnor invariant to be trivial, by modifying the map $f$ in some way -- in the proof of  Theorem A we did this by stabilising. A similar proof gives the following corollary. 

\begin{corollary}\label{cor:simply-connected-pos-genus}
    Let $M$ be a $4$-manifold with $\pi_1(M)$ good and let $\Sigma$ be a connected, oriented surface with positive genus. Suppose we have a generic immersion $f\colon (\Sigma,\partial\Sigma)\looparrowright (M,\partial M)$ with vanishing self-intersection number and an algebraically dual sphere. Then~$f$ is regularly homotopic, relative to~$\partial \Sigma$, to an embedding.
\end{corollary}

For other applications, including to non-orientable surfaces and to slicing knots in general $4$-manifolds, see~\cite{KPRT:sigmet}.

\section{Embedding surfaces using surgery theory}
\label{sec:alex-poly-one}
In this section, we switch gears and describe a more indirect strategy to construct locally flat surfaces in a given $4$-manifold. The procedure described here can be effectively encapsulated in the so-called \emph{surgery sequence}, as we briefly describe later in \Cref{sec:surgery-sequence}. We will sketch the proof of Theorem B in \Cref{sec:proof-alex-poly-one}. Before that we need to recall a number of necessary ingredients -- namely a 0-surgery characterisation of sliceness in \Cref{sec:0-surgery-char-sliceness}, equivariant intersection and self-intersection numbers in \Cref{sec:lambda-mu}, and the sphere embedding theorem in \Cref{sec:sphere-embedding-theorem}. First we restate Theorem B, after recalling a relevant definition. 

\begin{definition}
    A knot $K\colon S^1\hookrightarrow S^3$ is \emph{(topologically) slice} if it extends to a proper locally flat embedding $\Delta$ of a disc in $B^4$. So we have 
    \[
    \begin{tikzcd}
        S^1\ar[r,hook,"K"]\ar[d, hook]  &S^3\ar[d,hook]\\
        D^2\ar[r,hook,"\Delta"] &B^4,        
    \end{tikzcd}
    \]
    where the vertical maps are the inclusions. The disc $\Delta$ is called a \emph{(topological) slice disc} for $K$.  
\end{definition}

Slice knots were first introduced by Fox and Milnor in the 1950s. Since then they have become an active area of study. For more details on slice knots see, e.g., \cite{livingston-slice-survey,winterbraids-slice-survey}. 

We now recall the statement of Theorem B. For a knot $K\subseteq S^3$ with Seifert matrix $V$, the Alexander polynomial $\Delta_K(t)\in \Z[t,t^{-1}]$ is defined as $\det(tV-V^T)$. For more details, see e.g.~\cite{gordon-knot-survey,rolfsen-book}. It can be easily computed from a genus one Seifert surface that all (untwisted) Whitehead doubles have Alexander polynomial one, so the following theorem gives numerous examples of nontrivial slice knots. 

\begin{theoremB}
    Every knot $K\colon S^1\hookrightarrow S^3$ with Alexander polynomial one is (to\-po\-lo\-gi\-ca\-lly) slice. 
\end{theoremB}

A proof of Theorem B using surgery theory was given in \cite[Theorem~7]{freedman-icm} and \cite[Theorem~11.7B]{FQ} (see also~\cite[{Theorem~1.14}]{DET-book-context}). An alternative, more direct proof is given in \cite{garoufalidis-teichner}, using a single application of the disc embedding theorem for a finite collection of discs, where the ambient manifold has infinite cyclic fundamental group.

\begin{remark} 
    The results~\cite[{Theorems~1.13 and~1.14}]{F} are commonly, but erroneously, cited for Theorem B. In fact, neither of these results match Theorem B. \cite[{Theorem~1.13}]{F}
    states that every Alexander polynomial one knot bounds an embedded, locally homotopically unknotted disc in $B^4$, but this was not shown to be locally flat. (Local flatness follows from later work of Quinn~\cite[{Theorem~9.3A}]{FQ} (see also the correction in \cite{venema:1alg}.)) The second result often cited,~\cite[{Theorem~1.14}]{F}, only asserts that the untwisted Whitehead double of a knot with Alexander polynomial one is topologically slice. Moreover, the proofs of~\cite[{Theorems~1.13 and~1.14}]{F} both rely on \cite[{Lemma~2}]{Freedman-Alex}, and  a counterexample to this lemma was presented in~\cite{garoufalidis-teichner}.

    Indeed, Theorem B above was never claimed by Freedman in~\cite{F}. The first proof that Alexander polynomial one knots are slice was given in \cite[{Theorem~7}]{freedman-icm} and makes crucial use of Quinn's work in~\cite{quinn:endsIII}, by working purely in the topological setting. Therefore we choose to attribute the result to both Freedman and Quinn.
\end{remark}

In \Cref{sec:proof-alex-poly-one} we will give a substantially expanded version of the proof of Theorem B given in \cite[Theorem~7]{freedman-icm}, \cite[Theorem~7]{FQ}, and \cite[Theorem~1.14]{DET-book-context}, unpacking the surgery technology. The proof will require some additional background, which we provide in the next three subsections.  

\subsection{Characterising sliceness using the 0-surgery}\label{sec:0-surgery-char-sliceness}
Suppose that $K$ is a slice knot with a slice disc $\Delta$. Let $\mathring{\nu} \Delta$ denote an open tubular neighbourhood of $\Delta$ in $B^4$. Observe that $\partial(B^4\sm \mathring{\nu} \Delta)$ is the result of 0-framed Dehn surgery on $S^3$ along $K$, denoted by $S^3_0(K)$ (Exercise~\ref{ex:slice-disc-complement-zero-surgery}). So when $K$ is slice, the 0-surgery $S^3_0(K)$ is the boundary of $W:=B^4\sm \mathring{\nu} \Delta$, where we can further check that the inclusion induced map $\Z\cong H_1(S^3_0(K);\Z)\to H_1(W;\Z)$ is an isomorphism; the fundamental group $\pi_1(W)$ is normally generated by the meridian $\mu_K$ of $K$, considered to lie in $S^3_0(K)$; and $H_2(W;\Z)=0$. It turns out that the converse is also true, yielding the following characterisation of sliceness. We leave the proof as an exercise (Exercise~\ref{ex:0-surgery-char-sliceness}).

\begin{theorem}\label{thm:0-surgery-char-sliceness}
    A knot $K\subseteq S^3$ is $($topologically$)$ slice if and only if the 0-framed Dehn surgery $S^3_0(K)$ is the boundary of some compact, connected $4$-manifold $W$ such that
    \begin{enumerate}
        \item the inclusion induced map $\Z\cong H_1(S^3_0(K);\Z)\to H_1(W;\Z)$ is an isomorphism;
        \item the fundamental group $\pi_1(W)$ is normally generated by the meridian $\mu_K$ of $K$, considered to lie in $S^3_0(K)$; and 
        \item the second homology group $H_2(W;\Z)=0$.
    \end{enumerate}
\end{theorem}

\subsection{Equivariant intersection and self-intersection numbers}\label{sec:lambda-mu}
In this subsection we briefly describe the equivariant intersection and self-intersection numbers $\lambda$ and $\mu$ respectively. For a more detailed account, see e.g.~\cite[Chapter~5]{wall-surgery-book}, \cite[Section~1.7]{FQ} and \cite[Section~11.3]{DET-book-DETintro}.

\begin{figure}[tb]
    \centering
    \begin{tikzpicture}
        \node[anchor=south west,inner sep=0] at (-0.25,0){\includegraphics[width=11cm]{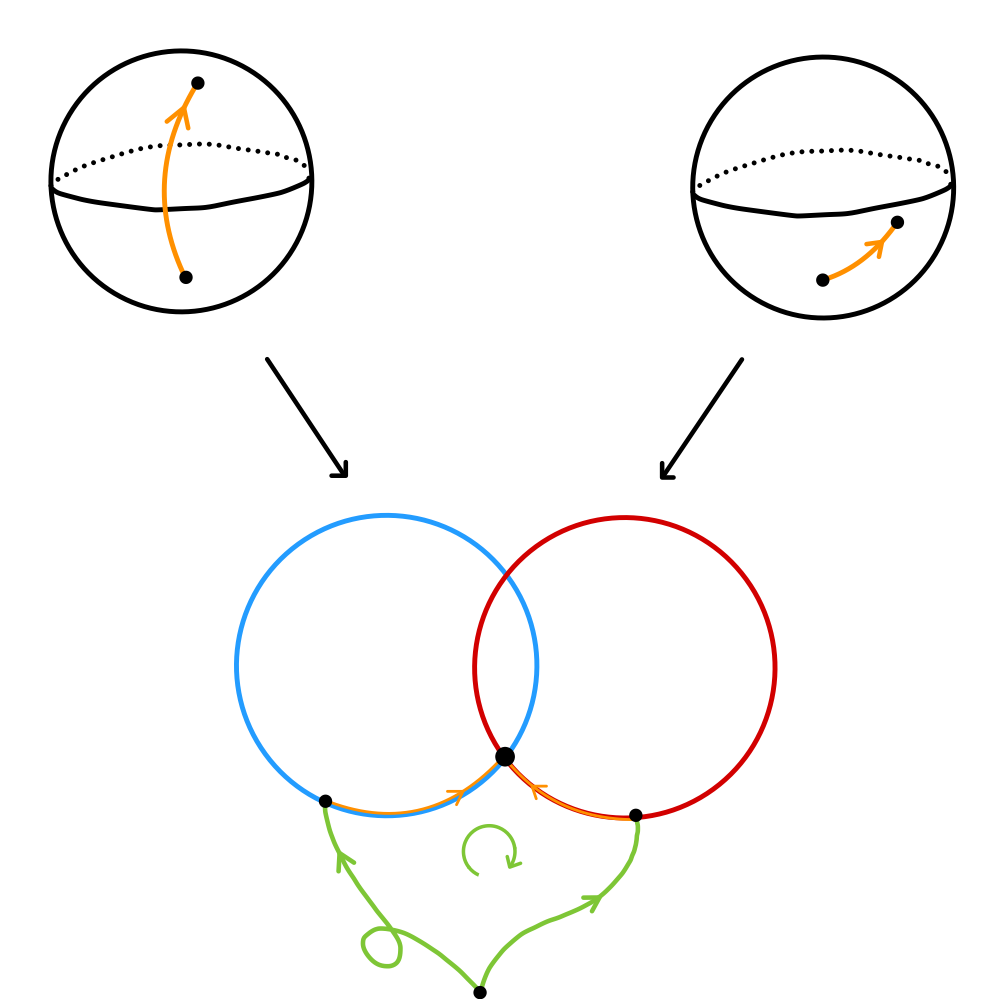}};
        \node at (5,-0.25) {$m$};
        \node at (3.3,6.5) {$f$};
        \node at (7.3,6.5) {$g$};
        \node at (2,7.9) {$x$};
        \node at (8.55,7.88) {$x$};
        \node at (1.65,10.1) {$p_f$};
        \node at (9.9,8.5) {$p_g$};
        \node at (3.5,1) {$w_f$};
        \node at (6.6,1) {$w_g$};
        \node at (5.6,2.65) {$p$};
        \node at (4.15,2.3) {$\alpha^p_f$};
        \node at (6.1,2.3) {$\alpha^p_g$};
        \node at (5.25,1.2) {$\gamma_p$};
        \node at (2.9,2) {$f(x)$};
        \node at (7.2,1.8) {$g(x)$};
	\end{tikzpicture} 
    \caption{Defining the equivariant intersection number. Two generically immersed spheres $f$ and $g$ are shown, in blue and red respectively, in an ambient $4$-manifold $M$ with basepoint $m$. The sphere has basepoint $x$. The whiskers $w_f$ and $w_g$ are shown in green. The immersions $f$ and $g$ intersect at a point $p\in M$ with $f(p_f)=p=g(p_g)$. The arcs $\alpha^p_f$ and $\alpha^p_g$, as well as their preimages in $S^2$, are shown in orange. The element in $\pi_1(M,m)$ associated with $p$, denoted by $\gamma_p$ is indicated as a circular arrow.}
    \label{fig:lambda}
\end{figure}

Let $M$ be a connected, oriented, topological $4$-manifold and choose a basepoint~$m\in M$. Consider two generic immersions $f,g\colon S^2\looparrowright M$ which intersect each other transversely. Choose a basepoint $x\in S^2$ and an orientation for $S^2$. Choose paths $w_f,w_g\colon [0,1]\to M$ with initial points $w_f(0)=w_g(0)=m$ and terminal points $w_f(1)=f(x)$ and $w_g(1)=g(x)$. These paths are called \emph{whiskers} for $f$ and~$g$. Define the following sum
\[
    \lambda(f,g):=\sum_{p\in f \pitchfork g} \varepsilon_p  \gamma_p \in \Z[\pi_1(M,m)],
\]
where
\begin{itemize}
    \item the path $\alpha_f^p$ is the image under $f$ in $M$ of a path in $S^2$ from $x$ to the preimage of~$p$ under $f$; 
    \item the path $\alpha^p_g$ is the image under $g$ in $M$ of a path in $S^2$ from $x$ to the preimage of~$p$ under $g$; 
    \item the element $\varepsilon_p\in\{\pm 1\}$ is the sign of the intersection point $p$; and 
    \item the loop $\gamma_p$ is the element of $\pi_1(M,m)$ given by the concatenation 
    \[
        w_f  \alpha_f^p  (\alpha_g^p)^{-1}  w_g^{-1}.
    \]
    \end{itemize}
The quantity $\lambda(f,g)$ is called the \emph{equivariant intersection number} of $f$ and $g$. Note that when $\pi_1(M)=1$, the sum $\lambda(f,g)$ is simply the signed count of intersections between $f$ and $g$. For more details on why $\lambda$ is well defined in general see Exercise~\ref{ex:lambda-mu-basic-properties}.

\begin{figure}[tb]
    \begin{tikzpicture}
        \node[anchor=south west,inner sep=0] at (0,0){\includegraphics[width=10cm]{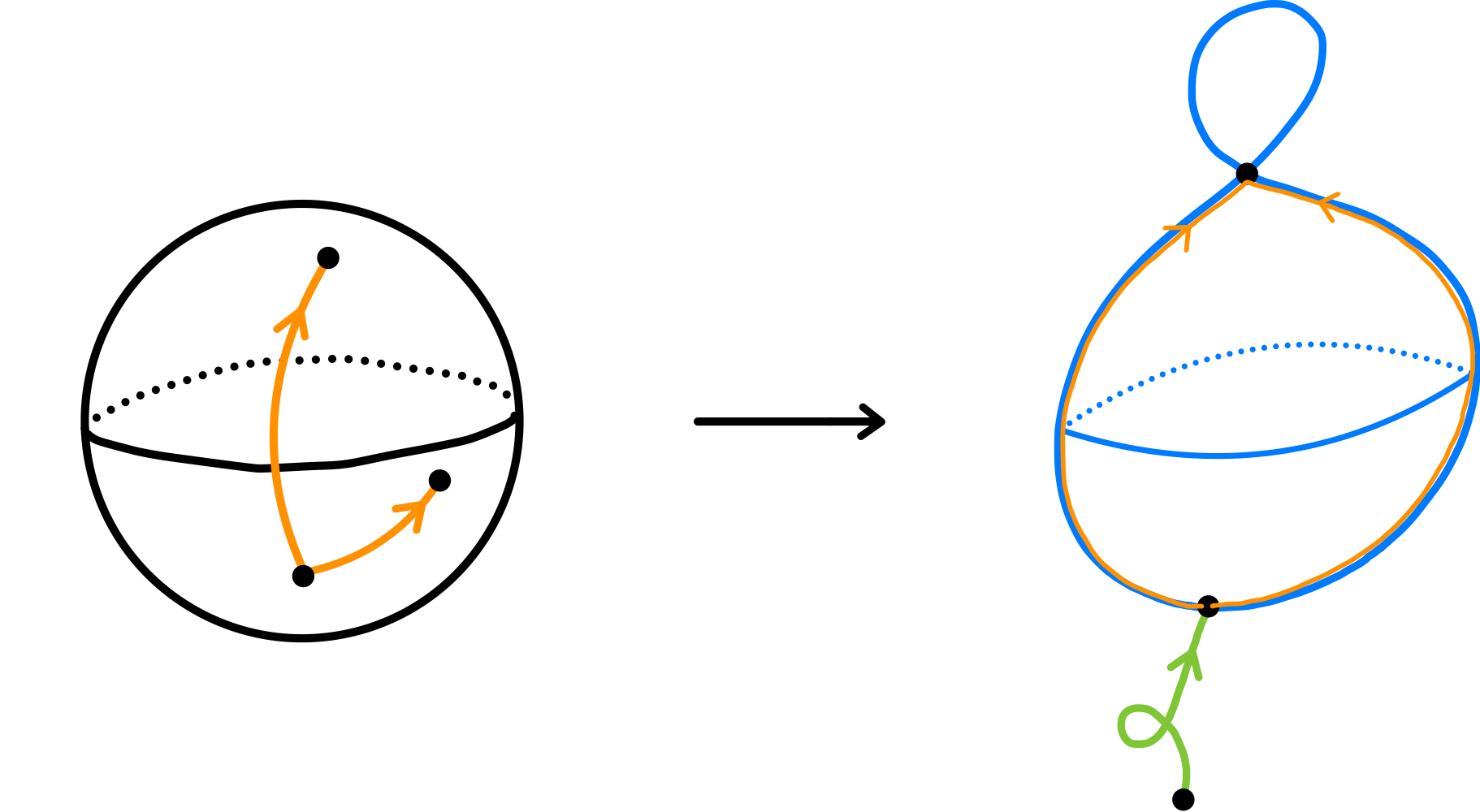}};
        \node at (8,-0.25) {$m$};
        \node at (5.25,3) {$f$};
        \node at (2,1.35) {$x$};
        \node at (8.6,1.1) {$f(x)$};
        \node at (8,4.25) {$p$};
        \node at (7.25,0.5) {$w_f$};
        \node at (2,3.8) {$p_1$};
        \node at (3.2,2.2) {$p_2$};
        \node at (6.85,2) {$\alpha^p_1$};
        \node at (9.6,1.7) {$\alpha^p_2$};
	\end{tikzpicture}
    \caption{Defining the equivariant self-intersection number. A generic immersion $f$ of a sphere is shown in blue in an ambient $4$-manifold $M$ with basepoint $m$. The sphere has basepoint $x$. The whisker $w_f$ is shown in green. The point $p=f(p_1)=f(p_2)$ is a transverse self-intersection in $f$. The two arcs $\alpha^p_1$ and $\alpha^p_2$, as well as their preimages are shown in orange.}
    \label{fig:mu}
\end{figure}

Similarly, we define 
\[
    \mu(f):= \sum_{p\in f \pitchfork f} \varepsilon_p \gamma_p\in \Z[\pi_1(M,m)],
\]
where 
\begin{itemize}
    \item the paths $\alpha_1^p$ and $\alpha_2^p$ are images under $f$ in $M$ of paths in $S^2$ from $x$ to the two distinct preimages of $p$; 
    \item the element $\varepsilon_p\in\{\pm 1\}$ is the sign of the intersection point $p$; and
    \item the loop $\gamma_p$ is the element of $\pi_1(M,m)$ given by the concatenation 
    \[
        w_f \alpha_1^p (\alpha_2^p)^{-1}  {w_f}^{-1}.
    \]
    \end{itemize}
The quantity $\mu(f)$ is called the \emph{equivariant self-intersection number} of $f$. Again, when $\pi_1(M)=1$, the sum $\mu(f)$ coincides with the signed count of self-intersections of $f$. One has to be slightly careful in the definition of $\mu$: as we indicate in Exercise~\ref{ex:lambda-mu-basic-properties}, the sum $\mu$ is well defined only in a quotient of $\Z[\pi_1(M,m)]$.

The vanishing of the intersection and self-intersection numbers has a nice characterisation in terms of the existence of generically immersed Whitney discs pairing up all double points (Exercise~\ref{ex:lambda-mu-zero-whitney}). 

\subsection{The sphere embedding theorem}\label{sec:sphere-embedding-theorem}

Now that we have defined the equivariant intersection and self-intersection numbers, we can finally state the most general known version of the disc embedding theorem. For the definition of a good group, see \Cref{rem:DET-nontrivial-pi1} and \cite{Freedman-Teichner:1995-1,DET-book-goodgroups}.

\begin{theorem}[Disc embedding theorem, general version; \cite{F}, \cite{freedman-icm}, {\cite[Theorem~5.1A]{FQ}}, and {\cite{Powell-Ray-Teichner:2018-1}}]\label{thm:DET-full}
    Let $M$ be a connected $4$-manifold such that $\pi_1(M)$ is a good group. Let
    \[F = (f_1, \dots, f_n)\colon (D^2 \sqcup \cdots \sqcup D^2,S^1  \sqcup \cdots \sqcup S^1)\looparrowright (M,\partial M)\]
    be a generically immersed collection of discs in $M$ with pairwise disjoint, locally flat, embedded boundaries.  Suppose there is a generically immersed collection of spheres
    \[
    G = (g_1,\dots,g_n) \colon S^2  \sqcup \cdots \sqcup S^2 \looparrowright M,
    \]
    which is \emph{algebraically dual} to $F$, i.e.~$\lambda(f_i,g_j)=\delta_{ij}$ for all~$i,j$.   Assume further that each $g_i$ has trivial normal vector bundle and $\lambda(g_i,g_j)=0=\mu(g_i)$ for all~${i,j=1,\dots,n}$.
    
    Then there exists a collection of pairwise disjoint, locally flat embedded discs
    \[
    \overline{F} = (\overline{f}_1,\dots,\overline{f}_n) \colon (D^2  \sqcup \cdots \sqcup D^2,S^1  \sqcup \cdots \sqcup S^1) \hookrightarrow (M,\partial M),
    \]
    and generically immersed spheres
    \[
    \overline{G} = (\overline{g}_1,\dots,\overline{g}_n) \colon S^2  \sqcup \cdots \sqcup S^2 \looparrowright M,
    \]
    which are geometrically dual to $\ol{F}$, i.e.~$\ol{f}_i\pitchfork \ol{g}_j$ is empty if $i\neq j$ and a single point otherwise. Moreover, for every $i$, the discs $\overline{f}_i$ and $f_i$ have the same framed boundary and $\overline{g}_i$ is homotopic to~$g_i$.
    \end{theorem}

We can apply the above to prove the \emph{sphere embedding theorem}, which we state next. We leave the proof as an exercise for highly motivated readers (Exercise~\ref{ex:sphere-emb-thm}). 

\begin{theorem}[{Sphere embedding theorem~\cite[{Theorem~5.1B}]{FQ}}]\label{thm:sphere-emb-thm}
    Let $M$ be a connected $4$-manifold such that~$\pi_1(M)$ is a good group. Let 
    \[
    F= (f_1,\dots,f_n) \colon S^2 \sqcup \cdots \sqcup S^2 \looparrowright M
    \]
    be a generically immersed collection of spheres in $M$ with $\lambda(f_i,f_j)=0$ for every~$i \neq j$ and $\mu(f_i)=0$ for all $i$. Suppose there is a generically immersed collection of spheres
    \[
    G= (g_1,\dots,g_n) \colon S^2 \sqcup \cdots \sqcup S^2 \looparrowright M
    \]
    which is algebraically dual to $F$, i.e.~$\lambda(f_i,g_j) = \delta_{ij}$ for all $i,j$. Assume further that each $g_i$ has trivial normal vector bundle. 
    
    Then there exists a collection of pairwise disjoint, locally flat embedded spheres
    \[
    \ol{F} = (\ol{f}_1,\dots,\ol{f}_n) \colon S^2 \sqcup \cdots \sqcup S^2 \hookrightarrow M,
    \]
    with each $\ol{f}_i$ regularly homotopic to $f_i$, and generically immersed spheres
    \[
    \ol{G}= (\ol{g}_1,\dots,\ol{g}_n) \colon S^2 \sqcup \cdots \sqcup S^2 \looparrowright M,
    \]
    which are geometrically dual to $\ol{F}$, i.e.~$\ol{f}_i\pitchfork \ol{g}_j$ is empty if $i\neq j$ and a single point otherwise. Moreover, for every $i$, the sphere $\ol{g}_i$ has trivial normal vector bundle and is homotopic to $g_i$.   
\end{theorem}

\subsection{Proof of Theorem~\ref{thm:alex-poly-one}}\label{sec:proof-alex-poly-one}

We can finally sketch the proof of Theorem B. 

\begin{proof}[Proof of Theorem B]
    We will use the 0-surgery characterisation of sliceness (\Cref{thm:0-surgery-char-sliceness}), i.e.~we will build a $4$-manifold $W$ with $\partial W=S^3_0(K)$ satisfying the conditions given in \Cref{thm:0-surgery-char-sliceness}. Since $H_1(S^3_0(K);\Z)\cong \Z$ with generator a meridian $\mu_K$ of $K$, there exists a map $f\colon S^3_0(K)\to S^1$ such that the induced map on fundamental groups sends $[\mu_K]\mapsto 1$. Note that it will suffice to build $W$ so that we have an extension to a homotopy equivalence
    \begin{equation}\label{eq:alex-poly-one-goal}
    \begin{tikzcd}
        S^3_0(K)\ar[r,"f"]\ar[d,hook]  &S^1,\\
        W\ar[ur,"\simeq"']
    \end{tikzcd}
    \end{equation}
    where the vertical map is the inclusion of the boundary. Once again we break up the proof into multiple steps. 

    \setcounter{step}{0}

    \begin{step}
    Find an arbitrary spin null bordism of $(S^3_0(K), f)$ over $S^1$.    
    \end{step}
    
    Consider $\Omega_3^\spin(S^1)$, the 3-dimensional spin bordism group over $S^1$. By definition, elements of this group are represented by pairs $(Y,\alpha\colon Y\to S^1)$, where $Y$ is a spin~$3$-manifold, and where two such elements $(Y,\alpha)$ and $(Y',\alpha')$ are identified if there is an extension
    \[
    \begin{tikzcd}
        Y\ar[dr,"\alpha"]\ar[d,hook]    &\\
        W\ar[r] &S^1,\\
        Y'\ar[ur,"\alpha'"']\ar[u,hook'] &
    \end{tikzcd}
    \]
    where $W$ is a spin $4$-manifold with boundary $\partial W=Y\sqcup Y'$ and the vertical maps are inclusions. Note that here we mean that the spin structure on $W$ induces the given spin structures on $Y$ and $Y'$.

    We will now need the Arf invariant of $K$, denoted by $\Arf(K)\in \Z/2$. While it is merely valued in $\Z/2$, this is a surprisingly powerful knot invariant with a number of equivalent definitions. For example, 
    \[
    \Arf(K)=\begin{cases}
        0,   &\text{ if } \Delta_K(-1) \equiv \pm 1 \mod{8};\\
        1,   &\text{ if }\Delta_K(-1) \equiv \pm 3 \mod{8};
    \end{cases}
    \]
    by~\cite[Proposition~3.4]{levine:polynomial-invariants-codimension-two} and~\cite[Theorem~2]{murasugi:arf-invariant}. This shows that since $K$ has Alexander polynomial one, $\Arf(K)=0$. We give another definition of $\Arf(K)$ using Whitney discs in \Cref{ex:arf-via-whitney}. 
    
    Let $\mathfrak{s}$ denote either of the two spin structures on $S^3_0(K)$. 
    There is an isomorphism
    \[{\Omega_3^\spin(S^1)\xrightarrow{\cong}\Z/2}\]
    under which $(S^3_0(K),\mathfrak{s})$ is mapped to $\Arf(K)$, which can be taken as another definition of $\Arf(K)$. The proof is a computation using the Atiyah--Hirzebruch spectral sequence, which shows that $\Omega^\spin(S^1)\cong \Omega^\spin(*)\cong \Z/2$, and uses yet another definition of $\Arf(K)$ in terms of the spin structures on an arbitrary, capped-off Seifert surface for $K$. 
    
    In our case, since $\Arf(K)=0$, we conclude there is a connected, spin $4$-manifold~$V$ with~$\partial V=S^3_0(K)$, inducing the given spin structure, and a map~${F\colon V\to S^1}$ such that we have the diagram
    \begin{equation}\label{eq:V-boundary}
    \begin{tikzcd}
        S^3_0(K)\ar[r,"f"]\ar[d,hook]  &S^1,\\
        V\ar[ur,"F"']
    \end{tikzcd}
    \end{equation}
    where the vertical map is inclusion. We cannot assume that the map $F$ is a homotopy equivalence, which would complete the proof. But we will see that we can modify~$F$ (and $V$), so that the end result is a homotopy equivalence. That is the content of the rest of the proof. 

    Recall that by Whitehead's theorem if we can arrange that $F_*\colon \pi_i(V)\to \pi_i(S^1)$ is an isomorphism for all $i$, then we can conclude that $F$ is a homotopy equivalence. By Poincar\'{e}--Lefschetz duality, it will suffice to arrange that $F_*$ is an isomorphism on $\pi_1$ and $\pi_2$. 

    \begin{step}
        Arrange that $F$ induces an isomorphism on $\pi_1$. 
    \end{step}
    
    The map $F_*\colon \pi_1(V)\to \pi_1(S^1)$ is already surjective by construction. We can modify $V$, and $F$, so that it is also injective, by performing \emph{surgery on circles}, as in Exercise~\ref{ex:surgery-on-circles}. As we will see in the exercise, there are two possible framing choices for each such surgery. We have to use the framing induced by the spin structure to ensure that we still have a diagram as in \eqref{eq:V-boundary}, where $V$ is spin and induces the given spin structure on $S^3_0(K)$. 
    
    We assume henceforth that we have already arranged that $F$ induces an isomorphism on fundamental groups. 
    
    \medskip 
    Most of the rest of the proof consists of showing that $F$ can be modified so that the result induces an isomorphism on $\pi_2$. Since $\pi_2(S^1)$ is trivial, we want to modify~$V$, while ensuring there is still a compatible map to $S^1$, so that $\pi_2(V)$ is also trivial. At present though, $\pi_2(V)$ is some unknown $\Z[\pi_1(V)]$-module.
    
    \begin{step}
        Replace $V$ with some spin $V'$ with \emph{hyperbolic} intersection form. 
    \end{step}
    
     We know that $\pi_1(V)\cong \pi_1(S^1)\cong \Z$ by the previous step in the proof. Recall that the Alexander polynomial of $K$ annihilates the Alexander module $H_1(S^3_0(K);\Z[\Z])$. Since the Alexander polynomial of $K$ is one, this means that $H_1(S^3_0(K);\Z[\Z])=0$, which in turn implies that the equivariant intersection form 
    \[
        \lambda_V\colon \pi_2(V)\times \pi_2(V)\to \Z[\Z]
    \]
    is nonsingular.\footnote{This is the equivariant analogue of the fact that the integral intersection form on a compact, connected $4$-manifold is nonsingular if the boundary is an integer homology 3-sphere (recall that~$\pi_2(V)\cong H_2(V;\Z[\pi_1(V)])$). In more detail, the intersection form \[\lambda_V\colon H_2(V;\Z[\Z])\times H_2(V;\Z[\Z])\to \Z[\Z]\] is given by $\langle x,y\rangle \mapsto \langle \mathrm{PD}^{-1}(\iota(y)),x\rangle$, where $\iota\colon H_2(V;\Z[\Z])\to H_2(V,\partial V;\Z[\Z])$ is induced by inclusion. There are two possible sources of singularity of $\lambda$: the kernel of $\iota$ and the kernel of the map $\kappa\colon H_2(V;\Z[\Z])\to \mathrm{Hom}_{\Z[\Z]}(H_2(V;\Z[\Z]),\Z[\Z])$. Both kernels are trivial when~$H_1(\partial V;\Z[\Z])=0$.} Since $V$ is spin, there is a unique regular homotopy class within each element of $\pi_2(V)$ with trivial normal vector bundle, by \Cref{thm:generic-immersions-bijection}. In more detail, we observe that for any $\delta\in\pi_2(V)$ represented by a generically immersed $2$-sphere $h$, the value $h\cdot h$ is even since $V$ is spin, so we can perform interior twisting on $h$ to arrange that the euler number of the normal vector bundle of $h$ is trivial. Two representatives $h$ and $h'$ of $\delta$ are regularly homotopic if and only if their normal vector bundles have equal euler numbers by \Cref{thm:generic-immersions-bijection}. The self-intersection number $\mu_V$ is well defined on regular homotopy classes. By evaluating $\mu_V$ on the unique representative with trivial normal vector bundle, we get a map $\mu_V$ on $\pi_2(V)$.\footnote{Giving away part of the answer to Exercise~\ref{ex:lambda-mu-basic-properties}\,~\eqref{item:mu-domain-codomain}, this map is valued in~${\Z[\pi_1(V)]/g\sim g^{-1}}$. } 

    We consider now the triple $(\pi_2(V),\lambda_V, \mu_V)$. One can check that this is a \emph{nonsingular quadratic form}, i.e.~$\lambda_V$ is a sesquilinear, Hermitian, nonsingular form on the finitely generated, free $\Z[\pi_1(V)]$-module $\pi_2(V)$ with quadratic refinement $\mu_V$. We do not explain these terms further, except to say that such non-singular quadratic forms are precisely the elements of the \emph{$L$-group} $L_4(\Z[\pi_1(V)])$, modulo so-called (stably) \emph{hyperbolic quadratic forms}. We will address hyperbolic forms presently. For now, we note that in our case we have $\pi_1(V)\cong \Z$ and $L_4(\Z[\Z])$ is well-understood. Indeed we know that $L_4(\Z[\Z])\cong 8\Z$~\cite{shaneson-splitting}, generated by the so-called \emph{$E8$-form}, with the isomorphism given by the signature. We do not need to know what the~$E8$-form is precisely, except to note that it is a major result of Freedman~\cite[{Theorem~1.7}]{F} that there is a closed, spin, topological $4$-manifold\footnote{Spin structures, and relatedly Stiefel--Whitney classes, can be defined for topological manifolds (see e.g.~\cite[Chapter~7]{FNOP:4dguide}). We note here that the $4$-manifold $V$ could have been chosen to be smooth, by choosing a smooth spin null-bordism for $S^3_0(K)$ and surgering on smoothly embedded circles. However, since $E8$ is not smoothable, we can no longer assume that $V'$ is smooth.} called the~\emph{$E8$-manifold}, denoted by $E8$, which realises this quadratic form as the intersection form. So if $(\pi_2(V),\lambda_V,\mu_V)\in L_4(\Z[\Z])\cong 8\Z$ corresponds to $n\in 8\Z$, we can replace $V$ with~$V':=V\# -nE8$. Under the isomorphism $L_4(\Z[\Z])\cong 8\Z$, we see that $(\pi_2(V'),\lambda_{V'},\mu_{V'})$ maps to $n-n=0$. (For the effect of connected sum on the intersection form, see~\Cref{ex:equiv-int-form-conn-sum}.)
    
    We have arranged that the element $(\pi_2(V'),\lambda_{V'}, \mu_{V'})$ is trivial in $L_4(\Z[\Z])$, which means by definition that $(\pi_2(V'),\lambda_{V'}, \mu_{V'})$ is a hyperbolic quadratic form, possibly after taking the connected sum of $V'$ with more copies of $S^2\times S^2$.
    
    By construction, $V'$ is spin, since it is a connected sum of spin $4$-manifolds, and moreover there is still a map 
    \[
    \begin{tikzcd}
        S^3_0(K)\ar[r,"f"]\ar[d,hook]  &S^1\\
        V'\ar[ur,"F'"']
    \end{tikzcd}
    \]

    \begin{step}
        Apply the sphere embedding theorem to realise half a basis of $\pi_2(V')$ by locally flat, pairwise disjoint, embedded spheres, which are equipped with a family of geometrically dual spheres.
    \end{step}

    We have arranged that the intersection form on $V'$ is hyperbolic, which by definition means that the second homotopy group $\pi_2(V')$ has a basis of generically immersed spheres $\{f_1,\dots, f_n, g_1, \dots, g_n\}$, for some $n$, where each $f_i$ and $g_i$ has trivial normal vector bundle and such that 
    \begin{enumerate}[label=(\roman*)]
        \item $\lambda(f_i,g_j)=\delta_{ij}$, for all $i,j$; 
        \item $\lambda(f_i,f_j)=\lambda(g_i,g_j)=0$, for all $i,j$; and 
        \item $\mu(f_i)=0=\mu(g_i)$, for all $i$.
    \end{enumerate}
    We also know that $\Z$ is a good group (see \Cref{rem:DET-nontrivial-pi1}), so we can apply the sphere embedding theorem (\Cref{thm:sphere-emb-thm}). The theorem replaces the collection~${\{f_1,\dots, f_n\}}$ with a collection $\{\ol{f}_1,\dots, \ol{f}_n\}$ of pairwise disjoint locally flat embeddings, with each $\ol{f}_i$ regularly homotopic to $f_i$. It also provides a collection of generically immersed geometrically dual spheres $\{\ol{g}_1,\dots, \ol{g}_n\}$, i.e.~$\ol{f}_i\pitchfork \ol{g}_j$ is empty if $i\neq j$ and a single point otherwise. Moreover, for every $i$, the sphere $\ol{g}_i$ has trivial normal vector bundle and is homotopic to $g_i$. 

    \begin{step}
        Perform \emph{surgery} on the embedded, disjoint half-basis of $\pi_2(V')$ found in the previous step. Check that the resulting manifold $W$ satisfies the conditions of \Cref{thm:0-surgery-char-sliceness}.
    \end{step}

    Each $\ol{f}_i$ has trivial normal vector bundle, since it is regularly homotopic to $f_i$ which has trivial normal vector bundle (\Cref{thm:generic-immersions-bijection}). So there is a tubular neighbourhood $\nu \ol{f}_i$ of each $\ol{f}_i$ which is homeomorphic to $S^2\times D^2$. Then we \emph{perform surgery on $\{\ol{f}_1,\dots \ol{f}_n\}$}, i.e.~for each $i$, we cut out the tubular neighbourhood of $\ol{f}_i$ and glue in a copy of$D^3\times S^1$. This results in the manifold
    \[
        W:= V'\sm \Big(\bigcup_i \nu\ol{f}_i\Big) \cup \Big(\bigcup_i (D^3\times S^1)\Big). 
    \]
    We leave it to the reader to verify that there is a homotopy equivalence $W\xrightarrow{\simeq} S^1$ as in \eqref{eq:alex-poly-one-goal}. The geometrically dual spheres $\{\ol{g}_i\}$ ensure that the fundamental group of $W$ is still $\Z$ (cf.~Exercise~\ref{ex:geometric-dual-pi1}). This completes the sketch of the proof of Theorem B.
\end{proof}

\subsection{The surgery sequence}\label{sec:surgery-sequence}
The proof strategy used in the previous subsection can be systematised greatly. We briefly describe this here, and refer the reader to e.g.~\cite{wall-surgery-book,CLM:surgery-book,kirby-taylor,DET-book-surgery} for more details. Let $X$ be a closed, oriented $4$-manifold. If $\pi_1(X)$ is a good group, then we have the following exact sequence of pointed sets, called the \emph{surgery exact sequence}. Indeed the sequence continues on the left, and the sets can be endowed with a group structure, but we ignore these matters for this brief treatment. 

\begin{equation}\label{eq:surgery-sequence}
\begin{tikzcd}
    \calS(X)\ar[r] &\calN(X)\ar[r,"\sigma"]    &L_4(\Z[\pi_1(X)])
\end{tikzcd}    
\end{equation}

We have already encountered the $L$-group $L_4(\Z[\pi_1(X)])$ in the previous subsection. The set $\calN(X)$ is the set of degree one \emph{normal maps}: its elements are degree one maps $V\to X$, compatible with the stable normal vector bundles, where $V$ is a closed $4$-manifold, modulo degree one normal bordism. The \emph{structure set} $\calS(X)$ is the set of homotopy equivalences $W\to X$, where $W$ is a closed $4$-manifold, modulo homeomorphism. The distinguished point in both $\calN(X)$ and $\calS(X)$ is given by the identity map $X\to X$. Since a homotopy equivalence is in particular a degree one normal map, we have a map $\calS(X)\to \calN(X)$. The map $\sigma$, called the \emph{surgery obstruction map}, is roughly defined as follows. Given an element $f\colon V\to X$ of $\calN(X)$, we can assume, by performing surgery on circles, that $f$ induces an isomorphism on fundamental groups (this includes checking that the original map and the result of surgery on circles are related via a normal bordism). The kernel of the map $f_*$ on $\pi_2$ is called the \emph{surgery kernel}. The image of $f$ under $\sigma$ is the intersection form on this surgery kernel. Given such a map $f$, the image under $\sigma$ is called the \emph{surgery obstruction} for $f$. Exactness of the surgery sequence at $\calN(X)$ means, in particular, that if the intersection form on the surgery kernel is stably hyperbolic, then the map~$f$ can be replaced, via a normal bordism, by a homotopy equivalence, i.e~an element of~$\calS(X)$.  

In dimensions $\geq 5$ the surgery sequence is exact regardless of fundamental group and applies in both the smooth and the topological settings~\cite{Browder,Novikov,Sullivan,wall-surgery-book,KS}. That the sequence is exact for topological $4$-manifolds with good fundamental group was shown by Freedman and Quinn in~\cite[{Theorem~11.3A}]{FQ} (see also~\cite{Powell-Ray-Teichner:2018-1}). The sphere embedding theorem is a key ingredient -- as in our proof sketch for Theorem B, once we have a degree one normal map $V\to X$ from a $4$-manifold inducing an isomorphism on fundamental groups and with hyperbolic intersection form on the surgery kernel (i.e.~such that the image in $L_4(\Z[\pi_1(X)])$ under $\sigma$ is trivial), one uses the sphere embedding theorem to realises a half-basis of the surgery kernel by pairwise disjoint, locally flat embedded spheres, equipped with geometrically dual generically immersed spheres, and then performs surgery. The result is an element of $\calS(X)$. The surgery sequence for smooth $4$-manifolds is not exact even for trivial fundamental groups, by work of Donaldson~\cite{donaldson-1983}.

There is also a version of the surgery exact sequence for compact $4$-manifolds with nonempty boundary. This is what we could have used in the previous subsection: the target $4$-manifold would have been $X=S^1\times D^3$ with $\pi_1(X)\cong \Z$, the spin null bordism $V$ provides an element of $\calN(X)$, the modified spin null bordism $V'$ is an element of $\calN(X)$ with trivial surgery obstruction, and using the exactness of the surgery sequence, we would have produced the final $4$-manifold $W$ with a homotopy equivalence to $X$, namely an element of $\calS(X)$.  

\subsection{More general results}
It is not too hard to see that the proof of  Theorem B also shows that every knot in an integer homology sphere $Y$ with Alexander polynomial one is topologically slice in the unique, compact, contractible, topological $4$-manifold $C$ with $\partial C=Y$. The existence of such a $C$ was proven by Freedman~\cite{F}. A similar slicing result for knots using surgery theory was proven by Friedl and Teichner in~\cite{friedl-teichner}. Davis showed in~\cite{davis:hopf} that every 2-component link with multi-variable Alexander polynomial one is (topologically) concordant to the Hopf link. 

In a different direction, one can consider the question of existence of locally flat embedded closed surfaces in more general $4$-manifolds. Recall that given a knot~$K$ and integer $n$, the corresponding \emph{knot trace} $X_n(K)$ is built by attaching an~$n$-framed~$2$-handle to $B^4$ along the knot~$K$ in $S^3=\partial B^4$ and then smoothing corners. Note that~$X_n(K)\simeq S^2$ for all $K$ and $n$. A knot is said to be (topologically) \emph{$n$-shake slice} if a generator of $\pi_2(X_n(K))\cong \Z$ can be represented by a locally flat embedded sphere. Of course, every slice knot is $n$-shake slice for all $n$. Surgery theoretic techniques can be used to construct many more $n$-shake slice knots, as in the following theorem. 

\begin{theorem}[{\cite[{Theorem~1.1}]{FMNOPR}}]\label{thm:Zn-shake-slice}
Let $K\subseteq S^3$ be a knot and let $n$ be an integer.
A generator of~$\pi_2(X_n(K))$ can be represented by a locally flat embedded sphere whose complement has
abelian fundamental group if and only if:
  \begin{enumerate}[label=(\roman*)]
    \item\label{item:trivialh1}
    $H_1(S_n^3(K);\Z[\Z/n])=0$; or  equivalently for $n \neq 0$,
     $\prod_{\{\xi  \mid \xi^n=1\} }\Delta_K(\xi) =1$;
    \item\label{item:trivialarf} $\Arf(K)=0$; and
    \item\label{item:signaturesobstruct} $\sigma_K(\xi) =0$ for every $\xi \in S^1$ such that $\xi^n=1$.
     \end{enumerate}
\end{theorem}

Using \Cref{thm:Zn-shake-slice}, one can construct, for each $n\neq 0$, infinitely many topological concordance classes of knots that are $n$-shake slice but neither smoothly $n$-shake slice nor topologically slice~\cite[Corollary~2.3]{FMNOPR}

We have already seen the Alexander polynomial $\Delta_K(t) \in \Z[t,t^{-1}]$. As mentioned before, given a Seifert matrix $V$ for the knot $K$ the Alexander polynomial is defined as $\Delta_K(t)=\det(tV-V^T)$. Also recall that $\Arf(K)\in \Z/2$ is $0$ if $\Delta_K(-1) \equiv \pm 1 \mod{8}$ and is $1$ if $\Delta_K(-1) \equiv \pm 3 \mod{8}$ by~\cite[Proposition~3.4]{levine:polynomial-invariants-codimension-two} and~\cite[Theorem~2]{murasugi:arf-invariant}. Finally the Tristram--Levine signature, for $\xi \in S^1 \subseteq \mathbb{C}$, denoted by $\sigma_K(\xi)$, is the signature of the Hermitian matrix $(1-\xi)V + (1-\ol{\xi})V^T$. For some choices of $n$, there are logical dependencies among the conditions \ref{item:trivialh1}, \ref{item:trivialarf}, and \ref{item:signaturesobstruct} in \Cref{thm:Zn-shake-slice}.  When $n=0$, condition \ref{item:trivialh1} states that $H_1(S^3_0(K); \Z[\Z])=0$, which is equivalent to $\Delta_K(t)=1$, which in turn implies both conditions \ref{item:trivialarf} and \ref{item:signaturesobstruct}. So, in the case $n=0$, \Cref{thm:Zn-shake-slice} coincides with Theorem B. When $n=\pm1$, conditions \ref{item:trivialh1} and \ref{item:signaturesobstruct} are automatically satisfied.

Surgery theoretic techniques can also be used to study \emph{uniqueness} questions. Most famously, this includes the following result of Freedman--Quinn. 

\begin{theorem}[{\cite[{Theorem~11.7A}]{FQ}}]
    Every 2-knot $K\colon S^2\hookrightarrow S^4$ with $\pi_2(S^4\sm K)\cong \Z$ is (topologically) unknotted. 
\end{theorem}

Other results concerning uniqueness up to isotopy for locally flat surfaces in $S^4$ with abelian fundamental group of the complement have been proven using surgery theory, both the classical version alluded to here and the modified theory due to Kreck~\cite{kreck:sad}, in~\cite{conway-powell:pi1-Z,conway-powell-piccirillo:pi1-Z,conway-orson-powell:nonorientable}. 

In \cite{orson-powell:isotopy-knot-traces}, Orson and Powell showed that locally flat embedded spheres representing a generator of the second homotopy group of any given knot trace, with abelian fundamental group of the complement, i.e.~those in \Cref{thm:Zn-shake-slice}, are ambiently isotopic, modulo orientation. One can also consider uniqueness of slice discs using surgery theoretic methods, such as in~\cite{conway-powell-discs,conway:discs-samegroup}. 

Uniqueness up to isotopy of locally flat surfaces in more general $4$-manifolds was considered in~\cite{LW93,HK93,boyer93,conway-powell:pi1-Z,conway-powell-piccirillo:pi1-Z,conway-dai-miller:cp2,conway-orson:cp2}.

\section{Conclusion}
\label{sec:conclusion}

We hope this survey gives the reader a sense of the different flavours of techniques that are used in the topological setting for $4$-manifolds, as well as pointers for what to read next. The disc embedding theorem is the key ingredient in both direct and indirect approaches to finding locally flat embedded surfaces in $4$-manifolds described here.  We refer those interested in more details about the proof of the disc embedding theorem to \cite{FQ,DETbook}. An introduction to surgery theory from a $4$-dimensional consumer's perspective is given in~\cite{DET-book-surgery}. A detailed discussion of open problems regarding the disc embedding theorem can be found in~\cite{DET-book-enigmata}. 

There is a growing number of researchers actively working on topological $4$-manifolds and locally flat surfaces therein. I hope that readers are encouraged to explore not only the landmark achievements in the past, e.g.~\cite{F,FQ,quinn:endsIII}, but also some of the very recent work in this area, such as~\cite{boyle-chen:equivariant,cha:lightbulb,conway-dai-miller:cp2,conway-orson:cp2,conway:discs-samegroup,conway-orson-powell:nonorientable,conway-powell-discs,conway-powell:pi1-Z,conway-powell-piccirillo:pi1-Z,FMNOPR,galvin:casson-sullivan,KPR:gluck,KPRT:sigmet,kasprowski-land:4dgroups,nagy-nicholson-powell:he-she,orson-powell:isotopy-knot-traces,pencovitch:non-orientable}.

\section{Exercises}
\label{sec:exercises}

The upcoming problems are separated into three levels. The introductory problems should be attempted if you are seeing all of this material for the first time. Prerequisites are courses in introductory geometric and algebraic topology. The moderate problems are for readers who are already comfortable with some of the terminology; they may require nontrivial external input, which we have tried to indicate as hints. The section ends with a list of challenge problems for advanced readers. 

\subsection{Introductory problems}
\begin{exercise}\label{ex:emb-not-loc-flat}
    Give an example of a surface in a $4$-manifold which is topologically embedded (i.e.~there is a continuous map $f\colon\Sigma_g\hookrightarrow M$ where $\Sigma_g$ is some closed surface, $M$ is some $4$-manifold, and $f$ is a homeomorphism onto its image), but not locally flatly embedded.

    \emph{Hint:} Given a knot $K\subseteq S^3$, consider the disc given by $\cone(K)\subseteq \cone(S^3)= B^4$. When is this disc locally flat? Recall from classical knot theory that a knot $K$ is the unknot if and only if $\pi_1(S^3\sm K)\cong \Z$.     
\end{exercise}

\begin{exercise}\label{ex:smooth-implies-loc-flat}
    Convince yourself that every smooth embedding of a surface in a smooth $4$-manifold is locally flat. Remind yourself of the smooth analogues of \Cref{thm:normal-bundles-transversality,thm:immersion-lemma} and the ideas of their proofs. Without going into the details, consider why those proofs fail in the purely topological setting.
\end{exercise}

\begin{exercise}\label{ex:hopf-link-at-intersection}
    Consider $\R^4$, given by $x,y,z,t$ coordinates, as in \Cref{fig:transverse-movie-1}. Let $B$ denote the $4$-ball of unit radius at the origin. 
    
    \begin{enumerate}
        \item Show that the intersection of the $xy$- and $zt$-planes with $\partial B=S^3$ is a Hopf link. 
        \item Give the $xy$- and $zt$-planes, as well as $\R^4$, the positive orientation. Then $B$ inherits an orientation from $\R^4$. Orient $S^3$ as the boundary of $B$. Which of the two possible (oriented) Hopf links is obtained in part 1?
    \end{enumerate}
    Now suppose that two generically immersed surfaces $f$ and $g$ in a $4$-manifold $M$ intersect transversely at a point $p\in M$. Let $C\subseteq M$ be a small $4$-ball at $p$. 
        \begin{enumerate}
        \item[(3)] Conclude by part 2 that one can choose $C$ to be small enough so that $\partial C\cap (f\cup g)$ is a Hopf link in $\partial C$.
        \item[(4)] Suppose that $M$, $f$, and $g$ are all oriented. How does the sign of the intersection point $p$ determine the orientation of the Hopf link in part 3?
    \end{enumerate}
\end{exercise}

\begin{exercise}
    Draw the Clifford torus at the transverse intersection point shown in \Cref{fig:transverse-movie-2}.
\end{exercise}

\begin{exercise}\label{ex:pi1-negligible}
    Let $\Sigma$ be a surface and let $M$ be a $4$-manifold. Suppose we have a generic immersion $f\colon \Sigma\looparrowright M$ with a \emph{geometrically dual} sphere, i.e.~there is some~$g\colon S^2\looparrowright M$ such that $f\pitchfork g$ is a single point. Show that the inclusion~${\iota\colon M\sm \nu f\to M}$ induces an isomorphism 
    \begin{equation}\label{eq:pi1-negligible}
    \pi_1(M\sm \nu f)\xrightarrow[\iota_*]{\cong} \pi_1(M), 
    \end{equation}
    where $\nu f$ is the normal vector bundle of $f$. A generic immersion $f$ satisfying \eqref{eq:pi1-negligible} is said to be \emph{$\pi_1$-negligible}.
\end{exercise}


\begin{exercise}\label{ex:surgery-on-circles}
    Let $C\colon S^1\hookrightarrow M$ be an embedded, orientation preserving loop in a $4$-manifold $M$. The procedure of \emph{surgery on $M$ along $C$} is as follows. Choose a tubular neighbourhood of $C$, call it $\nu C\homeo S^1\times D^3$. Cut out the interior $\mathring{\nu}C$, and glue in $D^2\times S^2$, via the identity map along the boundary $S^1\times S^2$. There are two possible identifications of $\partial \nu C$ with $S^1\times S^2$, and therefore there are two possible gluing maps. 
    
    Suppose we have a map $X\to Y$ of $4$-manifolds, such that the induced map on fundamental groups is a surjection. Use surgery on circles in $X$ to change $X$ to some~$X'$ with a map to $Y$ inducing an isomorphism on fundamental groups. 
\end{exercise}

\subsection{Moderate problems}

\begin{exercise}\label{ex:interior-twisting-sign}
    Let $f\colon S^2\looparrowright M$ be a generic immersion in an oriented $4$-manifold~$M$. Choose an orientation on $f$. Determine the sign of the intersection point created in $f$ by the procedure described in \Cref{fig:interior-twisting}. Does the sign depend on the original orientation of $f$?
\end{exercise}

\begin{figure}[tbp]
    \centering
    \begin{tikzpicture}
        \node[anchor=south west,inner sep=0] at (0,0){\includegraphics[width=6.125cm]{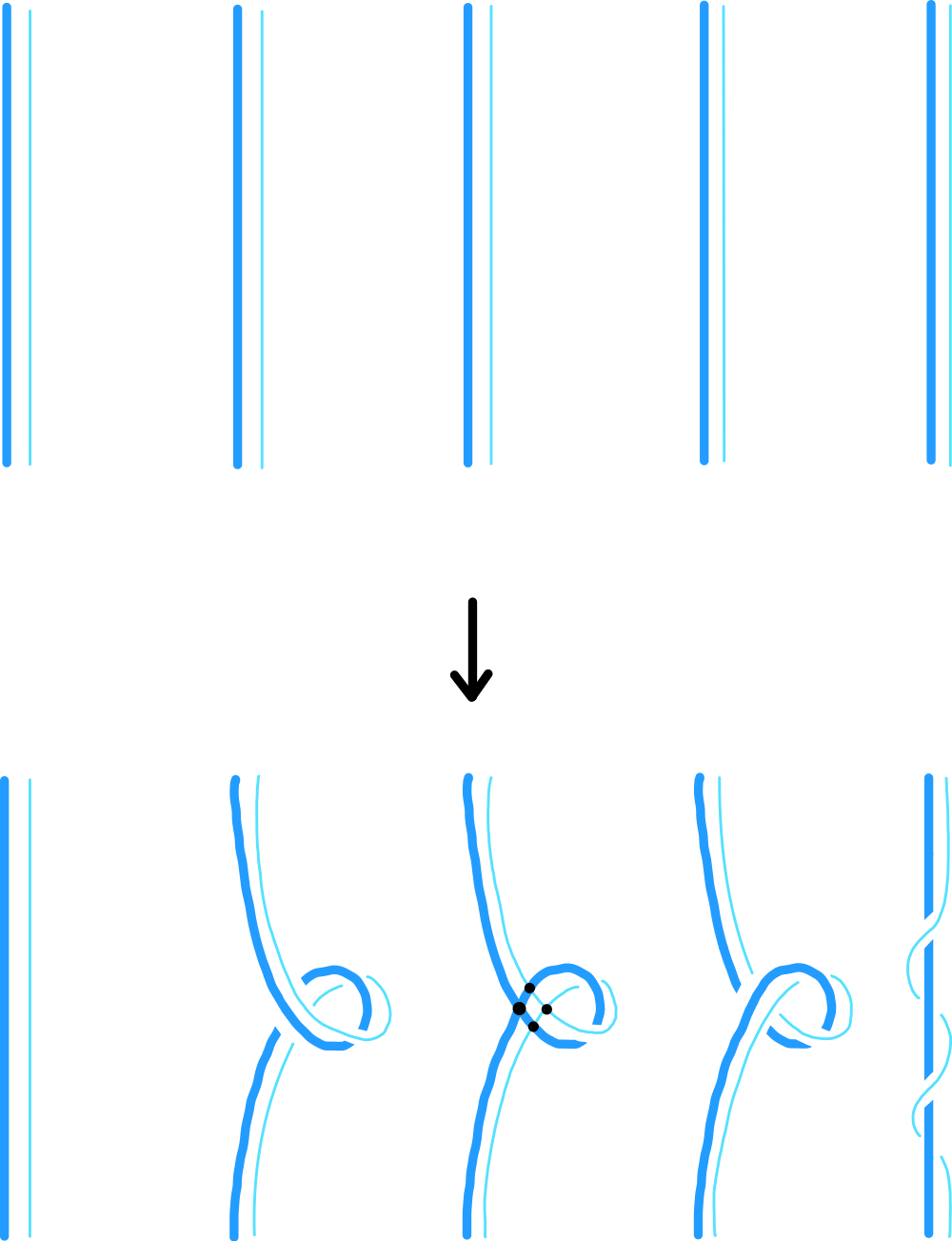}};
           \node at (0.1,-0.25) {$t=-\varepsilon$};
            \node at (1.6,-0.25) {$t=-\tfrac{\varepsilon}{2}$};
            \node at (3.05,-0.25) {$t=0$};
            \node at (4.55,-0.25) {$t=\tfrac{\varepsilon}{2}$};
            \node at (6.1,-0.25) {$t=\varepsilon$};
            \node at (0.1,4.75) {$t=-\varepsilon$};
            \node at (1.6,4.75) {$t=-\tfrac{\varepsilon}{2}$};
            \node at (3.05,4.75) {$t=0$};
            \node at (4.55,4.75) {$t=\tfrac{\varepsilon}{2}$};
            \node at (6.1,4.75) {$t=\varepsilon$};
	\end{tikzpicture} 
    \caption{Interior twisting changes the euler number of the normal vector bundle by $\mp2$.}
    \label{fig:interior-twisting-euler}
\end{figure}

\begin{exercise}\label{ex:twisting-euler}
    Prove the following statements. 
    \begin{enumerate}
        \item Let $f\colon S^2\looparrowright M$ be a generically immersed sphere in a $4$-manifold. By interior twisting, we can insert a double point in $f$ with sign $\pm 1$. Show that this changes the euler number of the normal vector bundle of $f$ by $\mp 2$. 
        \item Let $W$ be a generically immersed Whitney disc pairing intersections between generically immersed spheres $f,g\colon S^2\to M$ in a $4$-manifold $M$. We can do a boundary twist of $W$ about either $f$ or $g$ to introduce a new double point between $\mathring{W}$ and $f$ or $g$ respectively. Show that this changes $\tw(\partial W)$ by $\pm 1$.
    \end{enumerate}
    \emph{Hint:} For the case of interior twisting, consider and interpret \Cref{fig:interior-twisting-euler}. Draw a similar picture for boundary twisting. 
\end{exercise}

\begin{exercise}\label{ex:spin-primitive-sphere}
    Let $M$ be a closed, simply connected, spin $4$-manifold. Show that every primitive class in $H_2(M;\Z)$ can be represented by a locally flat, embedded sphere. Can it always be represented by a smoothly embedded sphere? 
\end{exercise}

\begin{exercise}\label{ex:slice-disc-complement-zero-surgery}
    Let $K\subseteq S^3$ be a knot bounding a topological slice disc~${\Delta\subseteq B^4}$. Let $\mathring{\nu} \Delta$ denote an open tubular neighbourhood of $\Delta$. Show that $\partial (B^4\sm \mathring{\nu} \Delta)$ is homeomorphic to $S^3_0(K)$, the result of 0-framed Dehn surgery on $S^3$ along $K$.  
\end{exercise}

\begin{exercise}\label{ex:lambda-mu-basic-properties}
    Consider the equivariant intersection and self-intersection numbers defined in~\Cref{sec:lambda-mu}.
    \begin{enumerate}
    \item What is the effect on $\lambda(f,g)$ of 
        \begin{itemize}
        \item changing the paths $\alpha_f^p$ and $\alpha_g^p$? 
        \item changing the whiskers $w_f$ and $w_g$?
        \item changing the basepoint $m$? (How might you get new whiskers?)
        \end{itemize}
    \item What is the effect on $\mu(f)$ of 
        \begin{itemize}
        \item changing the paths $\alpha_1^p$ and $\alpha_2^p$? 
        \item changing the whisker $w_f$?
        \item changing the basepoint $m$? (How might you get a new whisker?)
        \end{itemize}  
    \item Conclude from the above two parts that there is a well defined \emph{equivariant intersection number} 
    \begin{align}
    \lambda\colon \pi_2(M)\times \pi_2(M) \longrightarrow & \Z[\pi_1(M)]\\
    (f,g)\longmapsto    &\lambda(f,g)
    \end{align}
    \item\label{item:mu-domain-codomain} As above, try to define the self-intersection number. What should be the domain and codomain? \emph{Hint:} Was there any ambiguity in the definition of $\mu(f)$? Can we change the value of $\mu(f)$ by changing $f$ by a homotopy? (Which homotopies of surfaces in a $4$-manifold have we seen in this paper?) Recall that, generically, a homotopy of a surface $f$ in a $4$-manifold is some sequence of isotopies, interior twisting, finger moves, and Whitney moves along untwisted, embedded, and disjoint Whitney discs (\Cref{thm:generic-immersions-bijection}).
    \end{enumerate}
\end{exercise}

\begin{exercise}\label{ex:lambda-mu-zero-whitney}
    Let $f$ and $g$ be generically immersed spheres in some connected, oriented $4$-manifold $M$. Assume we have chosen a basepoint in $M$ and whiskers for~$f$ and $g$. Show the following. 
    \begin{enumerate}
        \item the intersection number $\lambda(f,g)=0$ if and only if all the intersections of $f$ and $g$ can be paired up by untwisted generically immersed Whitney discs in $M$, with disjointly embedded boundaries. 
        \item the self-intersection number $\mu(f)=0$ if and only if the self-intersections of $f$ can be paired up by untwisted generically immersed Whitney discs in $M$, with disjointly embedded boundaries. 
    \end{enumerate}
\end{exercise}

\begin{exercise}\label{ex:equiv-int-form-conn-sum}
    Let $M$ and $N$ be $4$-manifolds. 
    \begin{enumerate}
        \item Express the integral intersection form of $M\#N$ in terms of those of $M$ and $N$. 
        \item Suppose that $M=S^1\times S^3$ and $N$ is simply connected. What is the equivariant intersection form of $M\#N$? \emph{Hint: Consider the universal cover of $M\#N$.}
        \item More generally, suppose that $\pi_1(M)\cong \Z$ and $N$ is simply connected. What is the equivariant intersection form of $N$?
    \end{enumerate}
    For a further challenge, consider the intersection form of $M\# N$ where $\pi_1(M)$ and $\pi_1(N)$ are arbitrary. 
\end{exercise}
\begin{exercise}\label{ex:geometric-dual-pi1}
    Let $M$ be a simply connected $4$-manifold, and let $S\subseteq M$ be an embedded 2-sphere with trivial normal vector bundle. Let $M'$ denote the result of surgery on $M$ along $S$. 
    \begin{enumerate}
    \item What can you say about the fundamental group of $M'$? 
    \item Can you think of a condition on $S$ to ensure that $M'$ is simply connected? 
    \item Find an example of $S$ and $M$ such that $M'$ is simply connected.
    \item Find an example of $S$ and $M$ such that $M'$ is not simply connected. 
    \item Find an example of $S$ and $M$ such that $\pi_1(M')$ is nontrivial but $H_1(M';\Z)$ is trivial.
    \end{enumerate}
\end{exercise}


\subsection{Challenge problems}
\begin{figure}[tbp]
    \centering
    \begin{tikzpicture}
        \node[anchor=south west,inner sep=0] at (0,0){\includegraphics[width=11cm]{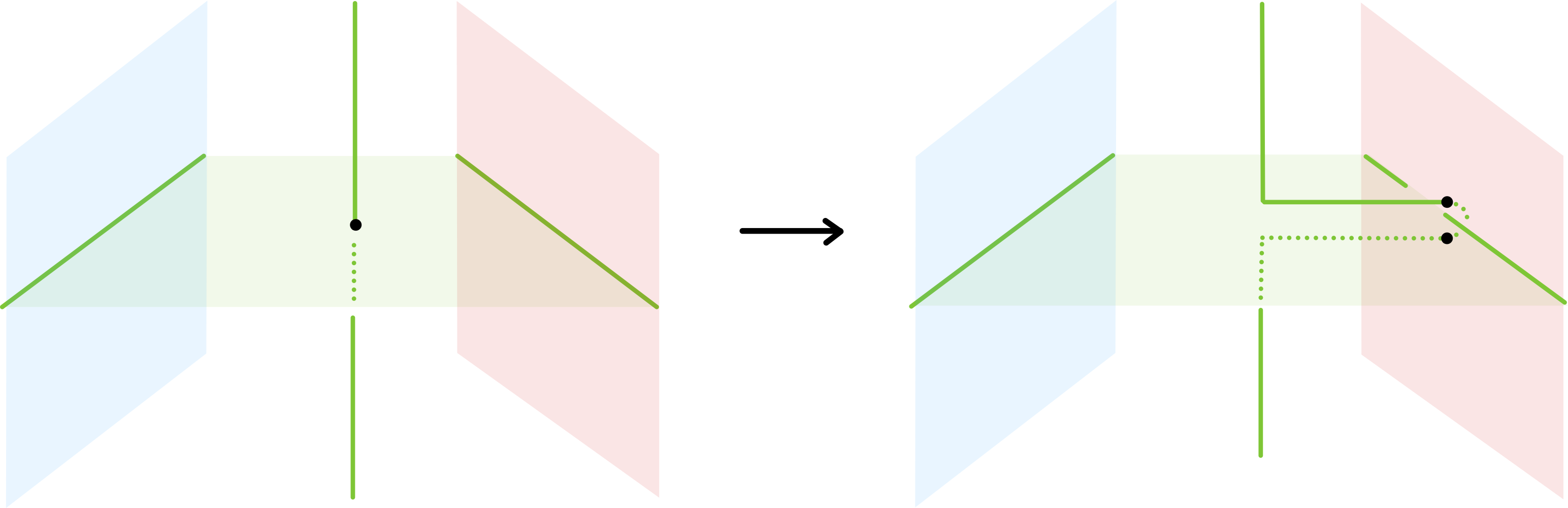}};
            \node at (0.2,0.4) {$f$};
            \node at (4.5,0.4) {$g$};
            \node at (2,2) {$W$};
            \node at (2.3,0.4) {$W$};            \node at (6.6,0.4) {$f$};
            \node at (10.9,0.4) {$g$};
            \node at (8.4,2) {$W$};
            \node at (8.7,0.4) {$W$};
	\end{tikzpicture} 
    \caption{A portion of a Whitney disc $W$ (green), pairing intersections between $f$ (blue) and $g$ (red). A self-intersection of $W$ is pictured on the left. A finger move on~$W$ in the direction of $g$ removes the self-intersection at the expense of creating two new, algebraically cancelling intersection points between $W$ and $g$.}
    \label{fig:geom-casson-hint}
\end{figure}

\begin{exercise}\label{ex:geometric-Casson-lemma}
    Prove the \emph{geometric Casson lemma}: Let~$f$ and $g$ be transversely intersecting generic immersions of compact surfaces in a connected $4$-manifold $M$. Assume that the intersection points $\{p,q\}\subseteq f\pitchfork g$ are paired by a generically immersed Whitney disc $W$. Then there is a regular homotopy from $f \cup g$ to $\ol f \cup \ol{g}$ such that  $\ol f \pitchfork \ol{g} = (f\pitchfork g) \smallsetminus \{p,q\}$, that is, the two paired intersections have been removed.
    
    The regular homotopy may create many new self-intersections of $f$ and $g$; however, these are algebraically cancelling. Moreover, the regular homotopy is supported in a small neighbourhood of $W$.

    \emph{Hint:} Recall the definition of a regular homotopy (\Cref{def:reg-htpy}). Draw a schematic picture of the Whitney disc, showing all the possible problems, as shown in \Cref{tab:table}. There are no problems of type 2 (why?) Use the manoeuvres from the proof of Theorem A to arrange that $W$ only has problems of type 1 and 4. For intersections of $W$ with itself, perform a finger move on $W$ towards its Whitney arcs, creating new intersections with $f$ or $g$, as shown in \Cref{fig:geom-casson-hint}. Try to do a similar move for intersections of $f$ and $g$ with $W$.
\end{exercise}

\begin{exercise}
    Let $\Sigma$ be a surface and let $M$ be a $4$-manifold. Let $f\colon \Sigma\looparrowright M$ be a generic immersion. Suppose a pair of double points of $f$ with opposite sign are paired by an untwisted generically immersed Whitney disc $W$. Show that the \emph{immersed Whitney move} on $f$ along $W$ is a regular homotopy. In other words, show that it can be expressed as a concatenation of isotopies, finger moves, and Whitney moves along untwisted, embedded, disjoint Whitney discs, with interiors disjoint from $f$. 
\end{exercise}

\begin{exercise}\label{ex:arf-via-whitney}
    Let $K\subseteq S^3$ be a knot and let $\Delta\colon D^2\looparrowright B^4$ be a generically immersed disc bounded by $K$. Suppose that the signed count of self-intersections of $\Delta$ is trivial. By Exercise~\ref{ex:lambda-mu-zero-whitney}, the double points of $f$ can be paired up by untwisted generically immersed Whitney discs $\{W_i\}$ in $M$, with disjointly embedded boundaries. Assume that $\{W_i\}$ only meets $\Delta$ transversely in the interiors, except at the Whitney circles. Show that 
    \[
        \Arf(K) \equiv \sum_i \left\lvert \mathring{W}_i\pitchfork \Delta \right\rvert \mod{2}.
    \]
    
    If we do not assume that the Whitney discs are untwisted, or that they have disjoint, embedded boundaries, how would the count on the right hand side need to be changed? 
\end{exercise}

\begin{exercise}\label{ex:sphere-emb-thm}
    Prove the sphere embedding theorem (\Cref{thm:sphere-emb-thm}).     
\end{exercise}

\begin{exercise}\label{ex:0-surgery-char-sliceness}
    Prove the 0-surgery characterisation of topological sliceness (\Cref{thm:0-surgery-char-sliceness}). 
    
    \emph{Hint:} You might need the following form of the (proven) $4$-dimensional topological Poincar\'{e} conjecture: a homotopy $4$-ball with boundary $S^3$ is homeomorphic to $B^4$. 
\end{exercise}

\bibliographystyle{alpha}
\bibliography{ray}

@article {quinn:endsIII,
    AUTHOR = {Quinn, Frank},
     TITLE = {Ends of maps. {III}. {D}imensions {$4$} and {$5$}},
   JOURNAL = {J. Differential Geometry},
  FJOURNAL = {Journal of Differential Geometry},
    VOLUME = {17},
      YEAR = {1982},
    NUMBER = {3},
     PAGES = {503--521},
      ISSN = {0022-040X,1945-743X},
   MRCLASS = {57N13},
  MRNUMBER = {679069},
MRREVIEWER = {R.\ C.\ Kirby},
       URL = {http://projecteuclid.org/euclid.jdg/1214437139},
}

@MISC {lee:MSE-propervsproper,
    TITLE = {Definitions of a proper embedding},
    AUTHOR = {Lee, John M.\},
    HOWPUBLISHED = {Mathematics Stack Exchange},
    NOTE = {URL:https://math.stackexchange.com/q/2362048 (version: 2017-07-18)},
    EPRINT = {https://math.stackexchange.com/q/2362048},
    URL = {https://math.stackexchange.com/q/2362048}
}

@misc{ranicki:orsay-slides,
  title = {High dimensional manifold topology then and now},
  howpublished = {\url{https://www.maths.ed.ac.uk/~v1ranick/slides/orsay.pdf}},
  note = {Slides from 3 lectures at Larry Siebenmann's retirement conference, Orsay, December 2005. Accessed: 2025-04-18},
            AUTHOR={Ranicki, Andrew}
}

@article {donaldson:h-cobordism,
    AUTHOR = {Donaldson, S. K.},
     TITLE = {Irrationality and the {$h$}-cobordism conjecture},
   JOURNAL = {J. Differential Geom.},
  FJOURNAL = {Journal of Differential Geometry},
    VOLUME = {26},
      YEAR = {1987},
    NUMBER = {1},
     PAGES = {141--168},
      ISSN = {0022-040X,1945-743X},
   MRCLASS = {57R80 (32C10 32G13 32J15 57N13 58E15 58G30)},
  MRNUMBER = {892034},
MRREVIEWER = {David\ R.\ Morrison},
       URL = {http://projecteuclid.org/euclid.jdg/1214441179},
}

@article {friedman-morgan:algebraicII,
    AUTHOR = {Friedman, Robert and Morgan, John W.},
     TITLE = {On the diffeomorphism types of certain algebraic surfaces.
              {II}},
   JOURNAL = {J. Differential Geom.},
  FJOURNAL = {Journal of Differential Geometry},
    VOLUME = {27},
      YEAR = {1988},
    NUMBER = {3},
     PAGES = {371--398},
      ISSN = {0022-040X,1945-743X},
   MRCLASS = {57R55 (14F05 14F45 14J27 32G13 32J15)},
  MRNUMBER = {940111},
MRREVIEWER = {I.\ Dolgachev},
       URL = {http://projecteuclid.org/euclid.jdg/1214442001},
}

@article {friedman-morgan:algebraicI,
    AUTHOR = {Friedman, Robert and Morgan, John W.},
     TITLE = {On the diffeomorphism types of certain algebraic surfaces.
              {I}},
   JOURNAL = {J. Differential Geom.},
  FJOURNAL = {Journal of Differential Geometry},
    VOLUME = {27},
      YEAR = {1988},
    NUMBER = {2},
     PAGES = {297--369},
      ISSN = {0022-040X,1945-743X},
   MRCLASS = {57R55 (14F05 14F45 14J27 32G13 32J15)},
  MRNUMBER = {925124},
MRREVIEWER = {I.\ Dolgachev},
       URL = {http://projecteuclid.org/euclid.jdg/1214441784},
}

@article {okonek-vandeven:infinite-log-transform,
    AUTHOR = {Okonek, C. and Van de Ven, A.},
     TITLE = {Stable bundles and differentiable structures on certain
              elliptic surfaces},
   JOURNAL = {Invent. Math.},
  FJOURNAL = {Inventiones Mathematicae},
    VOLUME = {86},
      YEAR = {1986},
    NUMBER = {2},
     PAGES = {357--370},
      ISSN = {0020-9910,1432-1297},
   MRCLASS = {32G13 (14F05 14J27 32J25 32L99 57R55)},
  MRNUMBER = {856849},
MRREVIEWER = {R.\ R.\ Simha},
       DOI = {10.1007/BF01389075},
       URL = {https://doi.org/10.1007/BF01389075},
}

@article {fintushel-stern:knot-surgery,
    AUTHOR = {Fintushel, Ronald and Stern, Ronald J.},
     TITLE = {Knots, links, and {$4$}-manifolds},
   JOURNAL = {Invent. Math.},
  FJOURNAL = {Inventiones Mathematicae},
    VOLUME = {134},
      YEAR = {1998},
    NUMBER = {2},
     PAGES = {363--400},
      ISSN = {0020-9910,1432-1297},
   MRCLASS = {57R57 (57M25 57N13 57R55)},
  MRNUMBER = {1650308},
MRREVIEWER = {Paolo\ Lisca},
       DOI = {10.1007/s002220050268},
       URL = {https://doi.org/10.1007/s002220050268},
}

@article {levine:polynomial-invariants-codimension-two,
    AUTHOR = {Levine, J.},
     TITLE = {Polynomial invariants of knots of codimension two},
   JOURNAL = {Ann. of Math. (2)},
  FJOURNAL = {Annals of Mathematics. Second Series},
    VOLUME = {84},
      YEAR = {1966},
     PAGES = {537--554},
      ISSN = {0003-486X},
   MRCLASS = {55.20},
  MRNUMBER = {200922},
MRREVIEWER = {R.\ H.\ Fox},
       DOI = {10.2307/1970459},
       URL = {https://doi.org/10.2307/1970459},
}

@article {murasugi:arf-invariant,
    AUTHOR = {Murasugi, Kunio},
     TITLE = {The {A}rf invariant for knot types},
   JOURNAL = {Proc. Amer. Math. Soc.},
  FJOURNAL = {Proceedings of the American Mathematical Society},
    VOLUME = {21},
      YEAR = {1969},
     PAGES = {69--72},
      ISSN = {0002-9939,1088-6826},
   MRCLASS = {55.20},
  MRNUMBER = {238301},
MRREVIEWER = {L.\ Neuwirth},
       DOI = {10.2307/2036859},
       URL = {https://doi.org/10.2307/2036859},
}

@article {F,
    AUTHOR = {Freedman, Michael Hartley},
     TITLE = {The topology of four-dimensional manifolds},
   JOURNAL = {J. Differential Geometry},
  FJOURNAL = {Journal of Differential Geometry},
    VOLUME = {17},
      YEAR = {1982},
    NUMBER = {3},
     PAGES = {357--453},
      ISSN = {0022-040X,1945-743X},
   MRCLASS = {57N12 (57R80 57R99)},
  MRNUMBER = {679066},
MRREVIEWER = {John\ J.\ Walsh},
       URL = {http://projecteuclid.org/euclid.jdg/1214437136},
}

@article {LW97,
    AUTHOR = {Lee, Ronnie and Wilczy\'{n}ski, Dariusz M.},
     TITLE = {Representing homology classes by locally flat surfaces of
              minimum genus},
   JOURNAL = {Amer. J. Math.},
  FJOURNAL = {American Journal of Mathematics},
    VOLUME = {119},
      YEAR = {1997},
    NUMBER = {5},
     PAGES = {1119--1137},
      ISSN = {0002-9327,1080-6377},
   MRCLASS = {57N35 (57N13)},
  MRNUMBER = {1473071},
MRREVIEWER = {Stefan\ A.\ Bauer},
       URL =
              {http://muse.jhu.edu/journals/american_journal_of_mathematics/v119/119.5lee.pdf},
}

@ARTICLE{KPRT:sigmet,
    AUTHOR = {Kasprowski, Daniel and Powell, Mark and Ray, Arunima and
              Teichner, Peter},
     TITLE = {Embedding surfaces in 4-manifolds},
   JOURNAL = {Geom. Topol.},
  FJOURNAL = {Geometry \& Topology},
    VOLUME = {28},
      YEAR = {2024},
    NUMBER = {5},
     PAGES = {2399--2482},
      ISSN = {1465-3060,1364-0380},
   MRCLASS = {57K40 (57N35)},
  MRNUMBER = {4793644},
MRREVIEWER = {Allison\ N.\ Miller},
       DOI = {10.2140/gt.2024.28.2399},
       URL = {https://doi.org/10.2140/gt.2024.28.2399},
}

@incollection {DET-book-goodgroups,
    AUTHOR = {Kim, Min Hoon and Orson, Patrick and Park, JungHwan and Ray, Arunima},
    TITLE = {Good groups},
    BOOKTITLE = {The disc embedding theorem},
    PUBLISHER = {Oxford University Press},
    YEAR = {2021},
   label={KOPR}
}

@incollection {DET-book-DETintro,
    AUTHOR = {Powell, Mark and Ray, Arunima},
     TITLE = {Intersection numbers and the statement of the disc embedding
              theorem},
 BOOKTITLE = {The disc embedding theorem},
     PAGES = {155--170},
 PUBLISHER = {Oxford Univ. Press, Oxford},
      YEAR = {2021},
      ISBN = {978-0-19-884131-9},
   MRCLASS = {57K40 (57R40)},
  MRNUMBER = {4519509},
}

@incollection {DET-book-enigmata,
    AUTHOR = {Kim, Min Hoon and Orson, Patrick and Park, JungHwan and Ray, Arunima},
    TITLE = {Open problems},
    BOOKTITLE = {The disc embedding theorem},
    PUBLISHER = {Oxford University Press},
    YEAR = {2021},
   label={KOPR}
}

@article{FMNOPR, 
    AUTHOR = {Feller, Peter and Miller, Allison N. and Nagel, Matthias and
              Orson, Patrick and Powell, Mark and Ray, Arunima},
     TITLE = {Embedding spheres in knot traces},
   JOURNAL = {Compos. Math.},
  FJOURNAL = {Compositio Mathematica},
    VOLUME = {157},
      YEAR = {2021},
    NUMBER = {10},
     PAGES = {2242--2279},
      ISSN = {0010-437X,1570-5846},
   MRCLASS = {57K40 (57K10 57N35 57N70 57R67)},
  MRNUMBER = {4328000},
       DOI = {10.1112/S0010437X21007508},
       URL = {https://doi.org/10.1112/S0010437X21007508},
	label={FMNOPR}
}

@article {conway-powell-discs,
    AUTHOR = {Conway, Anthony and Powell, Mark},
     TITLE = {Characterisation of homotopy ribbon discs},
   JOURNAL = {Adv. Math.},
  FJOURNAL = {Advances in Mathematics},
    VOLUME = {391},
      YEAR = {2021},
     PAGES = {Paper No. 107960, 29},
      ISSN = {0001-8708,1090-2082},
   MRCLASS = {57K10 (57K40 57N35)},
  MRNUMBER = {4300918},
MRREVIEWER = {Stephan\ Rosebrock},
       DOI = {10.1016/j.aim.2021.107960},
       URL = {https://doi.org/10.1016/j.aim.2021.107960},
}

@ARTICLE{conway-orson:cp2,
    AUTHOR = {Conway, Anthony and Orson, Patrick},
     TITLE = {Locally flat simple spheres in {$\mathbb{C}P^2$}},
   JOURNAL = {Bull. Lond. Math. Soc.},
  FJOURNAL = {Bulletin of the London Mathematical Society},
    VOLUME = {57},
      YEAR = {2025},
    NUMBER = {1},
     PAGES = {150--163},
      ISSN = {0024-6093,1469-2120},
   MRCLASS = {57K45 (57N35 57R67)},
  MRNUMBER = {4849496},
       DOI = {10.1112/blms.13188},
       URL = {https://doi.org/10.1112/blms.13188},
}

@article {davis:hopf,
    AUTHOR = {Davis, James F.},
     TITLE = {A two component link with {A}lexander polynomial one is
              concordant to the {H}opf link},
   JOURNAL = {Math. Proc. Cambridge Philos. Soc.},
  FJOURNAL = {Mathematical Proceedings of the Cambridge Philosophical
              Society},
    VOLUME = {140},
      YEAR = {2006},
    NUMBER = {2},
     PAGES = {265--268},
      ISSN = {0305-0041,1469-8064},
   MRCLASS = {57M25 (57M27)},
  MRNUMBER = {2212279},
MRREVIEWER = {Richard\ John\ Hadji},
       DOI = {10.1017/S0305004105008856},
       URL = {https://doi.org/10.1017/S0305004105008856},
}

@ARTICLE{nagy-nicholson-powell:he-she,
       author = {{Nagy}, Csaba and {Nicholson}, John and {Powell}, Mark},
        title = "{Simple homotopy types of even dimensional manifolds}",
      journal = {arXiv e-prints},
         year = 2023,
        month = nov,
          eid = {arXiv:2312.00322},
        pages = {arXiv:2312.00322},
          doi = {10.48550/arXiv.2312.00322},
archivePrefix = {arXiv},
       eprint = {2312.00322},
 primaryClass = {math.AT},
       adsurl = {https://ui.adsabs.harvard.edu/abs/2023arXiv231200322N}
}

@ARTICLE{cha:lightbulb,
       author = {{Cha}, Jae Choon and {Kim}, Byeorhi},
        title = "{Light bulb smoothing for topological surfaces in 4-manifolds}",
      journal = {arXiv e-prints},
         year = 2023,
        month = mar,
          eid = {arXiv:2303.12857},
        pages = {arXiv:2303.12857},
          doi = {10.48550/arXiv.2303.12857},
archivePrefix = {arXiv},
       eprint = {2303.12857},
 primaryClass = {math.GT},
       adsurl = {https://ui.adsabs.harvard.edu/abs/2023arXiv230312857C}
}

@article {kasprowski-land:4dgroups,
    AUTHOR = {Kasprowski, Daniel and Land, Markus},
     TITLE = {Topological 4-manifolds with 4-dimensional fundamental group},
   JOURNAL = {Glasg. Math. J.},
  FJOURNAL = {Glasgow Mathematical Journal},
    VOLUME = {64},
      YEAR = {2022},
    NUMBER = {2},
     PAGES = {454--461},
      ISSN = {0017-0895,1469-509X},
   MRCLASS = {57K40 (57N65)},
  MRNUMBER = {4404107},
MRREVIEWER = {Richard\ Stong},
       DOI = {10.1017/S0017089521000215},
       URL = {https://doi.org/10.1017/S0017089521000215},
}

@article {pencovitch:non-orientable,
    AUTHOR = {Pencovitch, Mark},
     TITLE = {Unknotting nonorientable surfaces of genus 4 and 5},
   JOURNAL = {Linear Algebra Appl.},
  FJOURNAL = {Linear Algebra and its Applications},
    VOLUME = {702},
      YEAR = {2024},
     PAGES = {195--217},
      ISSN = {0024-3795,1873-1856},
   MRCLASS = {57K40 (57K45 57N35)},
  MRNUMBER = {4792484},
       DOI = {10.1016/j.laa.2024.08.014},
       URL = {https://doi.org/10.1016/j.laa.2024.08.014},
}

@ARTICLE{galvin:casson-sullivan,
       author = {{Galvin}, Daniel A.~P.},
        title = "{The Casson-Sullivan invariant for homeomorphisms of 4-manifolds}",
      journal = {arXiv e-prints},
         year = 2024,
        month = may,
          eid = {arXiv:2405.07928},
        pages = {arXiv:2405.07928},
          doi = {10.48550/arXiv.2405.07928},
archivePrefix = {arXiv},
       eprint = {2405.07928},
 primaryClass = {math.GT},
       adsurl = {https://ui.adsabs.harvard.edu/abs/2024arXiv240507928G}
}

@article {KPR:gluck,
    AUTHOR = {Kasprowski, Daniel and Powell, Mark and Ray, Arunima},
     TITLE = {Gluck twists on concordant or homotopic spheres},
   JOURNAL = {Math. Res. Lett.},
  FJOURNAL = {Mathematical Research Letters},
    VOLUME = {30},
      YEAR = {2023},
    NUMBER = {6},
     PAGES = {1787--1811},
      ISSN = {1073-2780,1945-001X},
   MRCLASS = {57K40 (57N70)},
  MRNUMBER = {4779154},
}

@ARTICLE{conway-orson-powell:nonorientable,
       author = {{Conway}, Anthony and {Orson}, Patrick and {Powell}, Mark},
        title = "{Unknotting nonorientable surfaces}",
      journal = {arXiv e-prints (to appear: J.~Eur.~Math.~Soc.)},
         year = 2023,
        month = jun,
          eid = {arXiv:2306.12305},
        pages = {arXiv:2306.12305},
          doi = {10.48550/arXiv.2306.12305},
archivePrefix = {arXiv},
       eprint = {2306.12305},
 primaryClass = {math.GT},
       adsurl = {https://ui.adsabs.harvard.edu/abs/2023arXiv230612305C}
}

@article {kreck:sad,
    AUTHOR = {Kreck, Matthias},
     TITLE = {Surgery and duality},
   JOURNAL = {Ann. of Math. (2)},
  FJOURNAL = {Annals of Mathematics. Second Series},
    VOLUME = {149},
      YEAR = {1999},
    NUMBER = {3},
     PAGES = {707--754},
      ISSN = {0003-486X,1939-8980},
   MRCLASS = {57R67 (57R65)},
  MRNUMBER = {1709301},
MRREVIEWER = {Laurence\ R.\ Taylor},
       DOI = {10.2307/121071},
       URL = {https://doi.org/10.2307/121071},
}

@article {LW93,
    AUTHOR = {Lee, Ronnie and Wilczy\'nski, Dariusz M.},
     TITLE = {Representing homology classes by locally flat {$2$}-spheres},
   JOURNAL = {$K$-Theory},
  FJOURNAL = {$K$-Theory. An Interdisciplinary Journal for the Development,
              Application, and Influence of $K$-Theory in the Mathematical
              Sciences},
    VOLUME = {7},
      YEAR = {1993},
    NUMBER = {4},
     PAGES = {333--367},
      ISSN = {0920-3036,1573-0514},
   MRCLASS = {57N13 (57N35)},
  MRNUMBER = {1246281},
MRREVIEWER = {Ian\ Hambleton},
       DOI = {10.1007/BF00962053},
       URL = {https://doi.org/10.1007/BF00962053},
}

@ARTICLE{conway-powell-piccirillo:pi1-Z,
    AUTHOR = {Conway, Anthony and Piccirillo, Lisa and Powell, Mark},
     TITLE = {{$4$}-manifolds with boundary and fundamental group
              {$\mathbb{Z}$}},
   JOURNAL = {Comment. Math. Helv.},
  FJOURNAL = {Commentarii Mathematici Helvetici. A Journal of the Swiss
              Mathematical Society},
    VOLUME = {100},
      YEAR = {2025},
    NUMBER = {2},
     PAGES = {323--420},
      ISSN = {0010-2571,1420-8946},
   MRCLASS = {57K40 (57K41 57K43 57N35 57N65)},
  MRNUMBER = {4888080},
       DOI = {10.4171/cmh/587},
       URL = {https://doi.org/10.4171/cmh/587},
}

@article {conway-powell:pi1-Z,
    AUTHOR = {Conway, Anthony and Powell, Mark},
     TITLE = {Embedded surfaces with infinite cyclic knot group},
   JOURNAL = {Geom. Topol.},
  FJOURNAL = {Geometry \& Topology},
    VOLUME = {27},
      YEAR = {2023},
    NUMBER = {2},
     PAGES = {739--821},
      ISSN = {1465-3060,1364-0380},
   MRCLASS = {57K40 (57N35)},
  MRNUMBER = {4589564},
       DOI = {10.2140/gt.2023.27.739},
       URL = {https://doi.org/10.2140/gt.2023.27.739},
}

@ARTICLE{boyle-chen:equivariant,
    AUTHOR = {Boyle, Keegan and Chen, Wenzhao},
     TITLE = {Equivariant topological slice disks and negative amphichiral
              knots},
   JOURNAL = {Indiana Univ. Math. J.},
  FJOURNAL = {Indiana University Mathematics Journal},
    VOLUME = {73},
      YEAR = {2024},
    NUMBER = {5},
     PAGES = {1623--1637},
      ISSN = {0022-2518,1943-5258},
   MRCLASS = {57K10 (57K14)},
  MRNUMBER = {4845816},
}

@article {orson-powell:isotopy-knot-traces,
    AUTHOR = {Orson, Patrick and Powell, Mark},
     TITLE = {Simple spines of homotopy 2-spheres are unique},
   JOURNAL = {Proc. Lond. Math. Soc. (3)},
  FJOURNAL = {Proceedings of the London Mathematical Society. Third Series},
    VOLUME = {128},
      YEAR = {2024},
    NUMBER = {2},
     PAGES = {Paper No. e12583, 25},
      ISSN = {0024-6115,1460-244X},
   MRCLASS = {57K40 (57N35)},
  MRNUMBER = {4753775},
}

@ARTICLE{conway:discs-samegroup,
    AUTHOR = {Conway, Anthony},
     TITLE = {Homotopy ribbon discs with a fixed group},
   JOURNAL = {Algebr. Geom. Topol.},
  FJOURNAL = {Algebraic \& Geometric Topology},
    VOLUME = {24},
      YEAR = {2024},
    NUMBER = {8},
     PAGES = {4575--4587},
      ISSN = {1472-2747,1472-2739},
   MRCLASS = {57K10 (57N35 57N70 57R67)},
  MRNUMBER = {4843740},
       DOI = {10.2140/agt.2024.24.4575},
       URL = {https://doi.org/10.2140/agt.2024.24.4575},
}

@proceedings {DETbook,
     TITLE = {The disc embedding theorem},
    SERIES = {},
    VOLUME = {},
    EDITOR = {Behrens, Stefan and Kalm\'{a}r, Boldizs\'{a}r and Kim, Min Hoon and Powell, Mark and Ray, Arunima},
 PUBLISHER = {Oxford University Press},
      YEAR = {2021},
     PAGES = {496},
      ISBN = {9780198841319},
   MRCLASS = {},
  MRNUMBER = {},
       DOI = {},
       URL = {},
     Label = {BKKPR},
}

@article{Krushkal-Quinn:2000-1,
	Author = {Krushkal, Vyacheslav S. and Quinn, Frank},
	Fjournal = {Geometry and Topology},
	Issn = {1465-3060},
	Journal = {Geom. Topol.},
	Mrclass = {57N13},
	Mrnumber = {MR1796498 (2001i:57031)},
	Mrreviewer = {Darryl McCullough},
	Pages = {407--430},
	Title = {Subexponential groups in 4-manifold topology},
	Volume = {4},
	Year = {2000}}

@article{Freedman-Teichner:1995-1,
	Author = {Freedman, Michael H. and Teichner, Peter},
	Coden = {INVMBH},
	Fjournal = {Inventiones Mathematicae},
	Issn = {0020-9910},
	Journal = {Invent. Math.},
	Mrclass = {57N13},
	Mrnumber = {MR1359602 (96k:57015)},
	Mrreviewer = {Tatsuhiko Yagasaki},
	Number = {3},
	Pages = {509--529},
	Title = {{$4$}-manifold topology. \textup{{I}}. {S}ubexponential groups},
	Volume = {122},
	Year = {1995}}

@book {FQ,
    AUTHOR = {Freedman, Michael H. and Quinn, Frank},
     TITLE = {Topology of 4-manifolds},
    SERIES = {Princeton Mathematical Series},
    VOLUME = {39},
 PUBLISHER = {Princeton University Press, Princeton, NJ},
      YEAR = {1990},
     PAGES = {viii+259},
      ISBN = {0-691-08577-3},
   MRCLASS = {57N13 (57-02)},
  MRNUMBER = {1201584},
MRREVIEWER = {Ian\ Hambleton},
}

@book {GoGu,
    AUTHOR = {Golubitsky, M. and Guillemin, V.},
     TITLE = {Stable mappings and their singularities},
    SERIES = {Graduate Texts in Mathematics},
    VOLUME = {Vol. 14},
 PUBLISHER = {Springer-Verlag, New York-Heidelberg},
      YEAR = {1973},
     PAGES = {x+209},
   MRCLASS = {58C25},
  MRNUMBER = {341518},
MRREVIEWER = {G.\ R.\ Belitski\u i},
}

@article {marin:transversality,
    AUTHOR = {Marin, A.},
     TITLE = {La transversalit\'{e} topologique},
   JOURNAL = {Ann. of Math. (2)},
  FJOURNAL = {Annals of Mathematics. Second Series},
    VOLUME = {106},
      YEAR = {1977},
    NUMBER = {2},
     PAGES = {269--293},
      ISSN = {0003-486X},
   MRCLASS = {57A40 (57C50)},
  MRNUMBER = {470964},
MRREVIEWER = {C.\ L.\ Seebeck, III},
       DOI = {10.2307/1971096},
       URL = {https://doi.org/10.2307/1971096},
}

@ARTICLE{FNOP:4dguide,
       author = {{Friedl}, Stefan and {Nagel}, Matthias and {Orson}, Patrick and {Powell}, Mark},
        title = "{A survey of the foundations of four-manifold theory in the topological category}",
      journal = {arXiv e-prints (to appear: NYJM Monog.)},
         year = 2019,
        month = oct,
          eid = {arXiv:1910.07372},
        pages = {arXiv:1910.07372},
          doi = {10.48550/arXiv.1910.07372},
archivePrefix = {arXiv},
       eprint = {1910.07372},
 primaryClass = {math.GT},
       adsurl = {https://ui.adsabs.harvard.edu/abs/2019arXiv191007372F}
}

@article {quinn:transversality,
    AUTHOR = {Quinn, Frank},
     TITLE = {Topological transversality holds in all dimensions},
   JOURNAL = {Bull. Amer. Math. Soc. (N.S.)},
  FJOURNAL = {American Mathematical Society. Bulletin. New Series},
    VOLUME = {18},
      YEAR = {1988},
    NUMBER = {2},
     PAGES = {145--148},
      ISSN = {0273-0979,1088-9485},
   MRCLASS = {57N13 (57N15 57N75)},
  MRNUMBER = {929089},
MRREVIEWER = {A.\ Marin},
       DOI = {10.1090/S0273-0979-1988-15629-7},
       URL = {https://doi.org/10.1090/S0273-0979-1988-15629-7},
}

@incollection {DET-book-surgery,
    AUTHOR = {Orson, Patrick and Powell, Mark and Ray, Arunima},
     TITLE = {Surgery theory and the classification of closed, simply
              connected 4-manifolds},
 BOOKTITLE = {The disc embedding theorem},
     PAGES = {331--351},
 PUBLISHER = {Oxford Univ. Press, Oxford},
      YEAR = {2021},
      ISBN = {978-0-19-884131-9},
   MRCLASS = {57K40 (57R65)},
  MRNUMBER = {4519520},
}

@incollection {livingston-slice-survey,
    AUTHOR = {Livingston, Charles},
     TITLE = {A survey of classical knot concordance},
 BOOKTITLE = {Handbook of knot theory},
     PAGES = {319--347},
 PUBLISHER = {Elsevier B. V., Amsterdam},
      YEAR = {2005},
      ISBN = {0-444-51452-X},
   MRCLASS = {57M25},
  MRNUMBER = {2179265},
MRREVIEWER = {Swatee\ Naik},
       DOI = {10.1016/B978-044451452-3/50008-3},
       URL = {https://doi.org/10.1016/B978-044451452-3/50008-3},
}

@article {winterbraids-slice-survey,
    AUTHOR = {Ray, Arunima},
     TITLE = {Slice knots and knot concordance},
   JOURNAL = {Winter Braids Lect. Notes},
  FJOURNAL = {Winter Braids Lecture Notes},
    VOLUME = {8},
      YEAR = {2021},
     PAGES = {Exp. No. 2, 31},
      ISSN = {2426-0312},
   MRCLASS = {57K10 (57K40)},
  MRNUMBER = {4735598},
       DOI = {10.5802/wbln.39},
       URL = {https://doi.org/10.5802/wbln.39},
}

@incollection {gordon-knot-survey,
    AUTHOR = {Gordon, C. McA.},
     TITLE = {Some aspects of classical knot theory},
 BOOKTITLE = {Knot theory ({P}roc. {S}em., {P}lans-sur-{B}ex, 1977)},
    SERIES = {Lecture Notes in Math.},
    VOLUME = {685},
     PAGES = {1--60},
 PUBLISHER = {Springer, Berlin-New York},
      YEAR = {1978},
      ISBN = {3-540-08952-7},
   MRCLASS = {57M25},
  MRNUMBER = {521730},
MRREVIEWER = {Wilbur\ Whitten},
}

@book {rolfsen-book,
    AUTHOR = {Rolfsen, Dale},
     TITLE = {Knots and links},
    SERIES = {Mathematics Lecture Series},
    VOLUME = {7},
      NOTE = {Corrected reprint of the 1976 original},
 PUBLISHER = {Publish or Perish, Inc., Houston, TX},
      YEAR = {1990},
     PAGES = {xiv+439},
      ISBN = {0-914098-16-0},
   MRCLASS = {57M25},
  MRNUMBER = {1277811},
}

@article {garoufalidis-teichner,
    AUTHOR = {Garoufalidis, Stavros and Teichner, Peter},
     TITLE = {On knots with trivial {A}lexander polynomial},
   JOURNAL = {J. Differential Geom.},
  FJOURNAL = {Journal of Differential Geometry},
    VOLUME = {67},
      YEAR = {2004},
    NUMBER = {1},
     PAGES = {167--193},
      ISSN = {0022-040X,1945-743X},
   MRCLASS = {57M25},
  MRNUMBER = {2153483},
MRREVIEWER = {Jacob\ Mostovoy},
       URL = {http://projecteuclid.org/euclid.jdg/1099587731},
}

@inproceedings {freedman-icm,
    AUTHOR = {Freedman, Michael H.},
     TITLE = {The disk theorem for four-dimensional manifolds},
 BOOKTITLE = {Proceedings of the {I}nternational {C}ongress of
              {M}athematicians, {V}ol.\ 1, 2 ({W}arsaw, 1983)},
     PAGES = {647--663},
 PUBLISHER = {PWN, Warsaw},
      YEAR = {1984},
      ISBN = {83-01-05523-5},
   MRCLASS = {57N15 (57R65)},
  MRNUMBER = {804721},
MRREVIEWER = {Frank\ Quinn},
}

@article {Freedman-Alex,
    AUTHOR = {Freedman, Michael H.},
     TITLE = {A surgery sequence in dimension four; the relations with knot
              concordance},
   JOURNAL = {Invent. Math.},
  FJOURNAL = {Inventiones Mathematicae},
    VOLUME = {68},
      YEAR = {1982},
    NUMBER = {2},
     PAGES = {195--226},
      ISSN = {0020-9910},
   MRCLASS = {57M25 (57R65)},
  MRNUMBER = {666159},
MRREVIEWER = {Frank Quinn},
       DOI = {10.1007/BF01394055},
       URL = {https://doi.org/10.1007/BF01394055},
}

@incollection {DET-book-context,
    AUTHOR = {Behrens, Stefan and Powell, Mark and Ray, Arunima},
     TITLE = {Context for the disc embedding theorem},
 BOOKTITLE = {The disc embedding theorem},
     PAGES = {1--26},
 PUBLISHER = {Oxford Univ. Press, Oxford},
      YEAR = {2021},
      ISBN = {978-0-19-884131-9},
   MRCLASS = {57K40 (57K50 57R65 57R67 57R80)},
  MRNUMBER = {4519499},
}

@book {wall-surgery-book,
    AUTHOR = {Wall, C. T. C.},
     TITLE = {Surgery on compact manifolds},
    SERIES = {Mathematical Surveys and Monographs},
    VOLUME = {69},
   EDITION = {Second},
      NOTE = {Edited and with a foreword by A. A. Ranicki},
 PUBLISHER = {American Mathematical Society, Providence, RI},
      YEAR = {1999},
     PAGES = {xvi+302},
      ISBN = {0-8218-0942-3},
   MRCLASS = {57R67 (19J25 57-02)},
  MRNUMBER = {1687388},
       DOI = {10.1090/surv/069},
       URL = {https://doi.org/10.1090/surv/069},
}

@unpublished{Powell-Ray-Teichner:2018-1,
	Author = {Powell, Mark and Ray, Arunima and Teichner, Peter},
	Title = {The $4$-dimensional disc embedding theorem and dual spheres},
 Note = {Preprint, available at ar{X}iv:2006.05209},	
Year = {2020}
}

@book{CLM:surgery-book,
    AUTHOR = {L\"uck, Wolfgang and Macko, Tibor},
     TITLE = {Surgery theory---foundations},
    SERIES = {Grundlehren der mathematischen Wissenschaften (Fundamental
              Principles of Mathematical Sciences)},
    VOLUME = {362},
      NOTE = {With contributions by Diarmuid Crowley},
 PUBLISHER = {Springer, Cham},
      YEAR = {2024},
     PAGES = {xv+956},
      ISBN = {978-3-031-56333-1; 978-3-031-56334-8},
   MRCLASS = {57-02 (57R65)},
  MRNUMBER = {4789553},
       DOI = {10.1007/978-3-031-56334-8},
       URL = {https://doi.org/10.1007/978-3-031-56334-8},
}

@incollection{Sullivan,
	Author = {Sullivan, D. P.},
	Booktitle = {The {H}auptvermutung book},
	Doi = {10.1007/978-94-017-3343-4_3},
	Mrclass = {57Q15 (57Q25 57R10)},
	Mrnumber = {1434103},
	Mrreviewer = {Oliver Attie},
	Pages = {69--103},
	Publisher = {Kluwer Acad. Publ., Dordrecht},
	Series = {$K$-Monogr. Math.},
	Title = {Triangulating and smoothing homotopy equivalences and homeomorphisms. {G}eometric {T}opology {S}eminar {N}otes},
	Url = {https://doi.org/10.1007/978-94-017-3343-4_3},
	Volume = {1},
	Year = {1996}}

@article{Novikov,
	Author = {Novikov, S. P.},
	Fjournal = {Izvestiya Akademii Nauk SSSR. Seriya Matematicheskaya},
	Issn = {0373-2436},
	Journal = {Izv. Akad. Nauk SSSR Ser. Mat.},
	Mrclass = {57.10},
	Mrnumber = {0162246},
	Mrreviewer = {A. H. Wallace},
	Pages = {365--474},
	Title = {Homotopically equivalent smooth manifolds. {I}},
	Volume = {28},
	Year = {1964}}

@book{Browder,
	Author = {Browder, William},
	Mrclass = {57D65},
	Mrnumber = {0358813},
	Mrreviewer = {P. J. Kahn},
	Series = {Ergebnisse der Mathematik und ihrer Grenzgebiete}, 
	Volume={65},
	Pages = {ix+132},
	Publisher = {Springer-Verlag, New York-Heidelberg},
	Title = {Surgery on simply-connected manifolds},
	Year = {1972}}

@book{KS,
	Address = {Princeton, N.J.},
	Author = {Kirby, Robion C. and Siebenmann, Laurence C.},
	Mrclass = {57-02 (57AXX)},
	Mrnumber = {0645390 (58 \#31082)},
	Mrreviewer = {Ronald J. Stern},
	Note = {With notes by John Milnor and Michael Atiyah}, 
	Series={Annals of Mathematics Studies}, 
	Volume={88},
	Pages = {vii+355},
	Publisher = {Princeton University Press},
	Title = {Foundational essays on topological manifolds, smoothings, and triangulations},
	Year = {1977}}

@incollection {kirby-taylor,
    AUTHOR = {Kirby, Robion C. and Taylor, Laurence R.},
     TITLE = {A survey of 4-manifolds through the eyes of surgery},
 BOOKTITLE = {Surveys on surgery theory, {V}ol. 2},
    SERIES = {Ann. of Math. Stud.},
    VOLUME = {149},
     PAGES = {387--421},
 PUBLISHER = {Princeton Univ. Press, Princeton, NJ},
      YEAR = {2001},
      ISBN = {0-691-08814-4; 0-691-08815-2},
   MRCLASS = {57N13 (57R65)},
  MRNUMBER = {1818779},
}

@article {donaldson-1983,
    AUTHOR = {Donaldson, S. K.},
     TITLE = {An application of gauge theory to four-dimensional topology},
   JOURNAL = {J. Differential Geom.},
  FJOURNAL = {Journal of Differential Geometry},
    VOLUME = {18},
      YEAR = {1983},
    NUMBER = {2},
     PAGES = {279--315},
      ISSN = {0022-040X,1945-743X},
   MRCLASS = {57N13 (53C05 57R10 57R55 58E15)},
  MRNUMBER = {710056},
MRREVIEWER = {R.\ C.\ Kirby},
       URL = {http://projecteuclid.org/euclid.jdg/1214437665},
}

@incollection {DET-book-flowchart,
    AUTHOR = {Powell, Mark and Ray, Arunima},
     TITLE = {The development of topological 4-manifold theory},
 BOOKTITLE = {The disc embedding theorem},
     PAGES = {295--330},
 PUBLISHER = {Oxford Univ. Press, Oxford},
      YEAR = {2021},
      ISBN = {978-0-19-884131-9},
   MRCLASS = {57K40},
  MRNUMBER = {4519519},
}

@article {friedl-teichner,
    AUTHOR = {Friedl, Stefan and Teichner, Peter},
     TITLE = {New topologically slice knots},
   JOURNAL = {Geom. Topol.},
  FJOURNAL = {Geometry and Topology},
    VOLUME = {9},
      YEAR = {2005},
     PAGES = {2129--2158},
      ISSN = {1465-3060,1364-0380},
   MRCLASS = {57M25 (57M27 57N70)},
  MRNUMBER = {2209368},
MRREVIEWER = {Jacob\ Andrew\ Rasmussen},
       DOI = {10.2140/gt.2005.9.2129},
       URL = {https://doi.org/10.2140/gt.2005.9.2129},
}

@incollection {venema:1alg,
    AUTHOR = {Venema, Gerard A.},
     TITLE = {Local homotopy properties of topological embeddings in
              codimension two},
 BOOKTITLE = {Geometric topology ({A}thens, {GA}, 1993)},
    SERIES = {AMS/IP Stud. Adv. Math.},
    VOLUME = {2.1},
     PAGES = {388--405},
 PUBLISHER = {Amer. Math. Soc., Providence, RI},
      YEAR = {1997},
      ISBN = {0-8218-0654-8},
   MRCLASS = {57N35 (55Q52)},
  MRNUMBER = {1470738},
       DOI = {10.1090/amsip/002.1/21},
       URL = {https://doi.org/10.1090/amsip/002.1/21},
}

@ARTICLE{conway-dai-miller:cp2,
       author = {{Conway}, Anthony and {Dai}, Irving and {Miller}, Maggie},
        title = "{$\mathbb{Z}$-disks in $\mathbb{C} P^2$}",
      journal = {arXiv e-prints (to appear: Compos.~Math.)},
         year = 2024,
        month = mar,
          eid = {arXiv:2403.10080},
        pages = {arXiv:2403.10080},
          doi = {10.48550/arXiv.2403.10080},
archivePrefix = {arXiv},
       eprint = {2403.10080},
 primaryClass = {math.GT},
       adsurl = {https://ui.adsabs.harvard.edu/abs/2024arXiv240310080C}
}

@article {boyer93,
    AUTHOR = {Boyer, Steven},
     TITLE = {Realization of simply-connected {$4$}-manifolds with a given
              boundary},
   JOURNAL = {Comment. Math. Helv.},
  FJOURNAL = {Commentarii Mathematici Helvetici},
    VOLUME = {68},
      YEAR = {1993},
    NUMBER = {1},
     PAGES = {20--47},
      ISSN = {0010-2571,1420-8946},
   MRCLASS = {57N13 (57N10 57N35)},
  MRNUMBER = {1201200},
MRREVIEWER = {Yoshihisa\ Sato},
       DOI = {10.1007/BF02565808},
       URL = {https://doi.org/10.1007/BF02565808},
}

@article {HK93,
    AUTHOR = {Hambleton, Ian and Kreck, Matthias},
     TITLE = {Cancellation of hyperbolic forms and topological
              four-manifolds},
   JOURNAL = {J. Reine Angew. Math.},
  FJOURNAL = {Journal f\"ur die Reine und Angewandte Mathematik. [Crelle's
              Journal]},
    VOLUME = {443},
      YEAR = {1993},
     PAGES = {21--47},
      ISSN = {0075-4102,1435-5345},
   MRCLASS = {57N13 (57M50)},
  MRNUMBER = {1241127},
MRREVIEWER = {Shmuel\ Weinberger},
}

@article {shaneson-splitting,
    AUTHOR = {Shaneson, Julius L.},
     TITLE = {Wall's surgery obstruction groups for {$G\times Z$}},
   JOURNAL = {Ann. of Math. (2)},
  FJOURNAL = {Annals of Mathematics. Second Series},
    VOLUME = {90},
      YEAR = {1969},
     PAGES = {296--334},
      ISSN = {0003-486X},
   MRCLASS = {57.10},
  MRNUMBER = {246310},
MRREVIEWER = {R.\ Schultz},
       DOI = {10.2307/1970726},
       URL = {https://doi.org/10.2307/1970726},
}

\end{document}